\documentclass[smallextended]{svjour3}   

\usepackage{amsmath,amssymb,amsfonts,color}
\usepackage{algorithmic,algorithm}

\usepackage{graphicx,tikz}
\usepackage{todonotes}
\usepackage{epstopdf}
\usepackage{mathdots} 
\usepackage{graphics}
\usepackage{xargs,hyperref}   
\renewenvironment{example}[1][]{\refstepcounter{example}\par\medskip \noindent \textbf{Test~\theexample. #1} \rmfamily}{\medskip}

\newcommand{\cH}{\mathcal{H}}
\newcommand{\cL}{\mathcal{L}}
\newcommand{\cW}{\mathcal{W}}

\newcommand{\norm}[1]{\left\lVert #1\right\rVert}

\newcommand{\cC}{\mathcal{C}}

\newcommand{\R}{\mathbb{R}}

\newcommand{\cU}{\mathcal{U}}

\newcommand{\cJ}{\mathcal{J}}

\newcommand{\Htwo}{\mathcal{H}_2}
\newcommand{\Hinf}{\mathcal{H}_\infty}

\newcommand{\pt}[1]{\partial_{t} #1}
\newcommand{\pxi}[1]{\partial_{\xi} #1}

\begin{document}
	
\title{State-dependent Riccati equation feedback stabilization for nonlinear PDEs\thanks{Part of this work was supported by the Indam-GNCS 2019 Project ``Tecniche innovative e parallele per sistemi lineari e non lineari di grandi dimensioni, funzioni ed equazioni matriciali ed applicazioni''. VS is a member of the GNCS-Indam activity group. DK was supported by the UK Engineering and Physical Sciences Research Council (EPSRC) grants EP/V04771X/1, EP/T024429/1, and EP/V025899/1. AA was supported by the CNPq research grant  3008414/2019-1 and by a research grant from PUC-Rio.}
}


\author{Alessandro Alla \and
	Dante Kalise \and
	Valeria Simoncini 
}

\authorrunning{A. Alla, D. Kalise, and V. Simoncini} 

\institute{A. Alla \at
	Department of Mathematics, PUC-Rio, Rio de Janeiro, 22451-900 Brazil.\\
	\email{alla@mat.puc-rio.br}           
	\and
	D. Kalise \at
	School of Mathematical Sciences, University of Nottingham, University Park Campus, Nottingham NG7 2RD,  United Kingdom.\\
	\email{dante.kalise@nottingham.ac.uk}
	\and
	V. Simoncini \at
	AM$^2$ and Dipartimento di Matematica, Universit\`a di Bologna,
	Piazza di Porta San Donato  5, I-40127 Bologna, and IMATI-CNR, Pavia, Italy.\\
	\email{valeria.simoncini@unibo.it}	 
}

\date{Received: date / Accepted: date}

\maketitle

\begin{abstract}
The synthesis of suboptimal feedback laws for controlling nonlinear dynamics arising from semi-discretized PDEs is studied. An approach based on the State-dependent Riccati Equation (SDRE) is presented for $\Htwo$ and $\Hinf$ control problems. Depending on the nonlinearity and the dimension of the resulting problem, offline, online, and hybrid offline-online alternatives to the SDRE synthesis are proposed. The hybrid offline-online SDRE method reduces to the sequential solution of Lyapunov equations, effectively enabling the computation of suboptimal feedback controls for two-dimensional PDEs. Numerical tests for the Sine-Gordon, degenerate Zeldovich, and viscous Burgers' PDEs are presented, providing a thorough experimental assessment of the proposed methodology.
\end{abstract}


\section{Introduction} \label{Section1}
Feedback control laws are ubiquitous in modern science and engineering and can be found in autonomous vehicles, fluid flow control, and network dynamics, among many others. A distinctive feature of feedback laws is their ability to compensate external perturbations in real time. While an offline training phase is often affordable, an operational feedback law must be able to provide a control signal at a rate that is determined by the underlying time scales of the physical system to be controlled. This requirement poses a formidable computational constraint for the synthesis of feedback controls which require an online optimization procedure.

A natural approach to generate an optimal feedback law for real-time control is to follow a dynamic programming approach. Here, we solve a nonlinear Hamilton-Jacobi-Bellman (HJB) partial differential equation (PDE) for the value function of the optimal control problem under study. This is done in an offline phase, and once the value function has been computed, a feedback law is obtained as by-product. Once online, the cost of implementing an HJB-based control in real time, assuming the current state of the system is available, is reduced to the evaluation of a nonlinear feedback map. Unfortunately, the dynamic programming approach is not suitable for systems described by large-scale dynamics, as the computational complexity of approximating the associated high-dimensional HJB PDEs goes beyond the reach of traditional computational methods. Only very recently, the use of 
effective computational approaches such as sparse grids~\cite{axelsparse,kang1}, tree structure algorithms~\cite{allatree}, polynomial approximation ~\cite{poli1,KK17,AKK21} tensor decomposition methods~\cite{tensor2,tensor3,tensor4,tensor5}, and representation formulas \cite{chow1,chow3} have addressed the solution of high-dimensional HJB PDEs. Recent works making use of deep learning ~\cite{ml1,ml2,ml3,ml4,ml5,nakazim,nakazim2} anticipate that the synthesis of optimal feedback laws for large-scale dynamics can be a viable path in the near future.

An alternative to the dynamic programming approach is the use of a Nonlinear Model Predictive Control (NMPC) 
framework \cite{larsnmpc}. Here, an optimal open-loop control law is computed at every sampling instant of the dynamics. However, the control action is optimized over a prediction horizon, which is sufficiently large to guarantee asymptotic stability of the closed-loop, but short enough to ensure a computing time that is compatible with the rate at which the control law is called. The NMPC framework embodies the trade-off between optimality and real-time computability. It has been shown \cite{GR2008} that the NMPC corresponds to a relaxation of the dynamic programming approach, in the sense that the NMPC law is suboptimal with respect to the optimal stabilizing feedback law provided by the dynamic programming approach, although its suboptimality can be controlled by increasing the prediction horizon.  

There exists a third synthesis alternative, which incorporates elements from both dynamic programming and NMPC, known as the State-Dependent Riccati Equation (SDRE) approach \cite{sdrefirst1,sdrefirst2}. The SDRE method originates from the dynamic programming and the HJB PDE associated to infinite horizon optimal stabilization, however, it circumvents its solution by reformulating the feedback synthesis as the sequential solution of state-dependent Algebraic Riccati Equations (ARE), which are updated online along a trajectory. The SDRE feedback is implemented similarly as in NMPC, but the online solution of an optimization problem is replaced by a nonlinear matrix equation. Alternatives to the online formulation include the use of neural networks \cite{WANG98,ABK21} in supervised learning, and reformulating the SDRE synthesis as an optimization problem \cite{Astolfi2020}. However, in this work we will focus on addressing the SDRE synthesis from a numerical linear algebra perspective. The efficient solution of matrix equations arising in feedback control has been subject of extensive research over the last decades, leading to the design of solvers which can effectively be applied to large-scale 
dynamics such as those arising in the control of systems governed by PDEs (see, e.g.,
\cite{Antoulas.05,Benner2005a,Benner.Saak.survey13,Binietal.book.12,Simoncini.survey.16,Slowik2007}), and agent-based models \cite{HK18}. Moreover, under certain stabilizability conditions, this feedback law generates a locally asymptotically stable closed-loop and approximates the optimal feedback law from the HJB PDE.

The purpose of the present work is to study the design of SDRE feedback laws for the control of nonlinear PDEs, including nonlinear reaction and transport terms. This is a class of problems that are inherently large-scale, and where the presence of nonlinearities renders linear controllers underperformant. In this framework, the use of dynamic programming or NMPC methods is often prohibitively expensive, as for instance, in feedback control for multi-dimensional PDEs. Here instead, we propose and assess different alternatives for control design based on the SDRE approach which can be effectively implemented for two and three dimensional PDEs. The methodology studied in this paper is based on the parametrization of the SDRE synthesis proposed in \cite{BTB00,BLT07}. In these works, different SDRE feedback laws are developed based on the representation of the nonlinear terms in the system dynamics. This is particularly relevant to PDE dynamics, as unlike the nonlinear ODE world, there exists a clear taxonomy of physically meaningful nonlinearities. In some particular cases, the SDRE approach is reduced to a series of offline computations, and a real-time controller only requires a nonlinear feedback evaluation. We explore the limitations of such a parametrization. In the more general case, the 
SDRE synthesis requires the sequential solution of AREs at a high rate, a task that can is demanding for large-scale dynamics. Here, we propose a variant requiring the solution of a Lyapunov equation instead, thus mitigating the computational cost associated to the online synthesis.

The paper is structured as follows. In section~\ref{sec:review_synth} and its subsections
 we review the $\cH_2/\cH_{\infty}$ optimal feedback synthesis problem. 
In section~\ref{sec:suboptimal} and its subsections
 we present the suboptimal approximation to these feedback laws by means of the SDRE approach, 
including offline, online, and hybrid offline-online implementations, illustrated with 
an application to the control of the Sine-Gordon equation. 
Section~\ref{sec:nla} discusses numerical linear algebra aspects of the solution 
of the Algebraic Riccati and Lyapunov equations which constitute the core building blocks of the SDRE 
feedback synthesis. Finally, in section~\ref{numerics} we perform a computational assessment of the 
proposed methodology on the two-dimensional degenerate Zeldovich and viscous Burgers PDEs.

\section{Optimal feedback synthesis for nonlinear dynamics} \label{sec:review_synth}
In this section we revisit the use of dynamic programming and Hamilton-Jacobi-Bellman/Isaacs PDEs for the computation of optimal feedback controls for nonlinear dynamics. We begin by stating the problem of optimal feedback stabilization via $\cH_2$ synthesis, to then focus on robustness under perturbation in the framework of $\cH_{\infty}$ control.

\subsection{The $\cH_2$ synthesis and the Hamilton-Jacobi-Bellman PDE}
We consider the following infinite horizon optimal control problem:
\[\underset{u(\cdot)\in\cU}{\min}\;\cJ(u(\cdot);x):=\int\limits_0^\infty \Big(\norm{y(t)}_Q^2+\norm{u(t)}_R^2\Big)\, dt
\]
subject to the nonlinear dynamical constraint
\begin{equation}\label{dynamics}
	\dot y(t)= f(y(t))+g(y(t))u(t)\,,\quad y(0)=x,
\end{equation}
where $y(t)=(y_1(t),\ldots,y_d(t))^{\top}\in\R^d$ denotes the state of the system, and the control signal $u(\cdot)$ belongs to $\cU:=L^{\infty}(\R_+;\R^m)$. The running cost is given by $\norm{y}_Q^2:=y^{\top}Qy$ with $Q\in\R^{d\times d},\,Q\succ0$, and $\norm{u}_R^2=u^{\top}Ru$ with $R\in\R^{m\times m},\,R\succ0$. We assume the system dynamics $f(y):\R^d\rightarrow\R^d$ to be such that $f(0)=0$, and the control operator {$g(y):\R^d\rightarrow\R^{d\times m}$} to be $\cC^1(\R^d)$. This formulation of the $\cH_2$ synthesis corresponds to the asymptotic stabilization of nonlinear dynamics towards the origin. By the application of the Dynamic Programming Principle, it is well-known that the optimal value function
\[V(x)=\underset{u(\cdot)\in\cU}{\inf} \cJ(u(\cdot);x)\]
characterizing the solution of this infinite horizon control problem is the unique viscosity solution of the Hamilton-Jacobi-Bellman equation
\begin{equation}\label{hjb}
	\underset{u\in \R^m}{\min}\;\{ \nabla V(x)^{\top}(f(x)+g(x) u)+ \norm{x}_Q^2+\norm{u}_R^2\}=0\,,\quad V(0)=0\,.
\end{equation}
The explicit minimizer $u^*$ of eq. \eqref{hjb} is given by
\begin{equation}\label{optc}
	u^*(x)=-\frac{1}{2}R^{-1} g(x)^{\top} \nabla V(x)\,.
\end{equation}
By inserting \eqref{optc} into \eqref{hjb}, we obtain the HJB equation
\begin{equation}\label{hjb2}
	\nabla V(x)^{\top} f(x)-\frac{1}{4}\nabla V(x)^{\top}g(x)R^{-1}g(x)^{\top}\nabla V(x)+x^{\top}Qx=0\,,\tag{HJB}
\end{equation}
to be understood in the classical sense. We recall that under the additional linearity assumptions $f(x)=Ax$ with $A\in\R^{d\times d}$ and $g(x)=B\in\R^{d\times m}$, the value function is known to be a quadratic form, $V(x)=x^{\top}\Pi x$, with $\Pi\in\R^{d\times d}$ positive definite, and eq. \eqref{hjb2} becomes an Algebraic Riccati Equation for $\Pi$
\begin{align}\label{are}
	A^{\top}\Pi+\Pi A -\Pi B R^{-1}B^{\top}\Pi+Q=0\,.\tag{ARE}
\end{align}
Solving for $V(x)$ in eq. \eqref{hjb2} globally in $\R^d$ allows an online synthesis of the optimal feedback law \eqref{optc} by evaluating the gradient $\nabla V(x)$ and $g(x)$ at the current state $x$. This leads to an inherently robust control law in the sense that if the state of the system is perturbed to $x'=x+\delta x$, there still exists a stabilizing feedback action departing from the perturbed state, namely, $u^*(x')$. However, this control design neglects the modelling of the disturbance/uncertainties affecting the dynamics, with no general stabilization guarantees. We overcome this limitation by formulating an $\cH_{\infty}$ synthesis, which we describe in the following.

\subsection{The $\cH_{\infty}$ problem and the Hamilton-Jacobi-Isaacs PDE} 
We extend the previous formulation by considering a nonlinear dynamical system of the form
\begin{equation}\label{dynamicsp}
	\dot y(t)= f(y(t))+g(y(t))u(t)+h(y(t))w(t)\,,\quad y(0)=x\,,
\end{equation}
where an additional disturbance signal $w(\cdot)\in\cW$, with $\cW=L^2(\R_+;\R^p)$ enters the system through $h(y):\R^d\rightarrow\R^{d\times p}$. We assume that $y=0$ is an equilibrium of the system for $u=w=0$. The $\cH_{\infty}$ control objective is to achieve both internal stability of the closed-loop dynamics and disturbance attenuation through the design of a feedback law $u=u(y)$ such that for a given $\gamma>0$, and for all $T\geq 0$ and $w\in L_2([0,T[,\R^p)$, the inequality
\begin{equation}\label{eq:kk10}
	\int\limits_{0}^{T}\Big(\norm{y}^2_{Q}+\norm{u(t)}^2_R\Big)\, dt\leq\gamma^2\int\limits_{0}^{T}\norm{w(t)}^2_{P}\;dt
\end{equation}
holds. Here, $P\in\R^{p\times p},\, P\succ0\,,$ and $y$ is the solution to \eqref{dynamicsp}. 
The parameter $\gamma$ plays a crucial role in $\cH_{\infty}$ control design. We say that the control system \eqref{dynamicsp} has ${\cL}_2$-gain not greater than $\gamma$,
if \eqref{eq:kk10} holds. Finding the smallest $\gamma^*$ for which this inequality is verified, also known as the $\cH_{\infty}$-norm of the system, is a challenging problem in its own right. The simplest yet computationally expensive method to find the $\cH_{\infty}$-norm of a system is through a bisection algorithm, as described in \cite{BBK89}.
Applying dynamic programming techniques leads to a characterization of the value function for this problem in terms of the solution of a Hamilton-Jacobi-Isaacs PDE
\begin{equation}\label{hji}
	\underset{u\in\R^m}{\min}\,\,\underset{w\in \R^p}{\max}\{ \nabla V(x)^{\top}(f(x)+g(x) u+h(x)w)+ \norm{x}^2_Q+\norm{u}_R^2-\gamma^2\|w\|_P^2\}=0\,,
\end{equation}
valid for all $\gamma \geq \gamma^*$. The unconstrained minimizer and maximizer of \eqref{hji}, $u^*_{\gamma}$ and $w^*_{\gamma}$ respectively, are explicitly computed as
\begin{align}
	u^*_\gamma(x)&=-\frac12 R^{-1}g(x)^{\top}\nabla V_\gamma(x)\,,\label{uopt}\\
	w^*_\gamma(x)&=\frac{1}{2\gamma^2} P^{-1}h(x)^{\top}\nabla V_\gamma(x)\,,\label{wopt}
\end{align}
where $ V_\gamma(x)$ solves the Hamilton-Jacobi-Isaacs equation
\begin{align}\label{hjie}
	\nabla V_\gamma(x)^{\top}f(x)+ \frac14 \nabla V_\gamma(x)^{\top}S(x)\nabla V_\gamma(x)+x^{\top}Qx=0 \,,\tag{HJI}
\end{align}
with
\begin{align}
	S(x)=\frac{1}{\gamma^2}h(x)P^{-1}h(x)^{\top}-g(x)R^{-1}g(x)^{\top},
\end{align}
and $V_\gamma(0)=0$. For an initial condition which is not a steady state we have
\begin{equation}
	\int\limits_{0}^{\top}\Big(\norm{y}+\norm{u(t)}^2_R\Big)\, dt\leq\gamma^2\int\limits_{0}^{\top}\norm{w(t)}^2_{P}\;dt + V_\gamma(x),
\end{equation}
see e.g. \cite[Theorem 16]{V92}.
When there is no confusion, we denote $V_\gamma(x)$ by $V$.
Note that if the disturbance attenuation is neglected by taking $\gamma\to\infty$, we recover the solution of \eqref{hjb2} instead.
Moreover, under the linearity assumptions $f(x)=Ax$ with $A\in\R^{d\times d}$, $g(x)=B\in \R^{d\times m}$, and $h(x)=H\in \R^{d\times p}$, the value function is a quadratic form $V_{\gamma}(x)=x^{\top}\Pi x$, where $\Pi\in\R^{d\times d}$ is positive definite, and eq. \eqref{hjie} becomes the following Algebraic Riccati Equation for $\Pi$
\begin{equation}\label{arehinf}
	A^\top\Pi+\Pi A-\Pi\left( BR^{-1}B^\top- \dfrac{1}{2\gamma^2}HP^{-1}H^\top\right)\Pi+Q=0\,,\tag{\mbox{$\text{ARE}_{\infty}$}}
\end{equation}
solving the $\cH_{\infty}$ problem for full state feedback. In the following section we discuss the construction of computational methods to synthesize nonlinear feedback laws inspired by the solution of HJB/HJBI PDEs.

\section{Sub-optimal feedback control laws}\label{sec:suboptimal}
\renewcommand{\thealgorithm}{3.\arabic{algorithm}}
\setcounter{algorithm}{0}
Despite the extensive parametrization of the control problem, eqns. \eqref{hjb2} and \eqref{hjie} remain first-order, stationary nonlinear PDEs whose numerical approximation is a challenging task. These nonlinear PDEs are cast over $\R^d$, where $d$ relates to the dimension of the state space of the dynamics, which can be arbitrarily large. In particular, if the system dynamics correspond to nonlinear PDEs, the direct solution of the resulting 
infinite-dimensional eqns. \eqref{hjb2} and \eqref{hjie} remains an open computational problem. 
We explore different alternatives to circumvent the solution of high-dimensional HJB PDEs by resorting to the sequential solution of the Algebraic Riccati Equations \eqref{are} and \eqref{arehinf} respectively, providing a tractable alternative for feedback synthesis in large-scale nonlinear dynamics. The different techniques we propose trade the optimality associated to the HJB-based control for computability, and hence will be referred to as \textsl{suboptimal} feedback laws.

For the sake of simplicity, we focus on the $\Hinf$ synthesis. The $\Htwo$ feedback follows directly choosing $h(x(t)):=0$ in \eqref{dynamicsp}. We assume a quadratic cost, i.e. $\ell(x) = x^\top Qx$ with $Q\in\R^{d\times d}$ symmetric positive semidefinite.

\subsection{Linear-Quadratic Regulator (LQR) control law}
The simplest suboptimal control law for nonlinear systems uses the optimal feedback law for the linear-quadratic control problem arising from linearization around an equilibrium, which we assume to be $x=0$. For $\gamma>\gamma^*$, we solve eq. \eqref{arehinf} with $A_{ij}=\frac{\partial f_i(x)}{\partial x_j}|_{x=0}$, $B=g(0)$, and $H=h(0)$. From the solution $\Pi$, we obtain the linear feedback law
\begin{equation}\label{lqr}
u(x)=-R^{-1}B^\top\Pi x\,.
\end{equation}
Such a feedback law, applied to the original nonlinear system, can only guarantee stabilization in a neighbourhood around the origin. 

\subsection{State-Dependent Riccati Equation (SDRE)}\label{secsdre}
If we express the system dynamics through a space-dependent representation of the form
\begin{equation}\label{sdress1}
\dot x=A(x)x+B(x)u(t) +H(x)\omega(t)\,,
\end{equation} 
we can synthesize a suboptimal feedback law inherited from \eqref{hjie} by following an approach known as the State-dependent Riccati Equation (SDRE). This method has been thoroughly analyzed in the literature, 
see, e.g.\cite{C97,BLT07}, and despite being purely formal, it is extremely effective for stabilizing nonlinear dynamics. The SDRE approach is based on the idea that infinite horizon optimal feedback control for systems of the form \eqref{sdress1} are linked to a state-dependent ARE:
\begin{equation}\label{sdre1}
A^\top(x)\Pi(x)+\Pi(x) A(x)-\Pi(x)S(x)\Pi(x)+Q=0\,,
\end{equation}
where
\[
S(x)=B(x)R^{-1}B^\top(x)- \dfrac{1}{2\gamma^2}H(x)P^{-1}H(x)^\top\,.
\]
Solving this equation leads to a state-dependent Riccati operator $\Pi(x)$, from where we obtain a nonlinear feedback law given by
\begin{equation}\label{sdref1}
u(x):=-R^{-1}B^\top(x)\Pi(x)x\,.
\end{equation}

It is important to observe that the operator equation \eqref{sdre1} admits an analytical solution only in a limited number of cases. More importantly, even if this solution is computed for every state $x$, the closed-loop differs from the optimal feedback obtained from solving \eqref{hjie}, as the SDRE approach assumes that the value function is always locally approximated as $V(x)=x^{\top}\Pi(x)x$. From an optimal control perspective the SDRE can be interpreted as a model predictive control loop where at a given instant, the dynamics $(A(x),B(x),H(x))$ are frozen at the current state and an LQR feedback is numerically approximated. The procedure is illustrated in Algorithm \ref{Alg_sdre}. The resulting feedback is naturally different from the optimal nonlinear feedback law, and will remain different regardless of how fast the SDRE feedback is updated along a trajectory. The latter is also observed in \cite{Astolfi2020}, where it is shown that a direct derivation of the SDRE approach departing from the HJB PDE would lead to additional terms in \eqref{sdress1}. Notwithstanding, it is possible to show local asymptotic stability for the SDRE feedback. We recall the following result \cite{BLT07} concerning asymptotic stability of the SDRE closed-loop in the $\Htwo$ case.
\begin{proposition}Assume a nonlinear system
	\begin{equation}
		\dot x(t)=f(x(t))+B(x(t))u(t)\,,
	\end{equation}
	where $f(\cdot)$ is $\cC^1$ for $\|x\|\leq\delta$, and $B(\cdot)$ is continuous. If $f$ is parametrized in the form $f(x)=A(x)x$, and the pair $(A(x),B(x))$ is stabilizable for every $x$ in a non-empty neighbourhood of the origin. Then, the closed-loop dynamics generated by the feedback law \eqref{sdref1} are locally asymptotically stable.
\end{proposition}

\begin{algorithm}[ht]
	\caption{SDRE-MPC loop}
	\label{Alg_sdre}
	\begin{algorithmic}[1]
		\REQUIRE $\{t_0,t_1,\ldots\},$ model \eqref{sdress1}, $R,Q$,
		\FOR{$k=0,1,\ldots$}
		\STATE Compute $\Pi(x(t_k))$ from \eqref{sdre1}
		\STATE Set $K(x(t_k)) := R^{-1}B^\top(x(t_k))\Pi(x(t_k))$
		\STATE Set $u(t):=-K(x(t_k))x(t)$
		\STATE Integrate system dynamics to $x(t_{k+1})$
		\ENDFOR
	\end{algorithmic}
\end{algorithm}

Assuming the stabilizability hypothesis above, the main bottleneck in the implementation of Algorithm \ref{Alg_sdre} is the high rate of calls to an ARE solver for \eqref{sdre1}. Moreover, these ARE calls are expected to be sufficiently fast for real-time feedback control. This is a demanding computational task for the type of large-scale dynamics arising in optimal control of PDEs. In the following, we discuss two variations the SDRE-MPC algorithm which mitigate this limitation by resorting to offline and more efficient online computations.

\subsection{Offline approximation of the SDRE}
A first alternative for more efficient SDRE computations was proposed in \cite{BTB00}, inspired by a power series argument first discussed in \cite{Wernli1975}. Assuming that the state operator $A(x)$ can be decomposed into $A(x)=A_0+f_1(x)A_1\,,$ where $A_0$ and $A_1$ are constant matrices and $f(x)$ is a scalar function, $B(x)=B$ and $H(x)=H$, then the Riccati operator $\Pi(x)$ solving \eqref{sdre1} is approximated by
\[
\Pi(x)=\sum_{n=0}^{\infty}(f_1(x))^nL_n\,,
\]
where the matrices $L_n\in\R^{d\times d}$ solve 
\begin{align}
		L_{0} A_{0}+A_{0}^{\top} L_{0}-L_{0} S L_{0}+Q&=0\,,\label{are0}\\
		L_{1}\left(A_{0}-S L_{0}\right)+\left(A_{0}^{\top}-L_{0} S\right) L_{1}+L_{0} A_1+A_1^{\top} L_{0}&=0\,,\\
		L_{n}\left(A_{0}-S L_{0}\right)+\left(A_{0}^{\top}-L_{0} S\right)L_{n}+Q_n&=0\,,\\
		Q_n:=L_{n-1}A_1+A_1^{\top}L_{n-1}-\sum_{k=1}^{n-1}L_{k} SL_{n-k}&\notag\,.
\end{align}
After solving a single ARE and $N$ Lyapunov equations, an $N-$order approximation of $\Pi(x)$ yields the feedback
\begin{equation}
u^N(x)=-R^{-1}B^\top\left(\sum\limits_{n=0}^N(f_1(x))^nL_n\right)x
\end{equation}
Unfortunately, the reduction above is only possible for a scalar nonlinearity. If the nonlinearity that can be expressed as
\begin{equation}\label{sdrea}
A(x)=A_0+\sum_{j=1}^r f_j(x)A_j\,,
\end{equation}
where $A_j \in R^{d\times d}$ and the state dependence is restricted to $r$ scalar functions $f_j(x):R^d\rightarrow\R$, then a first order approximation of $\Pi(x)$ is given by
\begin{equation}\label{pbar}
\Pi(x)\approx \tilde\Pi(x)=\Pi_0+\sum_{j=1}^r\Pi_jf_j(x)
\end{equation}
where $\Pi_0$ solves \eqref{arehinf} for $A_0$ and the remaining $\Pi_j$ satisfy the Lyapunov equations
\begin{equation}\label{lyap_off}
	\Pi_jC_0+C_0^\top\Pi_j+Q_j=0\,,\quad j=1,\ldots, r\,,
\end{equation}	
with $C_0=A_0-S\Pi_0$, and $Q_j=\Pi_0A_j+A_j\Pi_0$. The resulting feedback law is given by
\begin{equation}\label{cont_off}
	u(x)=-R^{-1}B^\top\tilde\Pi(x)x\,.
\end{equation} 
 Overall, this approach requires the computation of the ARE associated to $A_0$ in addition to $r$ Lyapunov equations whose solution is fully parallelizable. Its implementation is summarized below.
\begin{algorithm}[H]
\caption{Offline SDRE}
\label{Alg_sdreoff}
\begin{algorithmic}[1]
\STATE Offline computations:
\STATE Compute $\Pi_0$ from \eqref{arehinf}
\FOR{$j=1,\ldots,r$}
\STATE Compute $\Pi_j$ from the Lyapunov equation \eqref{lyap_off}
\ENDFOR
\STATE Online process:
\STATE Set $u(x(t)):= -R^{-1}B^\top\tilde\Pi(x(t))x(t)$
\STATE Integrate system dynamics for $x(t)$
\end{algorithmic}
\end{algorithm}

\paragraph{\bf An offline-online SDRE approach} Although the previous approach is a valid variant to circumvent the online solution of Riccati equations at a high rate, it becomes unfeasible in cases where 
both $d$ and $r$ are large, as it requires the solution of $r$ Lyapunov equations \eqref{lyap_off}  and storage of the solution matrix with $d^2$ entries each. Such a large-scale setting arises naturally in feedback control of dynamics from semidiscretization of nonlinear PDEs and agent-based models. We present a variant of the offline SDRE approach which 
circumvents this limitation by resorting to an online phase requiring the solution of a single Lyapunov equation per step. 

Let us define the quantity $W(x) = \sum\limits_{j=1}^r\Pi_jf_j(x)$. Multiplying each equation in \eqref{lyap_off} 
by its corresponding $f_j(x)$, it follows that $W(x)$ satisfies the Lyapunov equation
\begin{equation}\label{lyap_on}
W(x) C_0+C_0^\top W(x) +\sum_{j=1}^r Q_jf_j(x) = 0.
\end{equation}
Therefore, the feedback law can be expressed as
\begin{equation}\label{con:var}
u(x)=-R^{-1}B^\top\left(\Pi_0+\sum_{j=1}^r\Pi_jf_j(x)\right)x =
-R^{-1}B^\top\left(\Pi_0+W(x) \right)x .
\end{equation} 
The feedback law \eqref{con:var} can be computed by solving an offline Riccati equation for $\Pi_0$ 
and an online Lyapunov equation for $W(x)$; {see 
section~\ref{sec:nla} for a discussion on the computational costs of the two approaches.} The offline-online SDRE approach is summarized in Algorithm \ref{Alg_sdreon} below.
\begin{algorithm}[H]
\caption{Offline-online SDRE}
\label{Alg_sdreon}
\begin{algorithmic}[1]
\STATE Offline phase:	
\STATE Compute $\Pi_0$ from \eqref{arehinf}
\STATE Online phase:
\FOR{$k=0,1,\ldots$}
\STATE Compute $W(x(t_k))$ from the Lyapunov equation \eqref{lyap_on}
\STATE Set $K(x(t_k)) := R^{-1}B^\top\left(\Pi_0 +W(x(t_k))\right)$
\STATE Set $u(t):=-K(x(t_k))x(t)$
\STATE Integrate system dynamics to $x(t_{k+1})$
\ENDFOR
\end{algorithmic}
\end{algorithm}

\begin{remark}
The approximation of the state space solution $x(t_{k+1})$ at time $t_{k+1}$ can be performed with both explicit or implicit time-stepping schemes. In both cases, the feedback gain remains frozen at $K(x(t_n))$. 
\end{remark}

\subsection{A preliminary test: the damped Sine-Gordon equation}\label{sinegordon}
We present a preliminary assessment of the effectiveness of Algorithm \ref{Alg_sdre} and Algorithm \ref{Alg_sdreon} for the $\mathcal{H}_2$ case. In section \ref{numerics} we will focus on large-scale, two-dimensional PDEs and $\mathcal{H}_\infty$ control. Given a domain $\Omega\subset\R$, we consider {the control of} the damped Sine-Gordon equation {(see e.g. \cite{MR04})} with homogeneous Dirichlet boundary conditions over $\Omega\times\R^+_0$:
\begin{align}\label{ex:wave}
\begin{aligned}
\partial_{tt} X(\xi,t)&= -\alpha \pt X(\xi,t) + \partial_{\xi\xi} X(\xi,t)-\beta \sin X(\xi,t)+\chi_{\omega_c}(\xi)u(t)\\
X(\xi,t)&=0 \quad\xi\in\partial\Omega, t>0\,,\\
 X(\xi,0)&={\mathbf x}_0(\xi)\,,\;\xi\in\Omega\,,\\
 \pt X (\xi,0)&={\mathbf x}_1(\xi)\,,\;\xi\in\Omega\,,
\end{aligned}
\end{align}
where the control variable $u(t)$ acts through an indicator function $\chi_{\omega_c}(\xi)$ supported over  $\omega_c\subset\Omega$.
{The cost functional to be minimized is given by:
\begin{equation}\label{costex_wave}
\cJ(u(\cdot);X(\cdot,0)):=\int\limits_0^{\infty} \sum_{i=1}^z \dfrac{1}{|\omega_{o_i}|}\left(\int_{\omega_{o_i}} X(\xi,t)\, d\xi\right)^2+R\,|u(t)|^2 \, dt\,
\end{equation}
where $\omega_o:=\cup_{i = 1}^z \omega_{o_i} \subset\Omega$ represents a collection of local patches where we average the state.
} 
Defining $y(t)=(X(\cdot,t), \dot X(\cdot,t))^{\top}$, we write the dynamics as a first-order abstract evolution system
$$
\dot{y}(t)=Ay(t) + f(y(t))+ Bu(t)\;,
$$
where
\begin{align}\label{ops}
    A=
          \begin{bmatrix}
           0 & I \\
            \partial^2_{\xi\xi} & -\alpha I
          \end{bmatrix} ,\,\,\,
          f(y(t))=
		  \begin{bmatrix}
         0 \\     
         -\beta \sin(X(\cdot,t))
          \end{bmatrix}, \,\,\,     
          Bu(t)=
		  \begin{bmatrix}
         0\\
         \chi_{\omega_c}(\xi) u(t)
          \end{bmatrix}.         
\end{align}

\begin{figure}[htbp]
\centering
\includegraphics[scale=0.3]{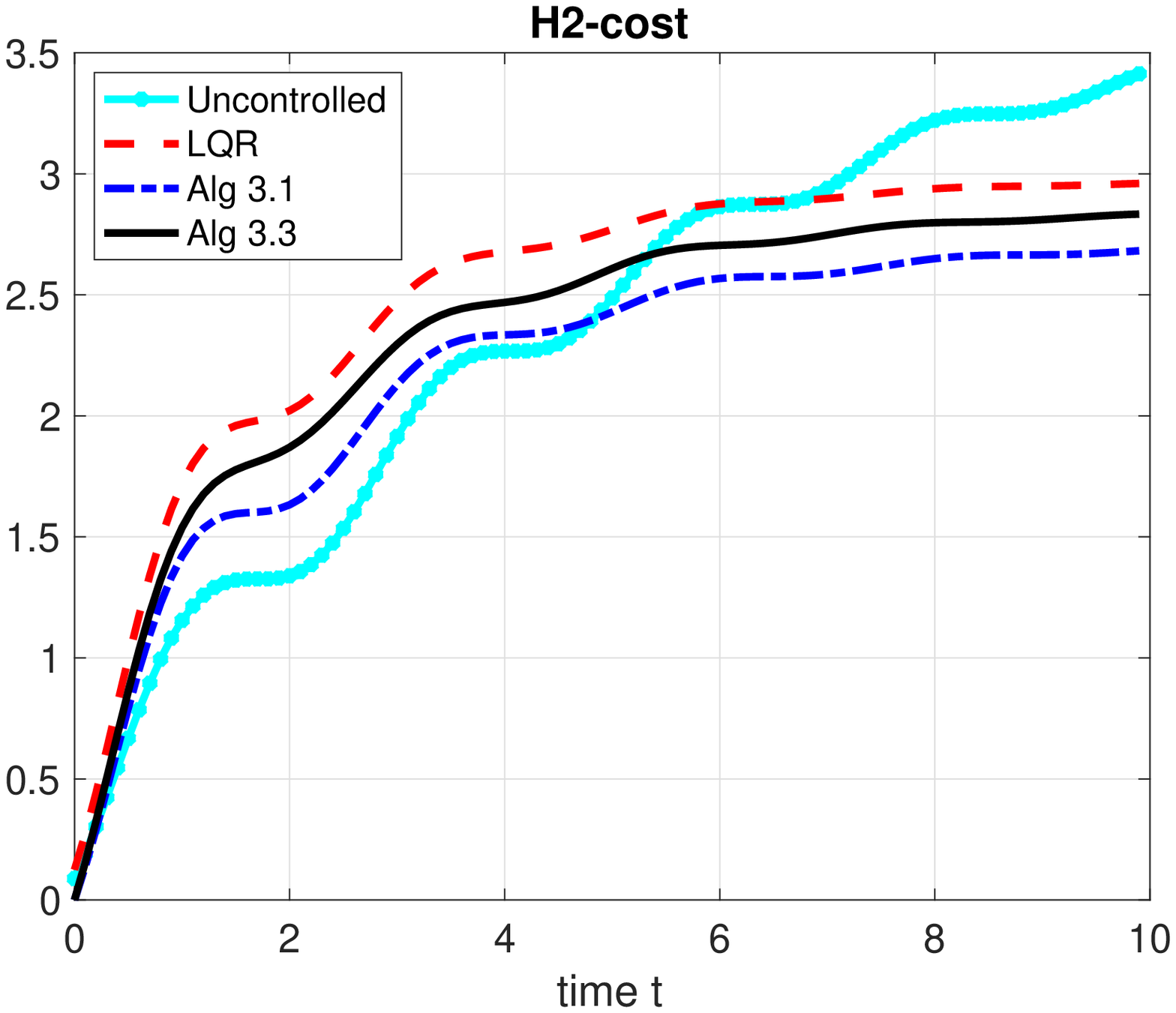}
\includegraphics[scale=0.3]{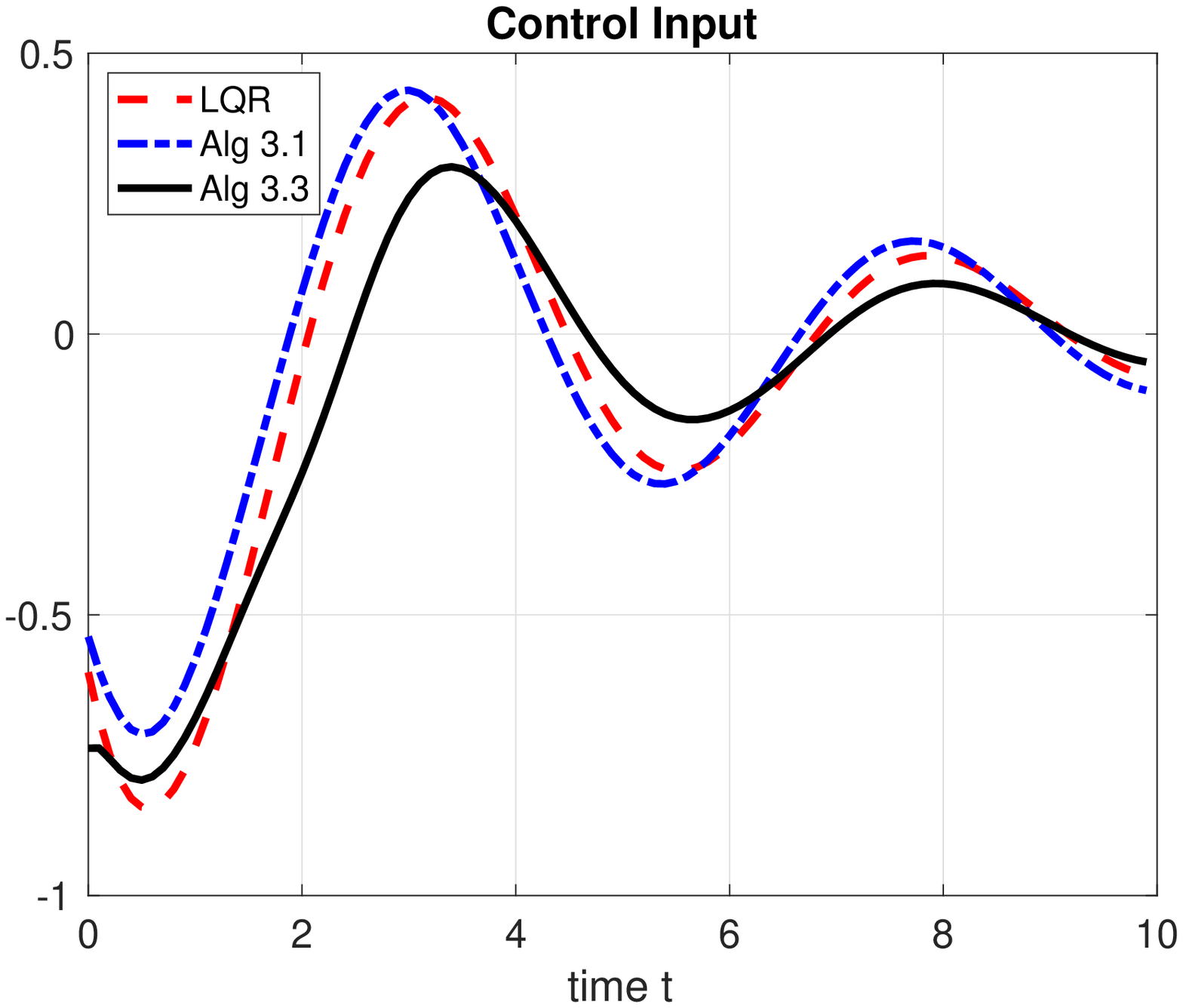}
\includegraphics[scale=0.29]{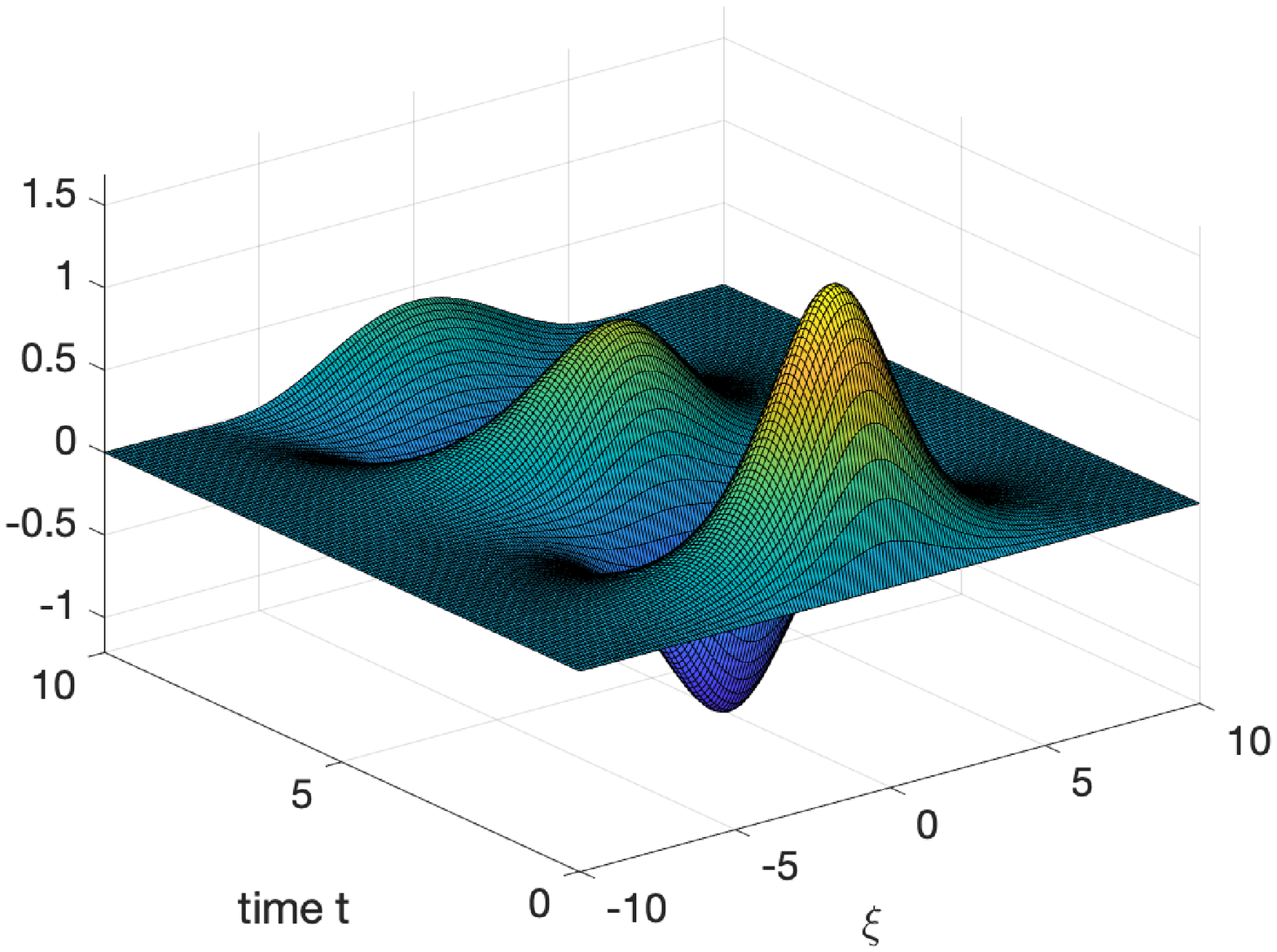}
\includegraphics[scale=0.29]{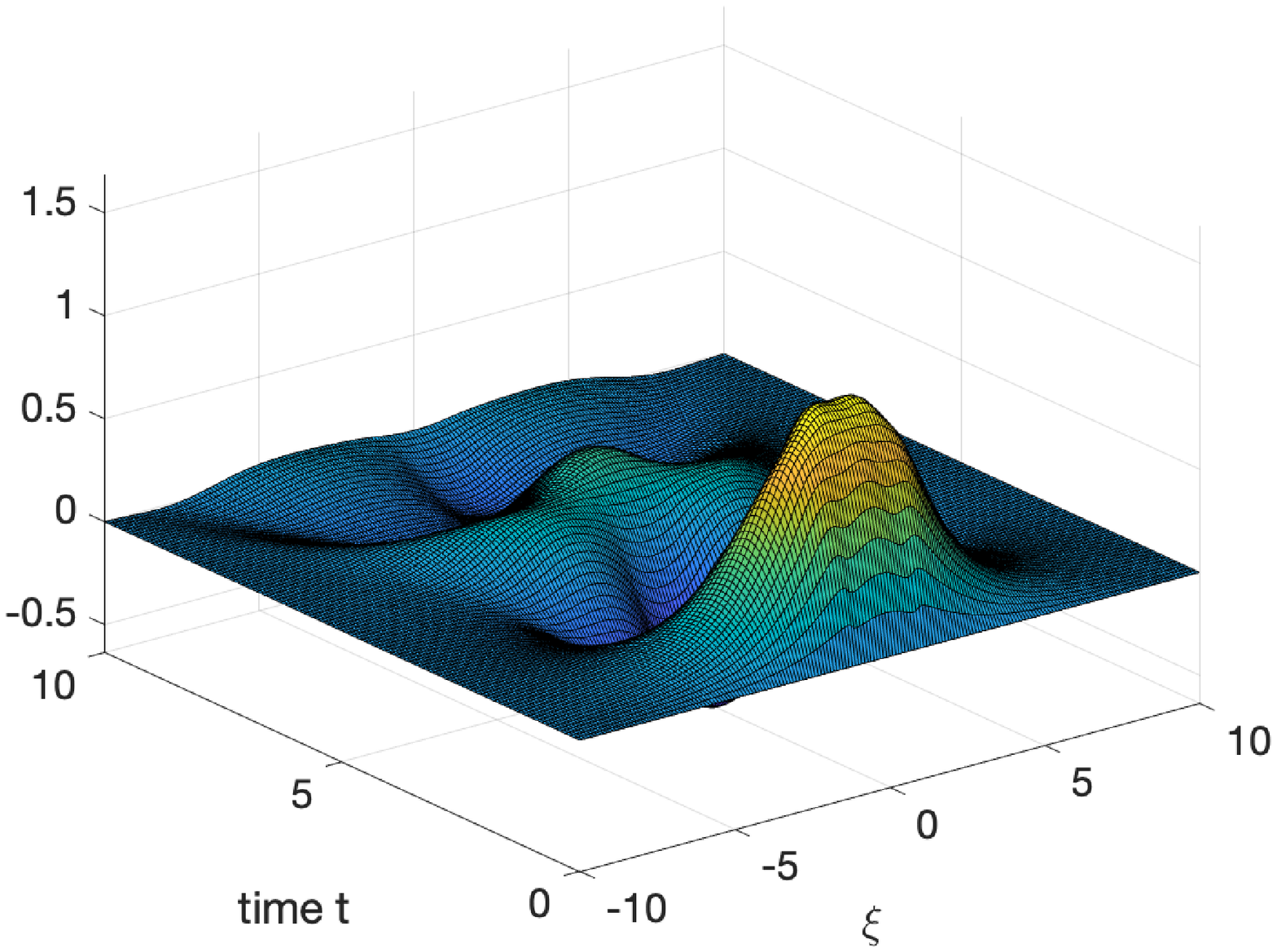}
\includegraphics[scale=0.29]{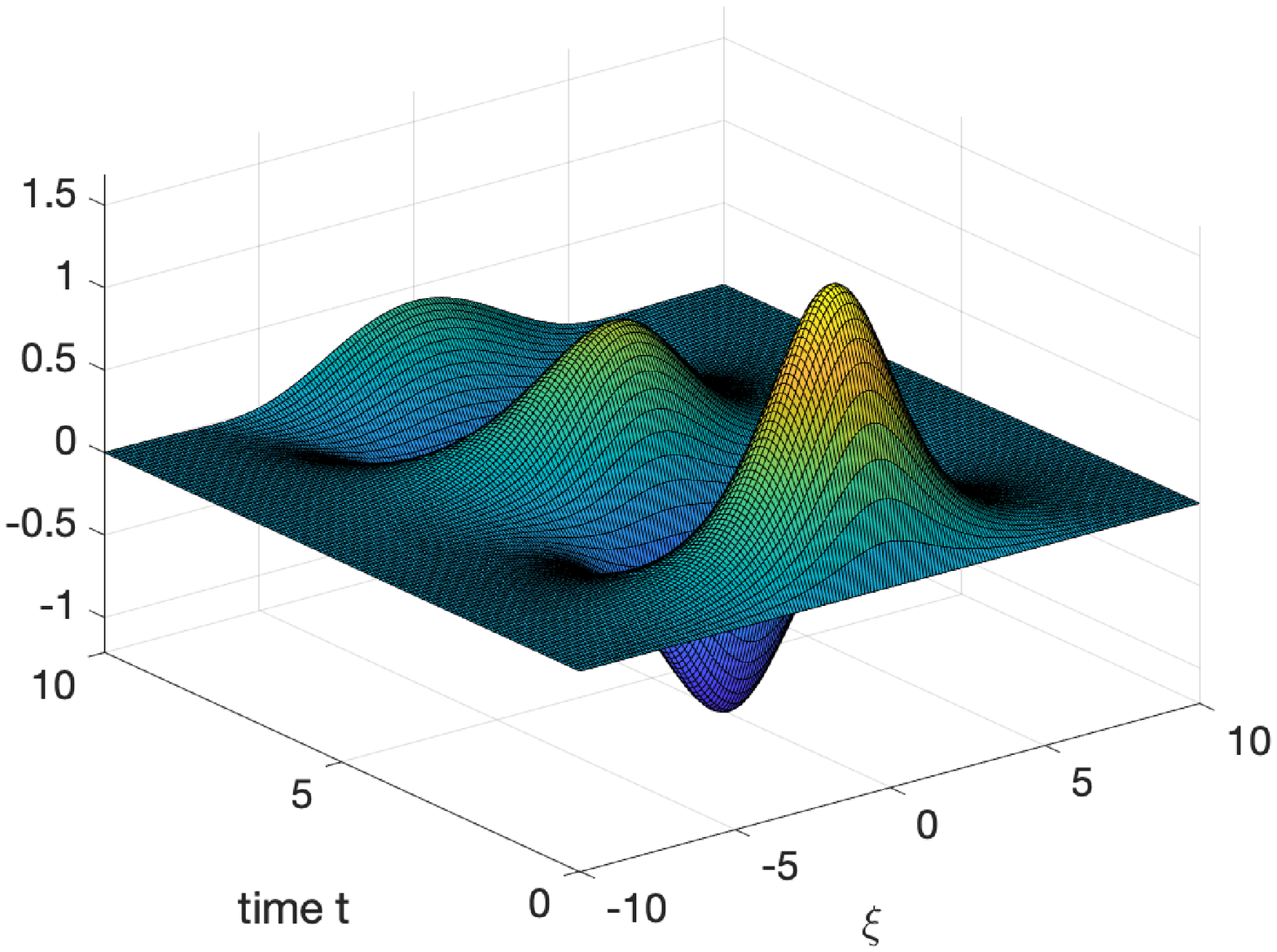}
\includegraphics[scale=0.29]{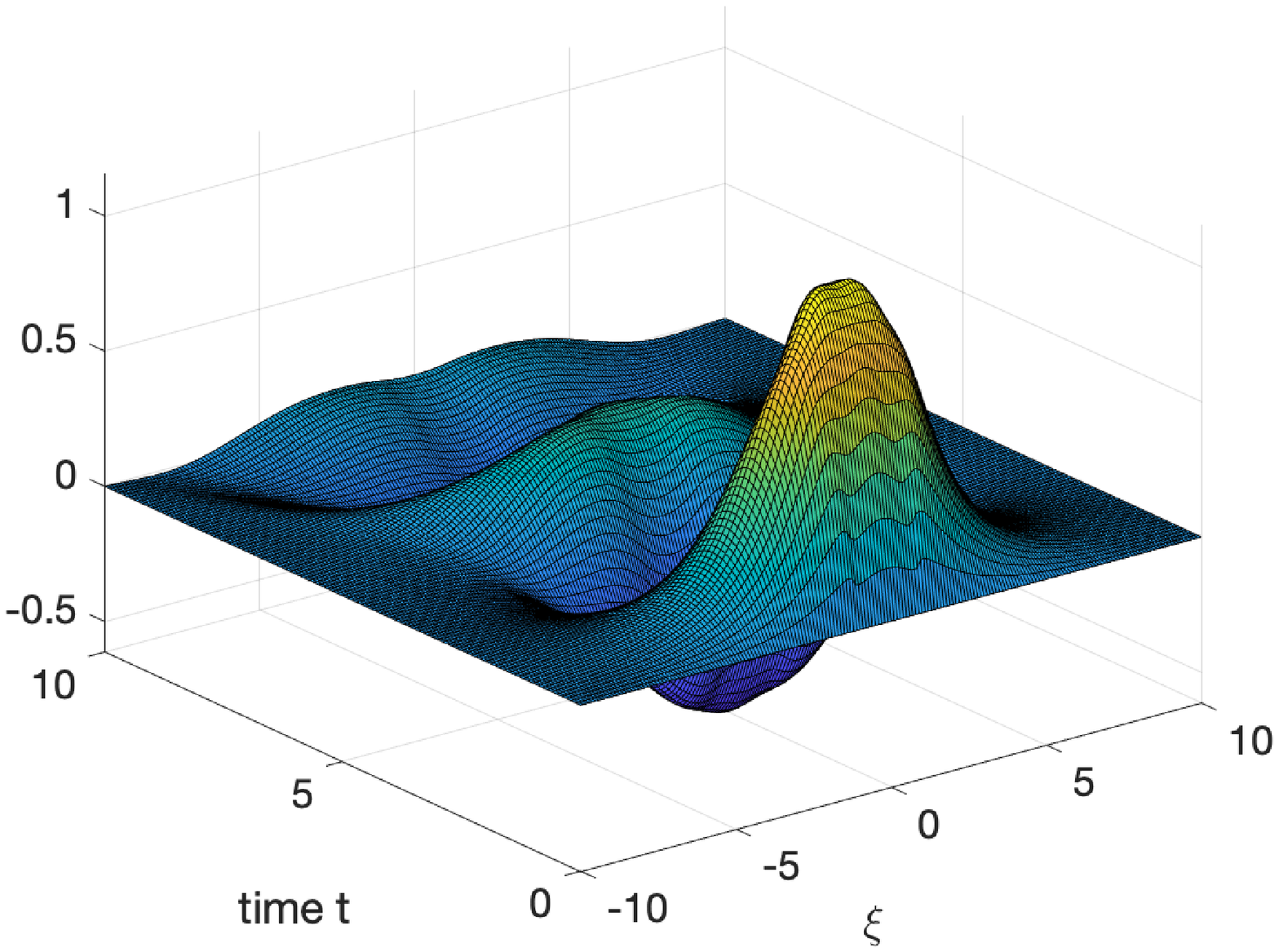}
\caption{Section \ref{sinegordon}. Damped sine-Gordon equation, controlled trajectories. Top:
accumulated running cost with different $\Htwo-$control (left) and corresponding control inputs (right). 
Middle: trajectories, uncontrolled (left), and linearized (right). 
Bottom: trajectories, $\mathcal{H}_2$-controlled solution with  SDRE-MPC Algorithm \ref{Alg_sdre} (left), 
and $\mathcal{H}_2$-controlled solution with SDRE offline-online Algorithm \ref{Alg_sdreon} (right). 
}
\label{fig1:wave}
\end{figure}

We approximate the operators above using a finite difference discretization in space. Given the domain $\Omega=[\xi_L,\xi_R]$, we construct the uniform grid $\xi_i=\xi_L+(i-1)\Delta \xi$ with $\Delta\xi=(\xi_R-\xi_L)/(d-1)$. We define the discrete state $X_i(t):=X(\xi_i,t)$, and the discrete augmented state $Y(t)=(X_1(t),\ldots,X_d(t),\dot X_1(t),\ldots, \dot X_d(t))^{\top}$. The discrete operators read
\begin{align}\label{opsd}
	A_d=
	\begin{bmatrix}
		\mathbf{0}_{d\times d} \;\; \mathbf{I}_{d} \\
		-\Delta_d \;\; \alpha \mathbf{I}_{d}
	\end{bmatrix} ,\quad
	Bu(t)=
	\begin{bmatrix}
		\mathbf{0}_{d\times 1}\\
		\{\chi_{\omega_c}(\xi_i)\}_{i=1}^d
	\end{bmatrix}u(t)\,,         
\end{align}
where $\mathbf{I}_d$ is the $d\times d$ identity matrix, and $\Delta_d$ is the discrete Laplace 
operator\footnote{The notation $\texttt{tridiag}([\begin{array}{ccc}a&b&c\end{array}],d)$ stands for 
a tridiagonal $d\times d$ matrix having the constant values $b\in\mathbb{R}$ on the main diagonal, 
$a\in\mathbb{R}$ on the lower diagonal and $c\in\mathbb{R}$ on the upper diagonal.},  $\Delta_d:=\Delta\xi^{-2}\texttt{tridiag}([\begin{array}{ccc}1&-2&1\end{array}],d)\in\mathbb{R}^{d\times d}$. {The quantity $\sum_{i=1}^z \dfrac{1}{|\omega_{o_i}|}\left(\int_{\omega_{o_i}} X(\xi,t)\, d\xi\right)^2$ in \eqref{costex_wave} is approximated by $X(t)^{\top}Q X(t)$ where $Q= C^{\top}C\in\R^{d\times d}$ and 
$$
C^{\top}=\Delta\xi\, 
\begin{bmatrix} 
\dfrac{\{\chi_{\omega_{o_1}}(\xi_i)\}_{i=1}^d}{|\omega_{o_1}|}, \ldots, 
\dfrac{\{\chi_{\omega_{o_z}}(\xi_i)\}_{i=1}^d}{|\omega_{o_z}|}
\end{bmatrix}\in\R^{d\times z}\,.
$$} 
To express the nonlinearity in a semilinear form consistent with \eqref{sdress1} we define 
\begin{equation}
  \tilde f(Y(t))=
	-\beta	  \begin{pmatrix}
         \mathbf{0}_{d\times 1} \\     
          \left\{\dfrac{\sin(X_i(t))}{X_i(t)}\right\}_{i=1}^d
          \end{pmatrix}Y(t)\,.    
\end{equation}
In our test we consider the following values for the given parameters,
$\Omega = [-10,10],\, t\in[0,10],\, \alpha =0.05,\, \beta = 2$,
$\omega_c = [-1,1]$, and
$$X(\xi,0) = 0,\quad \partial_t X(\xi,0) =\dfrac{8}{3}\mbox{sech}\left(\dfrac{2}{3\xi}\right).$$
In the cost functional {\eqref{costex_wave}} we set $R=1$, $\omega_o(x) = [-2.5,-1.5] \cup [1.5, 2.5]$ and $z=2$. In this test we take $d=402$ nodes in the finite difference discretization. Time-stepping is performed with an implicit Euler method with a time step of $0.1$. 
The small size Riccati and Lyapunov equations are solved using the Matlab functions {\tt icare} 
and {\tt lyap}, respectively. Controlled dynamics with different feedback controls are shown in 
Figure~\ref{fig1:wave}.  The presence of the damping term $\alpha \pt X(\xi,t)$ generates a stable trajectory for both the uncontrolled and LQR-controlled dynamics using the linearized feedback \eqref{lqr}. However, we still observe differences in the state and control variables with respect to the SDRE controllers, SDRE-MPC and SDRE offline-online, described in Algorithms \ref{Alg_sdre} and \ref{Alg_sdreon}, respectively.  The accumulated running costs 
in Figure \ref{fig1:wave} (top-left) indicate that the SDRE-MPC implementation 
achieves the best closed-loop performance, followed by SDRE offline-online, both outperforming linearized LQR and uncontrolled trajectories. However, the main difference between the SDRE-MPC and SDRE offline-online closed-loops is related to 
computational time. 
The SDRE-MPC solver requires the solution of multiple AREs in sequence, taking a total of 
23 minutes of CPU time for this test, whereas the online solution of Lyapunov equations of the 
SDRE offline-online implementation reduces this computation time to 45 seconds. 
A deeper study on the methods performance in the large scale setting will be reported
on in section~\ref{numerics}, while implementation aspects associated with the solution
of these algebraic equations are discussed in the next section.

\section{Solving large-scale Algebraic Riccati/Lyapunov equations} \label{sec:nla}
The time steps discussed in Section~\ref{secsdre} all use 
matrices that stem from the solution of algebraic matrix equations,
and more precisely the (quadratic) Riccati and the (linear) Lyapunov equations. 
The past few years have seen
a dramatic improvement in the effectiveness of numerical solution strategies
for solving these equations in the large scale setting. For a survey
in the linear case we refer the reader to the recent article \cite{Simoncini.survey.16}, while
for the algebraic Riccati equation we point 
the reader to, e.g., \cite{B08,SSM14,Binietal.book.12}.

In our derivation we found 
projection methods to be able to adapt particularly well
to the considered setting, with a similar reduction
framework for both linear and quadratic problems; other approaches are
reviewed for instance in \cite{Benner.Saak.survey13,Simoncini.survey.16}. We emphasize that
all considered methods require that the zero order term in the matrix equation, e.g., matrix $Q$ in
the Riccati equation \eqref{sdre1}, be low rank.

The general idea consists of first determining an approximation space ${\cal K}_k$ that can be naturally
expanded if needed, and then seeking an approximate solution in this space, by imposing a Galerkin condition on the matrix residual for computing the projected approximate solution.

Let ${\cal V}_k$ be such that ${\cal K}_k={\rm range}({\cal V}_k)$, with ${\cal V}_k$ having orthonormal columns. Recalling that in both
the Riccati and Lyapunov equations the solution matrix is symmetric, the approximate
solution can be written as  ${\cal V}_k Y {\cal V}_k^\top$. 

Consider the Riccati equation \eqref{sdre1} for a fixed $x=x(t_*)$, so that we set $A(x(t_*))=A_*$ and
$B(x(t_*))=B_*$. 
Let ${\cal R} = A_*^\top {\cal V}_k Y_X {\cal V}_k^\top + {\cal V}_k Y_X {\cal V}_k^\top A_* -
{\cal V}_k Y_X {\cal V}_k^\top B_* {\cal V}_k Y_X {\cal V}_k^\top + Q$.
Imposing the Galerkin condition on ${\cal R}$ means that the residual matrix ${\cal R}$ be orthogonal to 
the approximation space, in the matrix sense, that is
\begin{equation}\label{eqn:reduced_Riccati}
{\cal V}_k^\top {\cal R} {\cal V}_k =0 \quad 
\Leftrightarrow \quad
{\cal V}_k^\top A_*^\top {\cal V}_k Y_X +  Y_X {\cal V}_k^\top A_* {\cal V}_k -
 Y_X ({\cal V}_k^\top B_* {\cal V}_k) Y_X  + {\cal V}_k^\top Q {\cal V}_k=0,
\end{equation}
where the orthogonality of the columns of ${\cal V}_k$ was used.
If ${\cal V}_k$ has small dimensions, the {\it reduced} Riccati matrix equation 
on the right also has small dimensions and can be solved by a ``dense'' method
to determine $Y_X$; see, e.g., 
\cite{Binietal.book.12}.
The cost of solving the reduced quadratic equation with coefficient matrices of size $\hat k$ 
is at least $63 \hat k^3$ floating point operations with an invariant subspace approach
\cite{Laub.79}.
Note that the large matrix ${\cal V}_k Y_X {\cal V}_k^\top$ is never
constructed explicitly, since it would be dense even for sparse data.

Analogously, for the Lyapunov equation in \eqref{lyap_on}, we can write $W \approx {\cal V}_k Y {\cal V}_k^\top$
for some $Y=Y_W$ to be determined. Let ${\cal Q}_* = \sum_k Q_k f_k(x(t_*))$. As before, letting $\widetilde{\cal R} = C_0^\top {\cal V}_k Y {\cal V}_k^\top + {\cal V}_k Y {\cal V}_k^\top C_0
+ {\cal Q}_*$, the Galerkin condition leads to
$$
{\cal V}_k^\top \widetilde{\cal R} {\cal V}_k =0 \quad \Leftrightarrow \quad
 ({\cal V}_k^\top C_0^\top {\cal V}_k) Y  + Y {\cal V}_k^\top C_0 {\cal V}_k
+ {\cal V}_k^\top{\cal Q}_* {\cal V}_k = 0.
$$
This {\it reduced} Lyapunov equation can be solved by means of a ``dense'' method at a cost
of about 15 $\hat k^3$ floating point operations for coefficient matrices of size
$\hat k$, if the real Schur decomposition is employed; see, e.g.,
\cite{Simoncini.survey.16}. Note that the computational cost is significantly lower than
that of solving the reduced Riccati equation with matrices of the same size.

\subsection{On the selection of the approximation space}
 Choices as approximation space explored in the literature
include polynomial and rational Krylov subspaces
\cite{Simoncini.survey.16}. They both enjoy the property of being nested as they enlarge, that is
${\cal K}_k \subseteq {\cal K}_{k+1}$ where $k$ is associated with the space dimension.
Rational Krylov subspaces have emerged as the key choice because they are able to deliver
accurate  approximate solutions with a relatively small space dimension, compared with polynomial
spaces. Given a starting tall matrix $V_0$ and an invertible stable coefficient matrix
$A_*$, we have used two distinct rational spaces: the Extended Krylov subspace,
$$
{\cal E}{\cal K}_k = {\rm range}([V_0, A_*^{-1} V_0, A_* V_0, A_*^{-2} V_0, A_*^2 V_0, \ldots, 
A_*^{k-1} V_0, A_*^{-k} V_0]),
$$
which only involves matrix-vector products and solves with $A_*$, and the (fully) Rational Krylov subspace,
$$
{\cal R}{\cal K}_ k =
{\rm range}([V_0, (A_*- \sigma_2 I)^{-1} V_0, \ldots, \prod_{j=2}^k (A_*- \sigma_j I)^{-1} V_0]) .
$$
where $\sigma_j$ can be computed a-priori or adaptively.
In both cases, the space is expanded iteratively, one block of vectors at the time,
and systems with $A_*$ or with $(A_*- \sigma_j I)$ are solved by fast sparse methods.
For $A_*$ real valued and stable, the $\sigma_j$'s are selected to be in ${\mathbb C}^+$, so that
$A_*-\sigma_jI$ is nonsingular. The actual choice of the shifts is a key step and a rich
literature is available, yielding theoretically grounded effective strategies; see, e.g., the
discussion in \cite{Simoncini.survey.16}.

In our implementation we used the Extended Krylov subspace for solving the Lyapunov equation
in Algorithm \ref{Alg_sdreon}, which has several advantages, such as the computation of the sparse Cholesky factorization of $A_0$ once for all.
 On the other hand, we used the Rational Krylov subspace for the Riccati equation, which
has been shown to be largely superior over the Extended Krylov on this quadratic equation,
in spite of requiring the solution of a different (shifted) sparse coefficient matrix at each iteration 
\cite{Benner.Saak.survey13,SSM14}.
Nonetheless, in Section~\ref{sec:feedback_matrix} we report an alternative approach that
makes the Extended Krylov subspace competitive again for the Riccati problem with $A_0$ symmetric.
Except for the operations associated with the reduced problems, the computational costs per iteration of the
Riccati and Lyapunov equation solvers are very similar, if the same approximation space is used.

Although we refer to the specialized literature for the algorithmic 
details\footnote{See {\tt www.dm.unibo.it/$\widetilde{\mbox{ }}$simoncin/software} for some related software.},
we would like to include some important implementation details that are specific to our
setting. 
In particular, the matrix $C_0$ employed in Algorithm \ref{Alg_sdreon} is given by
$C_0 = A_0 - \left(BR^{-1} B^\top - \frac 1 {2\gamma^2} H P^{-1} H^\top\right) \Pi_0$, which is
not easily invertible if explicitly written, since it is dense in general. Note that
$A_0$ is in general sparse, as it stems from the discretization of a partial differential
operator. We can write
$$
C_0 = A_0 - \left([B, H] G [B,H]^\top\right) \Pi_0, \quad {\rm with} \quad
G = \begin{pmatrix} 
R^{-1}& 0\\
0 & - \frac 1 {2\gamma^2} P^{-1}
\end{pmatrix}.
$$
Using the classical Sherman-Morrison-Woodbury formula, the product $C_0^{-1} V$ for some
tall matrix $V$ can be obtained as 
$$
W:=C_0^{-1} V = A_0^{-1} V - A_0^{-1} [B, H] G_1^{-1} [B,H]^\top \Pi_0 A_0^{-1}V,
$$
with $G_1 = I + G [B,H]^\top \Pi_0 A_0^{-1} [B, H]$, which is assumed to be nonsingular. 
Therefore, $C_0^{-1} V$  can be obtained
by first solving sparse linear systems with $A_0$, and then using matrix-matrix products.
More precisely, the following steps are performed:

\vskip 0.1in

- Solve $A_0 W_1 = V$

- Solve $A_0 W_2 = [B, H]$

- Compute $G_1 = I + G [B,H]^\top \Pi_0 W_2$

- Compute $W = W_1 - W_2 ( \, G_1^{-1} ( [B,H]^\top \Pi_0 W_1) \, )$

\vskip 0.1in
We also recall that $\Pi_0$, the solution to the initial Riccati equation, is stored
in factored form, and this should be taken into account when computing matrix-matrix
products with $\Pi_0$.

While trying to employ the Rational Krylov space, we found that the structure of $C_0$ 
made the selection of optimal shifts $\{\sigma_j\}$ particularly challenging, resulting
in a less effective performance of the method. Hence, our preference went for
the Extended Krylov space above, with the above
enhancement associated with solves with $C_0$.

\subsection{Feedback matrix oriented implementation}\label{sec:feedback_matrix}
In Algorithm \ref{Alg_sdre} the Riccati equation needs to be solved at each
time step $t_n$. However, its solution $\Pi(x_n)$ is only used to compute the feedback matrix 
$K(x_n) : = - R^{-1} B^\top \Pi(x_n)$. Hence, it would be desirable to be able to
immediately compute $K(x_n)$ without first computing $\Pi(x_n)$. This approach
has been explored in the Riccati equation literature but also for
other problems based on Krylov 
subspaces, see, e.g., \cite{Palitta.Simoncini.18a,Kressner.proc.2008}.

In this section, in the case where the matrix $A(x_n)$ is symmetric,
 we describe the implementation of one of the projection methods described above,
that is able to directly compute $K(x_n)$ without computing $\Pi(x_n)$ (in factored form), and
more importantly, without storing and computing the whole approximation basis. 
The latter feature is particularly important for large scale problems, for which dealing
with the orthogonal approximation basis represents one of the major computational and
memory costs. To the best of our knowledge, this variant of the Riccati solver is new, while it is currently explored in \cite{Palitta.Pozza.Simoncini.21} for related control problems and the rational Krylov space.

Here we consider the Extended Krylov subspace. For $A(x_n)$ symmetric, the 
orthonormal basis of ${\cal E}{\cal K}_k$
can be constructed by explicitly orthogonalizing only with respect
to the previous two basis blocks. Hence, only two previous blocks of vectors need to
be stored in memory, and require explicit orthogonalization when the new block of vectors
is added to the basis \cite{Palitta.Simoncini.18a}. This is also typical of polynomial
Krylov subspaces constructed for symmetric matrices, giving rise to the classical Lanczos three-term recurrence \cite{Kressner.proc.2008}. 

With this procedure, in the reduced equation \eqref{eqn:reduced_Riccati} 
the matrices $A_k := {\cal V}_k^\top A_* {\cal V}_k$, 
$B_k := {\cal V}_k^\top B$ and $Q_k:={\cal V}_k^\top Q {\cal V}_k$
are computed as $k$ grows by updating the new terms at each iteration, and
the solution $Y_X$ can be obtained without the whole matrix ${\cal V}_k$ being available.
Note that the stopping criterion does not require the computation of the whole residual
matrix, so that also in the standard solver
the full matrix ${\cal V}_k Y_X {\cal V}_k^\top$ is never explicitly accessed.

However, to be able to compute 
$$
K = -R B^\top {\cal V}_k Y_X {\cal V}_k^\top = -R B_k^\top Y_X {\cal V}_k^\top,
$$
the basis ${\cal V}_k$ appearing on the right still seems to be required.
As already done in the literature, this problem can be overcome by a so-called ``two-step'' 
procedure: at completion, once the final $Y_X$ is available,
 the basis ${\cal V}_k$ is computed again one block at the
time, and the corresponding terms in the product $Y_X {\cal V}_k^\top$ are updated.
Since $A(x_n)$ is already factorized and the orthogonalization coefficients are already
available (they correspond to the non-zero entries of $A_k$), then the overall computational
cost is feasible; we refer the reader to \cite{Palitta.Simoncini.18a} and its references
for additional details for the two-step procedure employed for different purposes.

\section{Large-scale nonlinear dynamical systems}\label{numerics}

In this section we present a numerical assessment of the proposed methodology applied to the synthesis of feedback control for two-dimensional nonlinear PDEs. The first test is a nonlinear diffusion-reaction equation, known as the degenerate Zeldovich equation, where the origin is an unstable equilibrium and traditional linearization-based controllers fail. The second test studies the viscous Burgers' equation with a forcing term. We  discretize the control problem in space using finite differences, similarly as in Section \ref{sinegordon}. Controlled trajectories are integrated in time using an implicit Euler method, which is accelerated using a Jacobian--Free Newton Krylov method (see e.g. \cite{KK04}). The goal of all our tests is the optimal and robust stabilization of the dynamics to the origin, encoded in the optimization of the following cost:

\begin{align}\label{costex1}
\begin{aligned}
\cJ(u(\cdot),w(\cdot);X(\cdot,0))&:=\int\limits_0^{\infty} \sum_{i=1}^z \dfrac{1}{|\omega_{o_i}|}\left(\int_{\omega_{o_i}} X(\xi,t)\, d\xi\right)^2\\
&\quad\quad+R|u(t)|^2 -\gamma^2 P|w(t)|^2\, \, dt\,.
\end{aligned}
\end{align}
 {This expression is similar to \eqref{costex_wave}, but includes the $\Hinf$ term $-\gamma^2 P|w(t)|^2$.}

The reported numerical simulations were performed on a MacBook Pro 
with CPU Intel Core i7-6, 2,6GHz and 16GB RAM, using Matlab \cite{matlab7}.

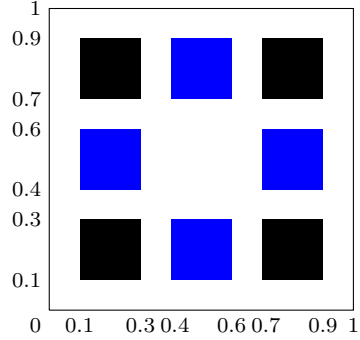
\begin{figure}[htbp]
\centering
\begin{tikzpicture}[scale=4]
\draw (0,0) -- (1,0) -- (1,1) -- (0,1) -- (0,0); 
\fill[black] (0.1,0.1) -- (0.3,0.1) -- (0.3,0.3) -- (0.1,0.3) -- (0.1,0.1);
\fill[black] (0.7,0.7) -- (0.9,0.7) -- (0.9,0.9) -- (0.7,0.9) -- (0.7,0.7);
\fill[black] (0.1,0.7) -- (0.3,0.7) -- (0.3,0.9) -- (0.1,0.9) -- (0.1,0.7);
\fill[black] (0.7,0.1) -- (0.9,0.1) -- (0.9,0.3) -- (0.7,0.3) -- (0.7,0.1);
\fill[blue]  (0.1,0.4) -- (0.3,0.4) -- (0.3,0.6) -- (0.1,0.6) -- (0.1,0.4);
\fill[blue] (0.1,0.4) -- (0.3,0.4) -- (0.3,0.6) -- (0.1,0.6) -- (0.1,0.4);
\fill[blue] (0.7,0.4) -- (0.9,0.4) -- (0.9,0.6) -- (0.7,0.6) -- (0.7,0.4);
\fill[blue] (0.4,0.1) -- (0.6,0.1) -- (0.6,0.3) -- (0.4,0.3) -- (0.4,0.1);
\fill[blue] (0.4,0.7) -- (0.6,0.7) -- (0.6,0.9) -- (0.4,0.9) -- (0.4,0.7);
\draw (0,0) node[below left] {$0$};
\draw (0.1,0) node[below ] {$0.1$};
\draw (0.3,0) node[below ] {$0.3$};
\draw (0.4,0) node[below ] {$\,\,0.4$};
\draw (0.6,0) node[below ] {$0.6$};
\draw (0.7,0) node[below ] {$\,\,0.7$};
\draw (0.9,0) node[below ] {$0.9$};
\draw (1,0) node[below ] {$1$};
\draw (0,0.1) node[left ] {$0.1$};
\draw (0,0.3) node[left ] {$0.3$};
\draw (0,0.4) node[left ] {$0.4$};
\draw (0,0.6) node[left] {$0.6$};
\draw (0,0.7) node[left] {$0.7$};
\draw (0,0.9) node[left ] {$0.9$};
\draw (0,1) node[left ] {$1$};
\end{tikzpicture}
\caption{Locations of the inputs 
$\omega_c(\xi)=\omega_d(\xi)$  (black) and outputs $\omega_o(\xi)$ (blue) in the region $\Omega$ for the degenerate Zeldovich equation.}
\label{fig:bc_coll}
\end{figure}

\subsection{Case study 1: the degenerate Zeldovich equation}\label{sec:test1}

We consider the control of a Zeldovich-type equation arising in combustion theory \cite{zeldo} over $\Omega\times\R^+_0$, with $\Omega\subset \R^2$ and Neumann boundary conditions:
\begin{align}\label{ex1}
\begin{aligned}
\pt X(\xi,t)&=\epsilon\Delta X(\xi,t)+\nu X(\xi,t)+\mu(X^2(\xi,t)-X^3(\xi,t))\\
&\qquad +\chi_{\omega_c}(\xi)u(t)+\chi_{\omega_d}(\xi)w(t)\\
\pxi X(\xi,t)&=0\,, \quad\xi\in\partial\Omega,\, t>0\,,\\
 X(\xi,0)&={\mathbf x}_0(\xi)\,,\quad \xi\in\Omega\,.
\end{aligned}
\end{align}
The scalar control and disturbance act, respectively, through functions $\chi_{\omega_c}(\xi)$ and $\chi_{\omega_d}(\xi)$ with support $\omega_c,\omega_d\subset\Omega$. The uncontrolled dynamics have three equilibrium points: $X\equiv 0$, $X\equiv\frac12\left(1\pm\sqrt{1+\frac{\nu}{\mu}}\right)$. Our goal is to stabilize the system to $X\equiv 0$, which is an unstable equilibrium point. A first step towards the application of the proposed framework is the space discretization of the system dynamics, leading to a finite-dimensional state-space representation. 
Following the setting presented in section~\ref{sinegordon}, using a finite difference discretization leads to
\begin{equation}\label{pol:sd}
\dot X(t) = \epsilon\Delta_dX(t)+\nu X(t)+\mu X(t)\circ X(t) \circ \left(\textbf{1}_{d\times 1}-X(t)\right) +B u(t) + H w(t)\,,
\end{equation}
where the discrete state $X(t)=(X_1(t),\ldots,X_d(t))^\top\in\R^d$ corresponds to the approximation of $X(\xi, t)$ at the grid points and the symbol $\circ$ denotes the Hadamard or component-wse product. The matrix $\Delta_d\in\R^{d\times d}$ is the finite difference approximation of the  
Neumann Laplacian and the matrices $B,H\in \R^d$ are the discretization of the indicator functions supported over $\omega_c$ and $\omega_d$, respectively. 
{The discretization of \eqref{costex1} follows similarly as in section~\ref{sinegordon}.}
Once the finite-dimensional state-space representation is obtained, we proceed to express the system in semilinear form \eqref{sdrea} and implement the proposed algorithms. 
To set Algorithm \ref{Alg_sdre} and Algorithm \ref{Alg_sdreon}, from \eqref{pol:sd} we define
$$A(X) := \epsilon\Delta_d +\nu\mathbf{I}_{d}+ \mu \texttt{diag}(X(t)-X(t)\circ X(t)),$$
where $\texttt{diag}(v)$ denotes a diagonal matrix with the components of the vector $v$ on the main
diagonal, and decompose $A(X)$ as
$$
A_0 = \epsilon\Delta_d +\nu\mathbf{I}_{d},\qquad [A_j]_{k,l}=\delta_{k,j}\delta_{l,j},\qquad 
f_j(X)=\mu(X_j-X_j^2),\qquad j=1\ldots,d,
$$
where $\Delta_d$ is the two-dimensional discrete Laplacian and $\delta_{i,j}$ denotes the Kronecker delta.
In this test we set
$$\Omega = [0,1]\times[0,1],\,\, \epsilon=0.2,\,\, \nu = 0.1,\,\, \mu=10,\,\, R = 0.1,$$  
and the initial condition ${\mathbf x}(\xi,0) = \sin (\xi_1)\sin(\xi_2)$, on a discretized space grid of $n_{\xi_1}\times n_{\xi_2}$ nodes with $n_{\xi_1} = n_{\xi_2} = 101$ ($d=10201$). For the matrices $B, H$ and $C$ we considered a collection of patches depicted in  Figure \ref{fig:bc_coll}, and given by
$$
\omega_d(\xi) = [0.1,0.3]^2 \cup [0.7,0.9]^2 \cup 
\left([0.1,0.3]\times[0.7,0.9]\right) \cup \left( [0.7,0.9]\times[0.1,0.3]\right),
$$
$\omega_c(\xi)=\omega_d(\xi)$, and
\begin{align*}
\begin{aligned}
\omega_o(\xi)  =&  \left([0.1,0.3]\times[0.4, 0.6]\right) \cup  \left([0.4,0.6]\times[0.1, 0.3]\right) \\
&\cup \left([0.4,0.3]\times[0.7,0.9]\right) \cup \left( [0.7,0.9]\times[0.4,0.6]\right)\,.
\end{aligned}
\end{align*}
In the following, we analyse results for the $\Htwo$ and  $\Hinf$ controls, i.e. $\gamma = 0$ and $\gamma\neq 0$, respectively, in \eqref{costex1}.

\begin{figure}[htbp]
\centering
\includegraphics[scale=0.29]{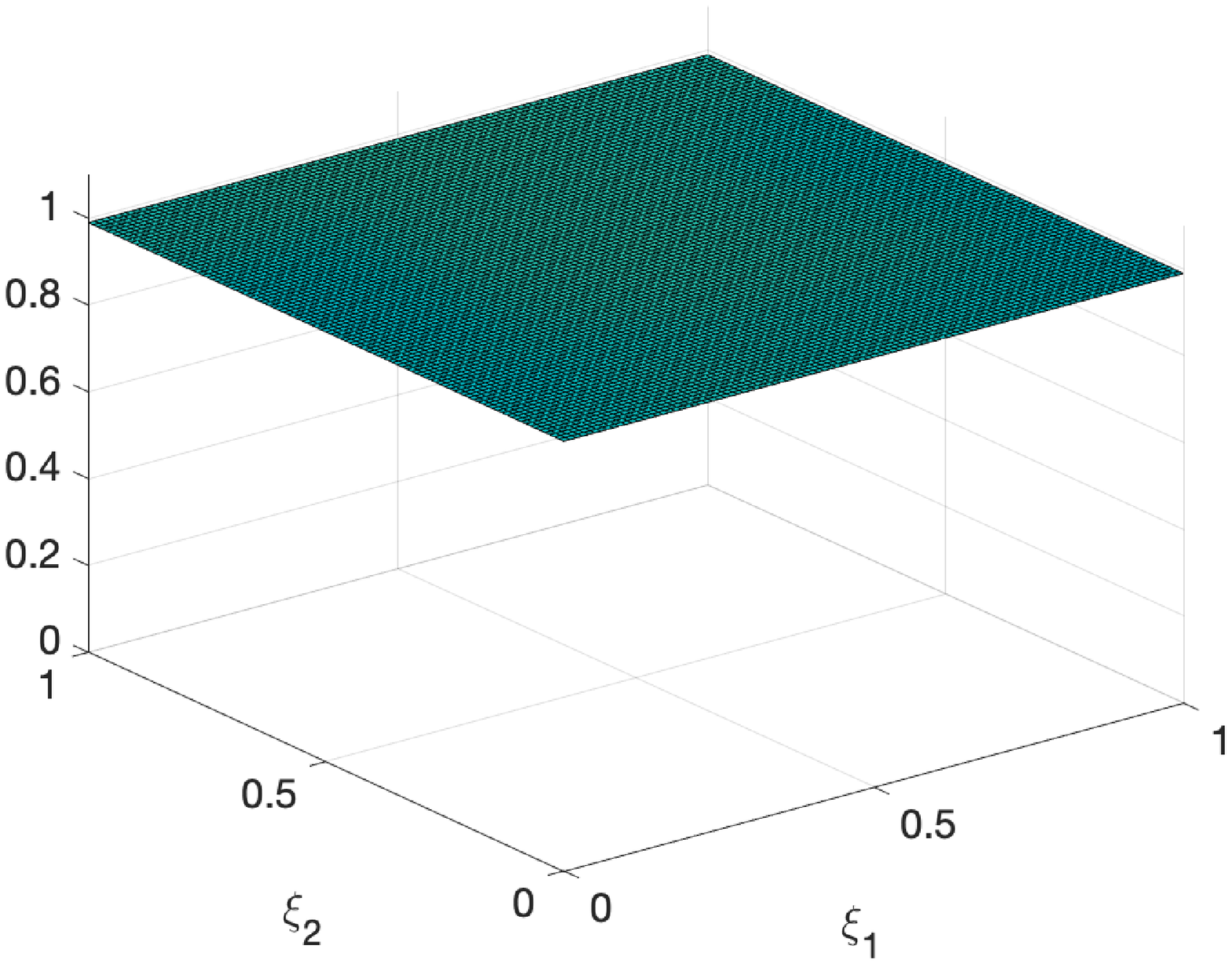}
\includegraphics[scale=0.29]{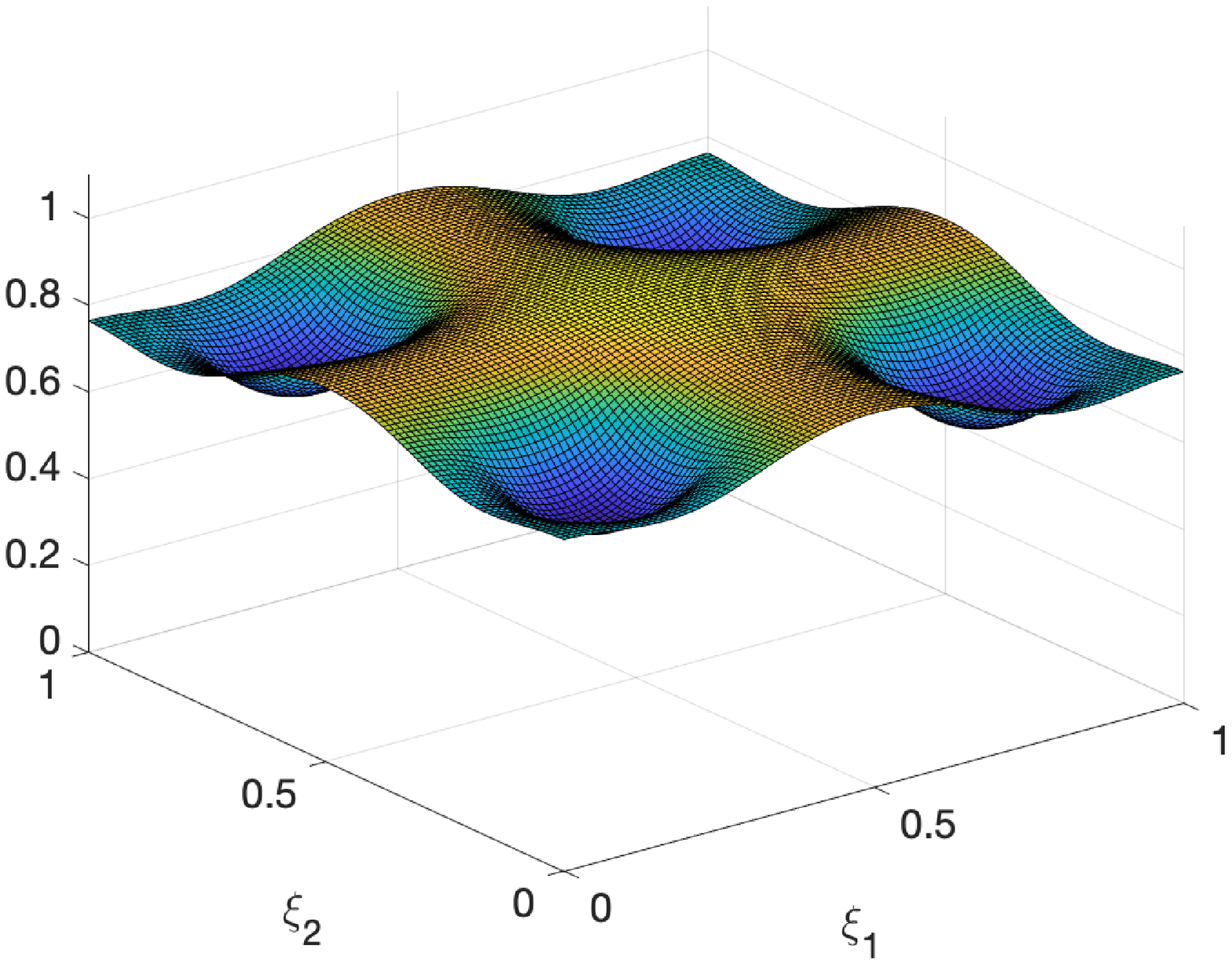}
\includegraphics[scale=0.29]{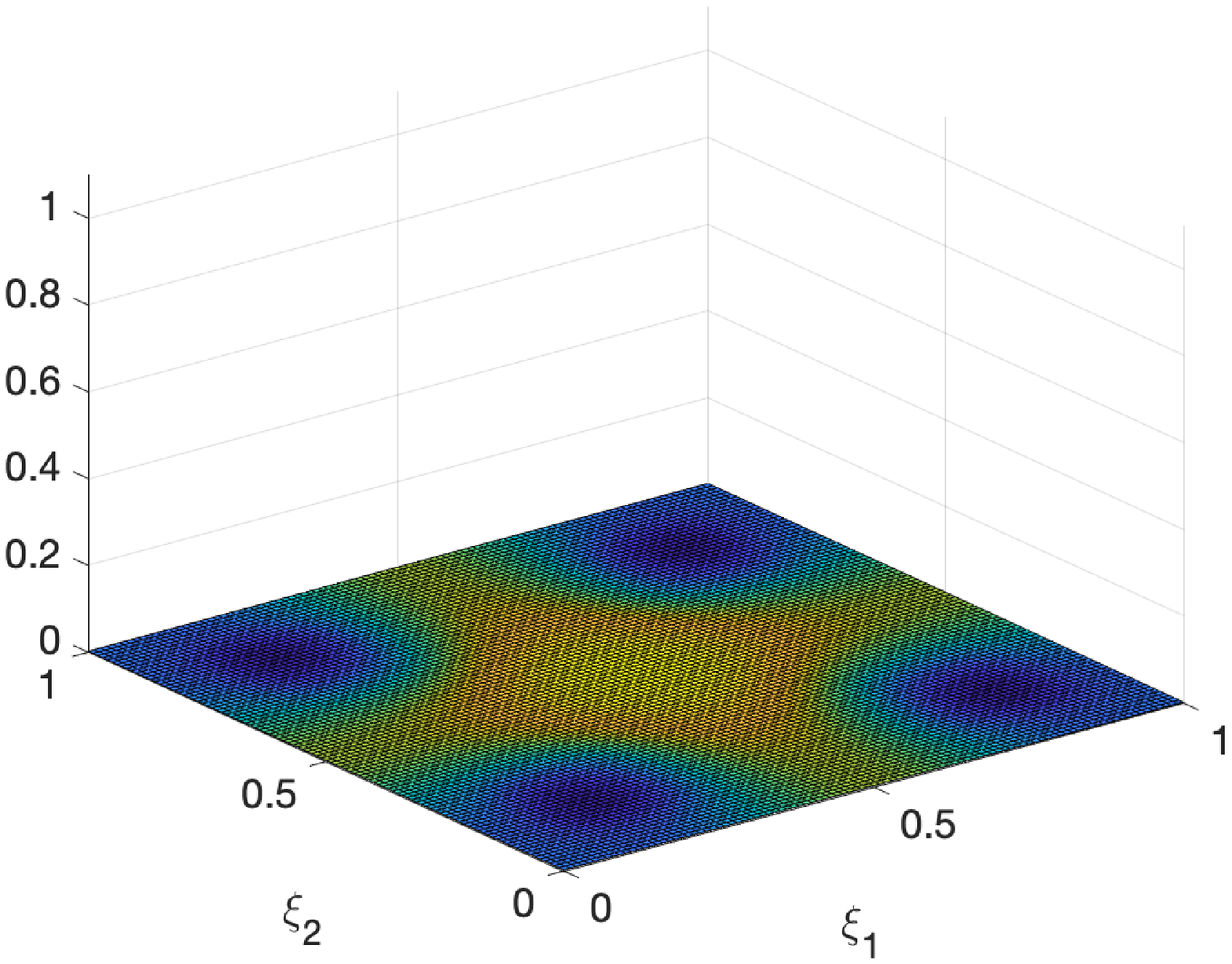}
\includegraphics[scale=0.29]{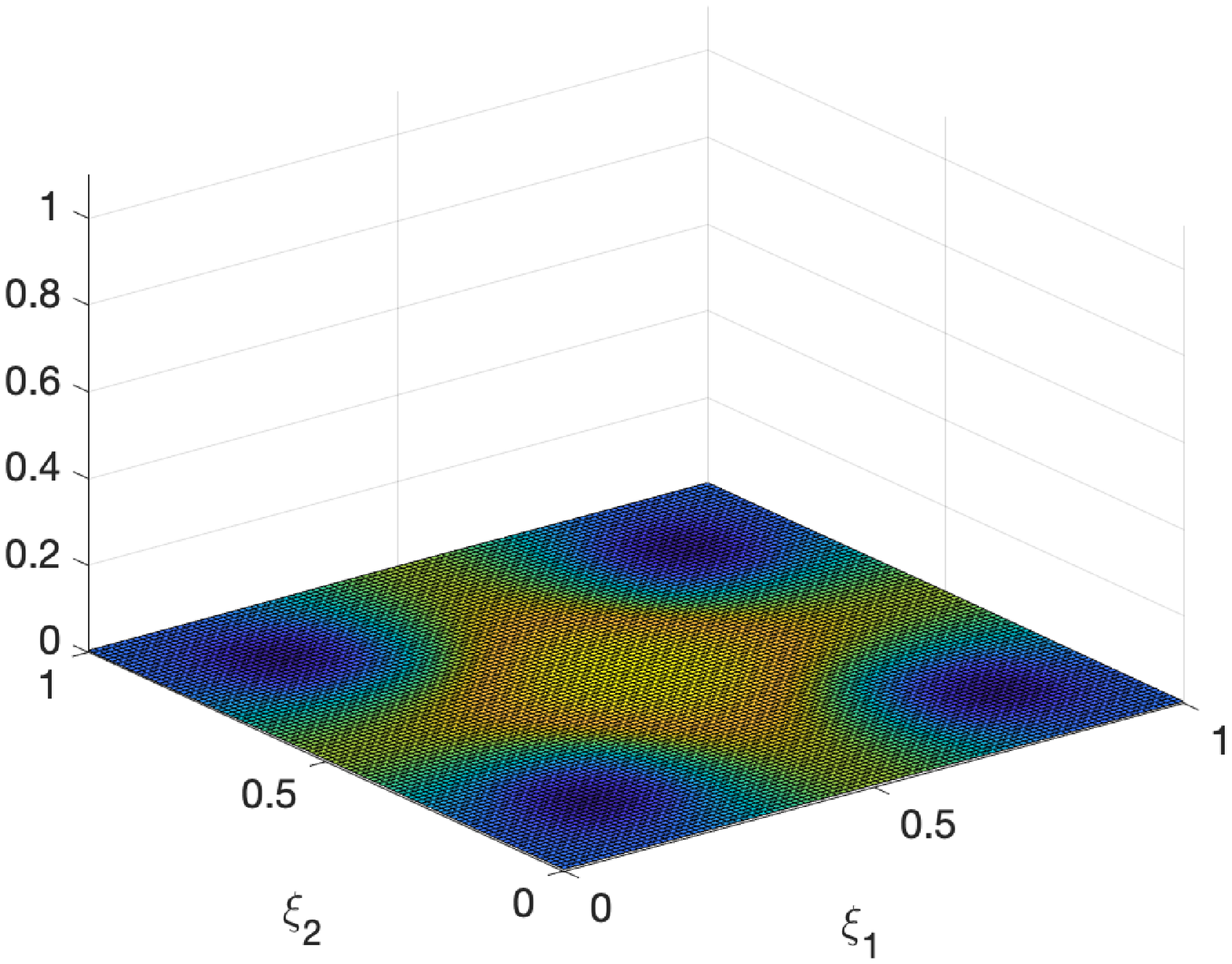}
\caption{Test \ref{sec:test1h2}: state of the system 
\eqref{ex1} at time $t=3$. Top: uncontrolled solution (left), $\mathcal{H}_2$-solution with LQR control (right). Bottom: $\mathcal{H}_2$-controlled solution with Algorithm \ref{Alg_sdre} (left), $\mathcal{H}_2$-controlled solution with Algorithm \ref{Alg_sdreon} (right).} 
 \label{fig1:sol}
\end{figure}
\begin{figure}[htbp]
\centering
\includegraphics[scale=0.3]{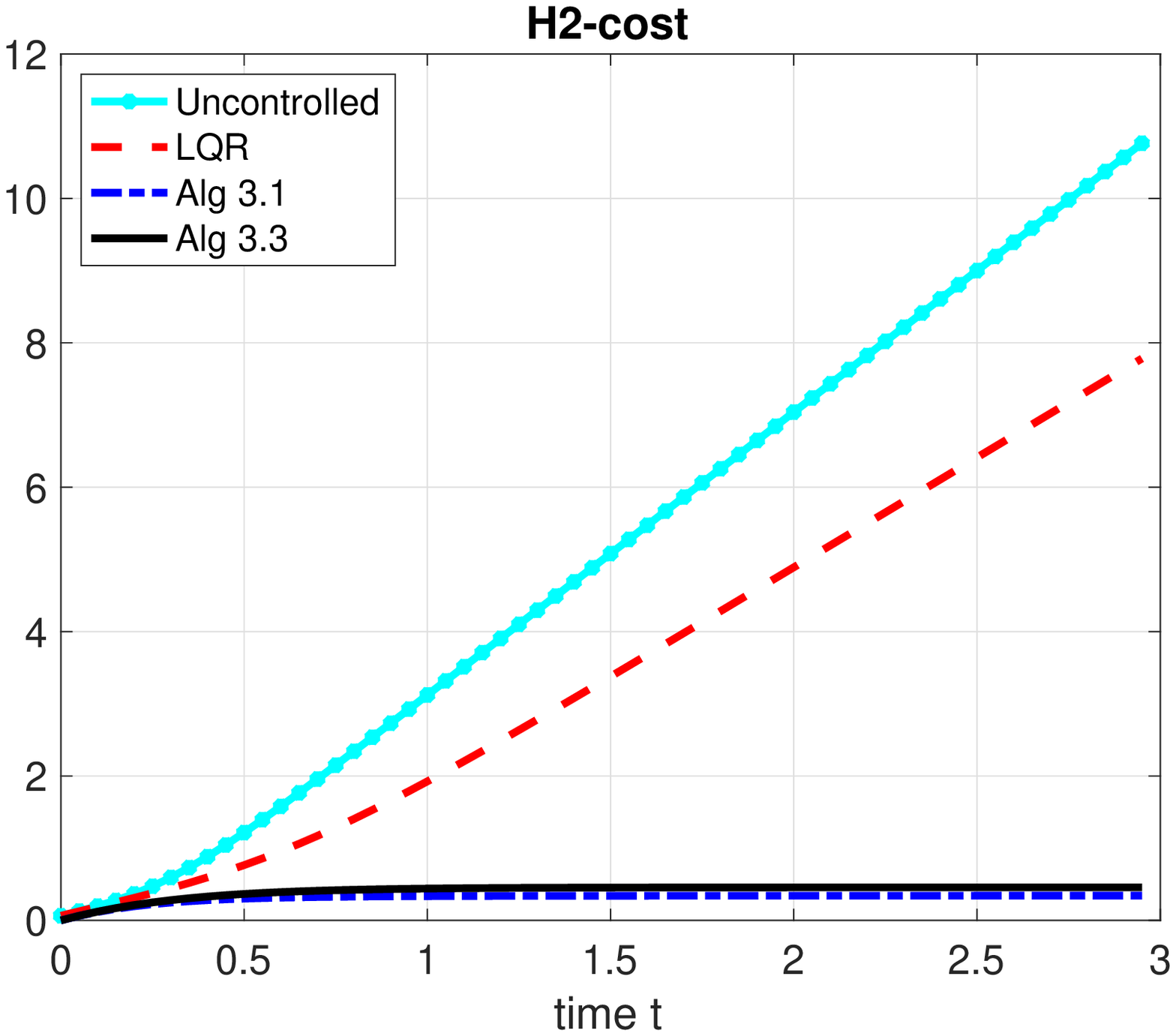}
\includegraphics[scale=0.3]{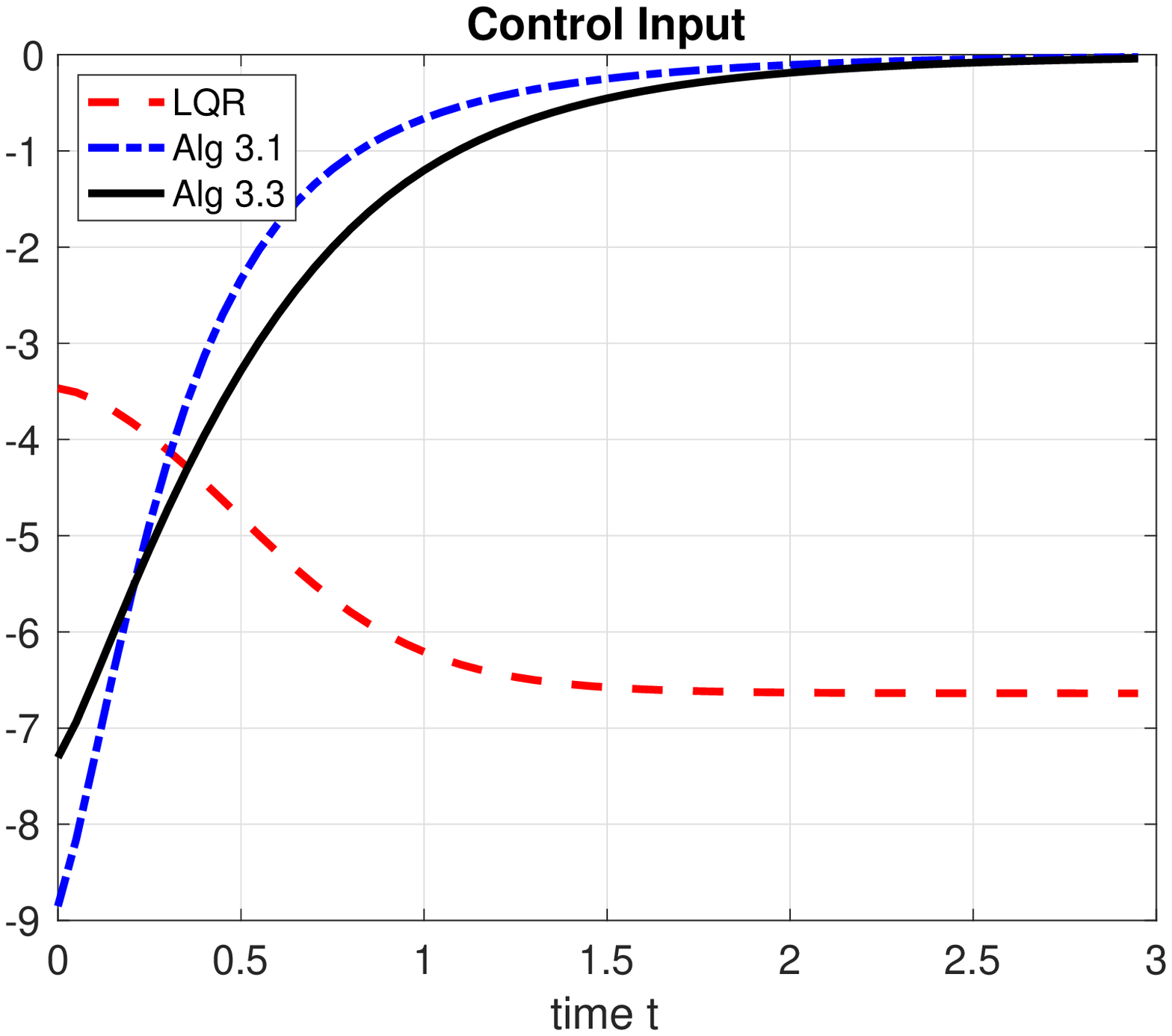}
\includegraphics[scale=0.3]{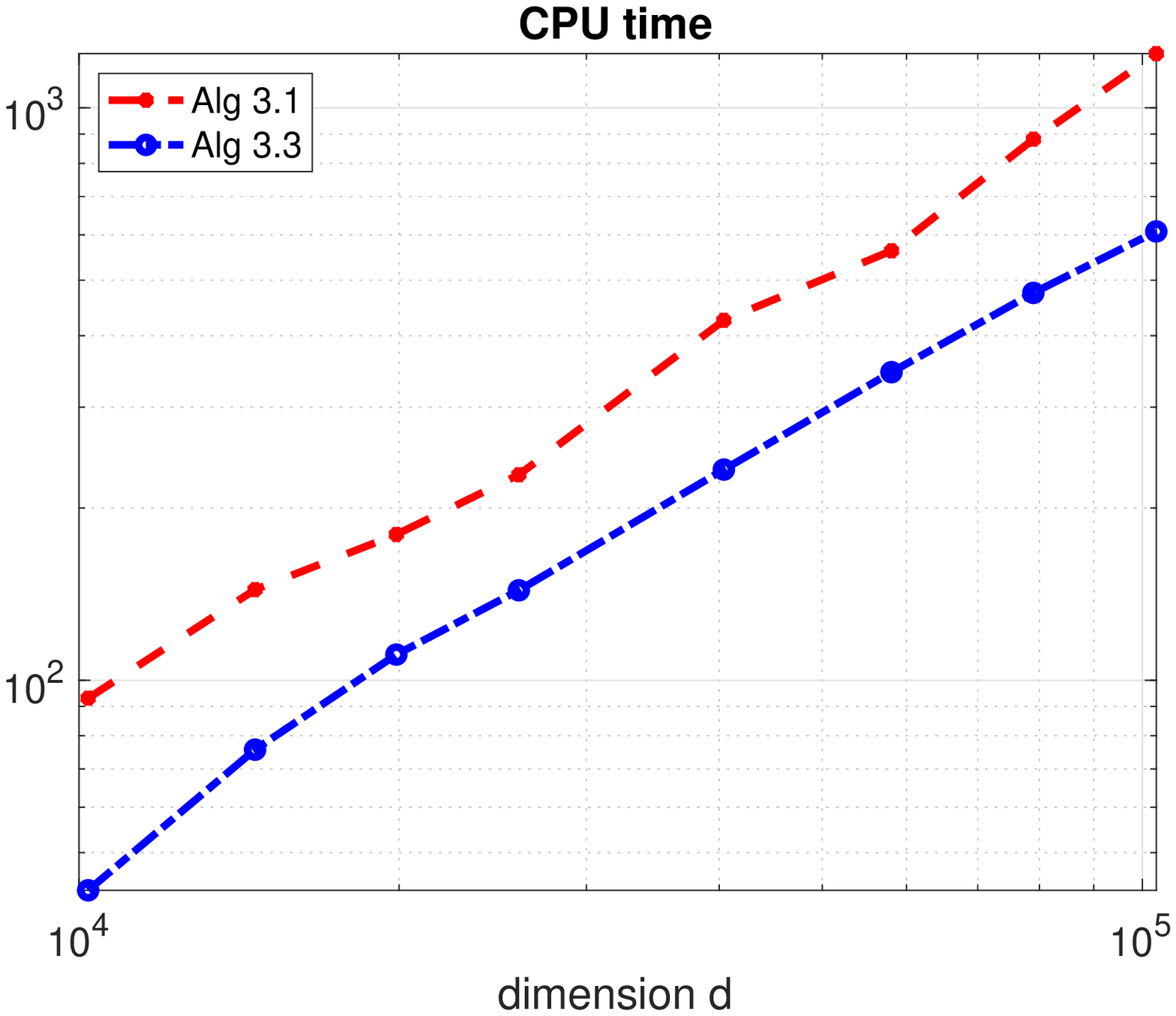}
\includegraphics[scale=0.3]{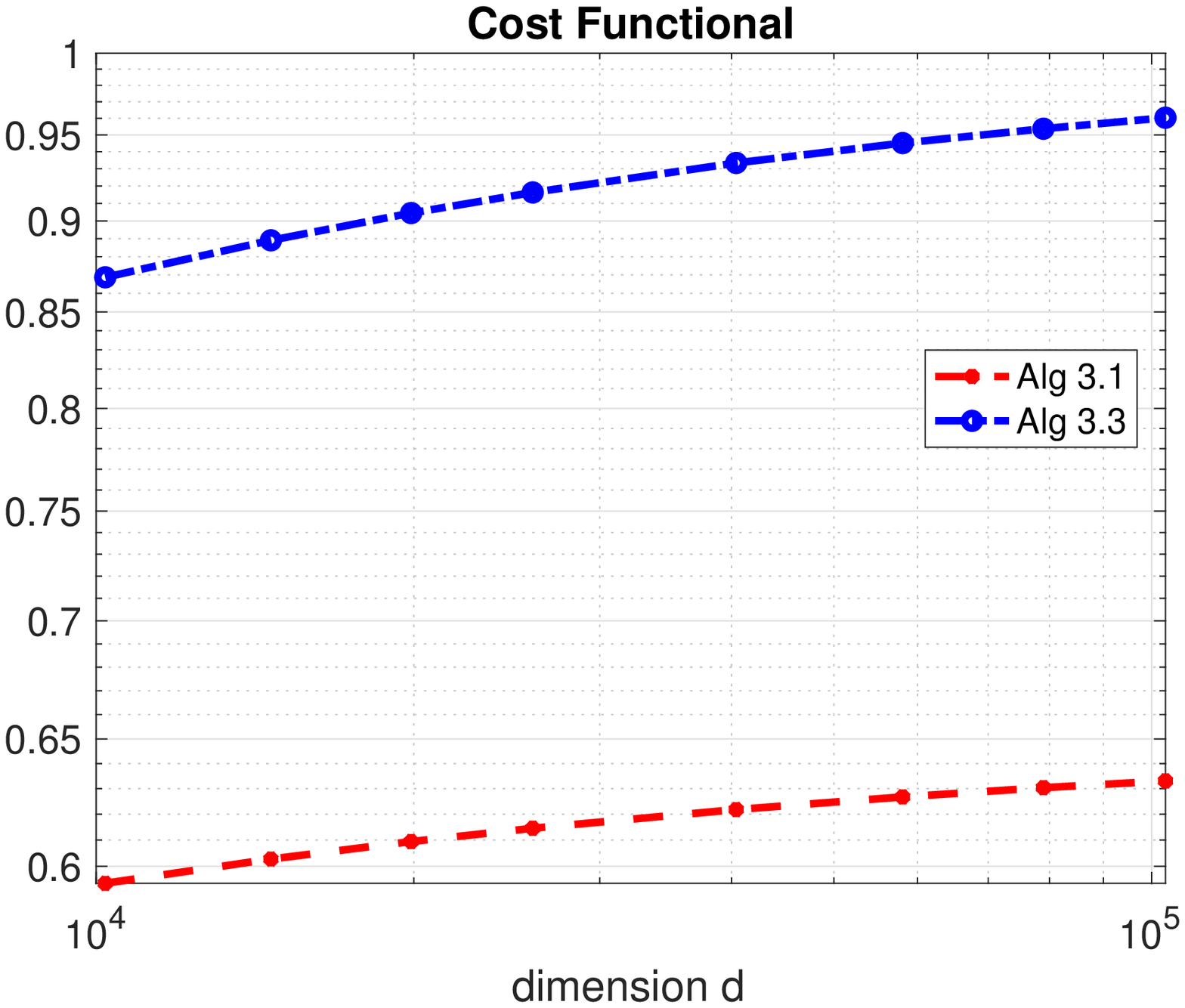}
\caption{Test \ref{sec:test1h2}. Top:  Cumulative cost functional  functional (left) with $\Htwo-$control and corresponding control inputs (right). Bottom: CPU time for Algorithm \ref{Alg_sdre} and Algorithm \ref{Alg_sdreon} (left), and convergence of the cost functional with respect to the dimension of the problem $d$ ($x$-axis) (right).}
\label{fig1:cost}
\end{figure}

\begin{example}{\it Experiments for $\Htwo$-control.}\label{sec:test1h2} 
We start by presenting results for $\Htwo$-control, i.e. $P \equiv 0$ in \eqref{costex1} and $H\equiv 0$ in \eqref{pol:sd}. In Figure \ref{fig1:sol} we show a snapshot of the controlled trajectories at $t=3$, a horizon sufficiently large for the dynamics to approach a stationary regime. In the top-left panel the uncontrolled problem reaches the stable equilibrium $X\approx 1.02$. In the top-right panel we show the results of the LQR control computed by linearizing equation \eqref{pol:sd} around the origin, which also fails to stabilize around the unstable equilibrium $X \equiv 0$. The controlled solutions with Algorithm \ref{Alg_sdre} and Algorithm \ref{Alg_sdreon} are shown at the bottom of the same figure: both algorithms reach the desired configuration. The corresponding control input is shown in the top-right panel of Figure \ref{fig1:cost}. We observe that the LQR control has a completely different behavior with respect to the control computed by Algorithm \ref{Alg_sdre} or Algorithm \ref{Alg_sdreon}.

The performance results of the different controlled trajectories is presented in Figure \ref{fig1:cost} where we show the evaluation of the cumulative cost functional in the top-left panel. As expected, Algorithm \ref{Alg_sdre} provides the best closed-loop performance among the proposed algorithms. However, in terms of efficiency Algorithm \ref{Alg_sdreon} is faster than Algorithm \ref{Alg_sdre} when increasing the dimension of the problem as shown in the bottom-left panel of Figure \ref{fig1:cost}. When the dimension $d$ increases ($x$-axis in the plot), the cost functional converges to $0.6$ for Algorithm \ref{Alg_sdre} and to $1$ for Algorithm \ref{Alg_sdreon}. Both methods are able to stabilize the problem.

\end{example}

\begin{figure}[htbp]
\centering
\includegraphics[scale=0.29]{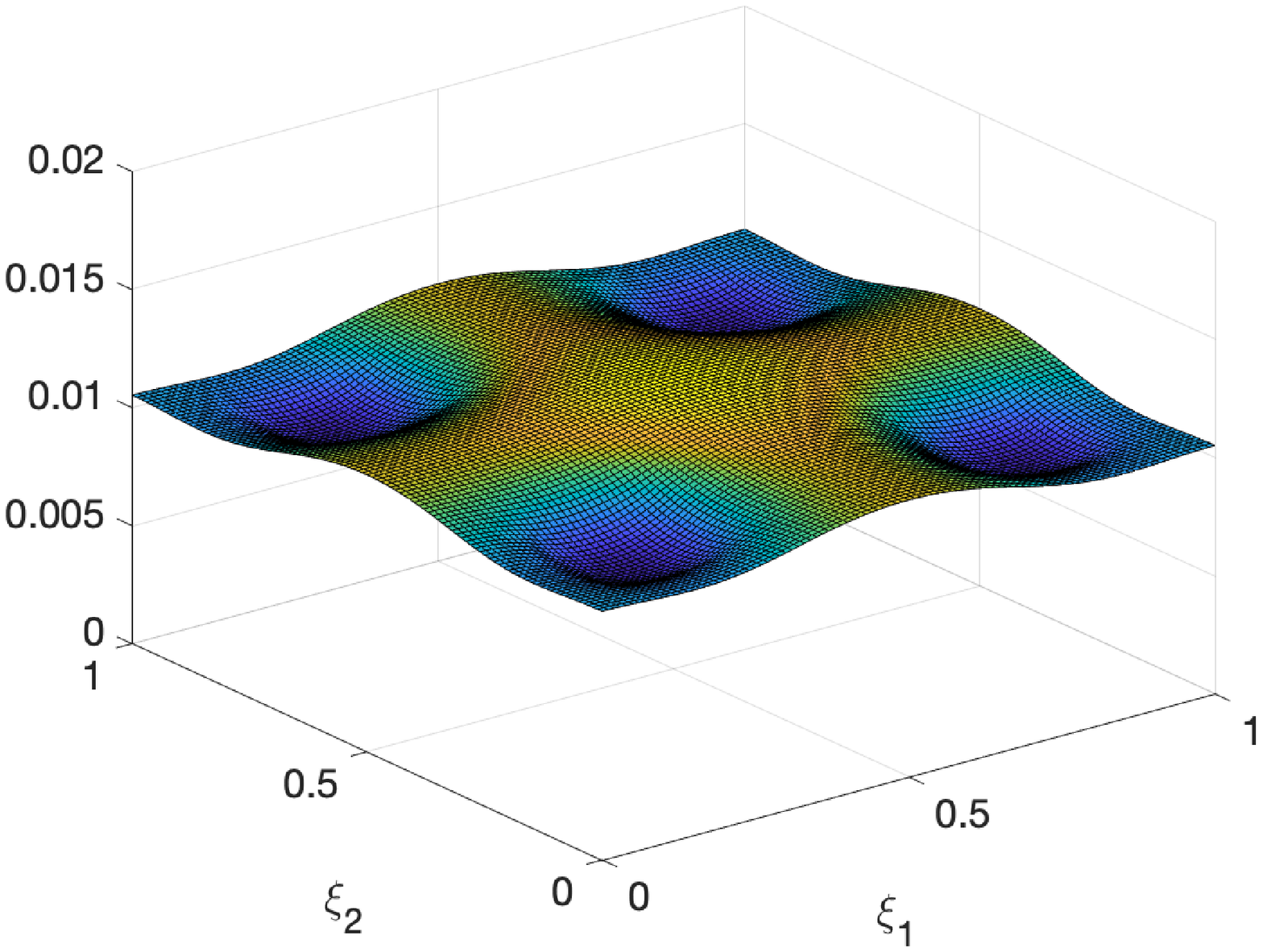}
\includegraphics[scale=0.29]{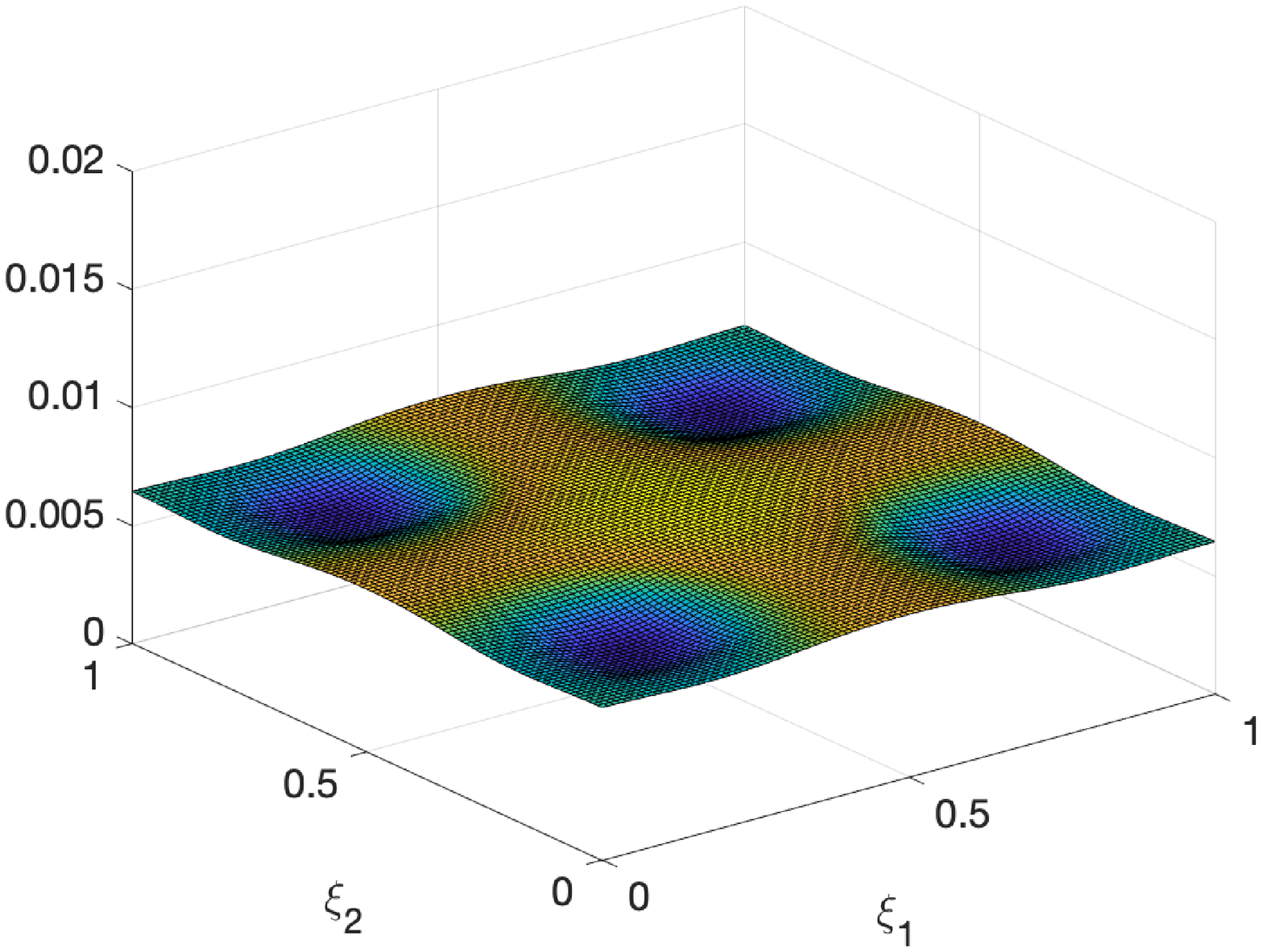}
\includegraphics[scale=0.29]{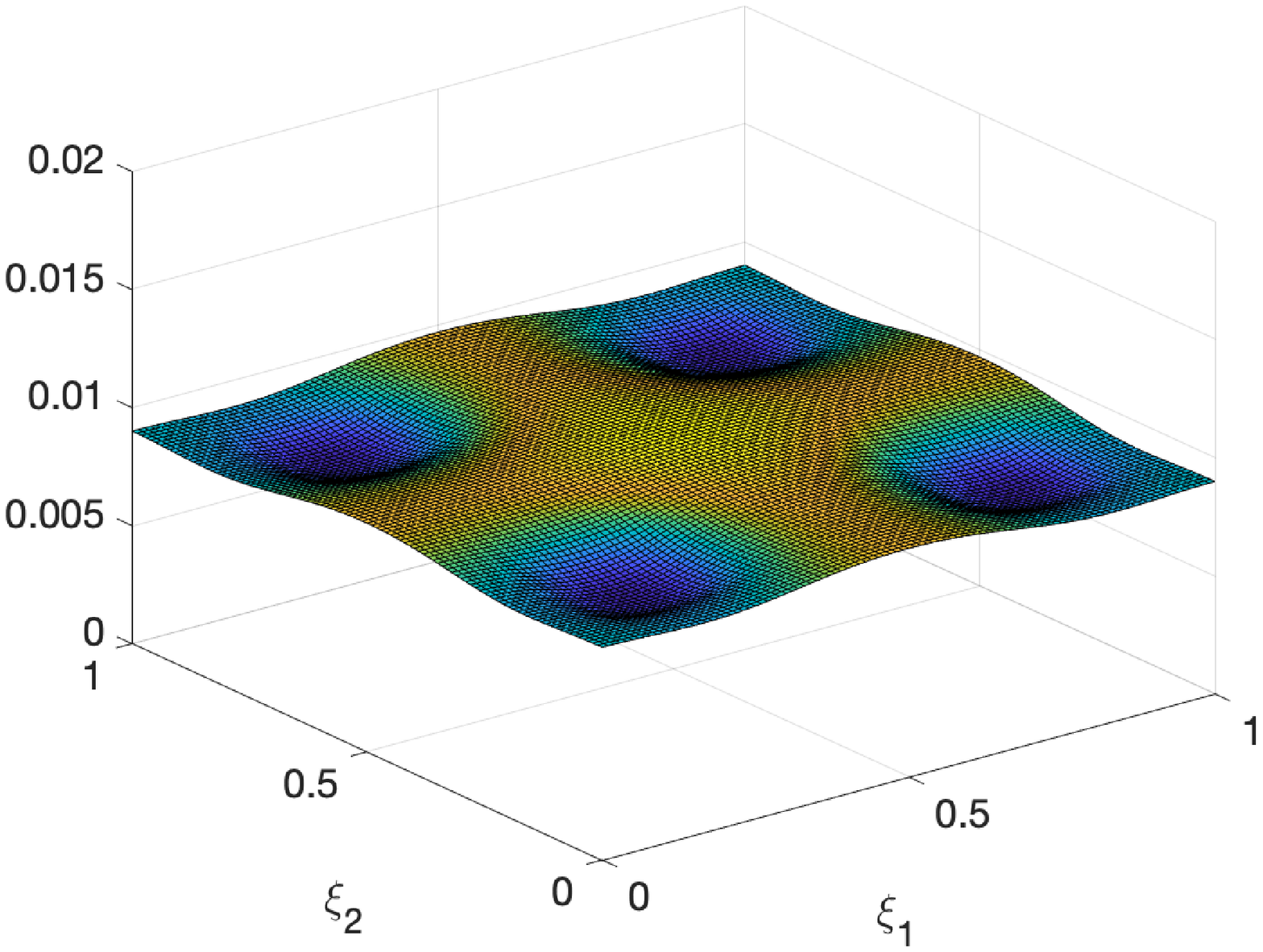}
\includegraphics[scale=0.29]{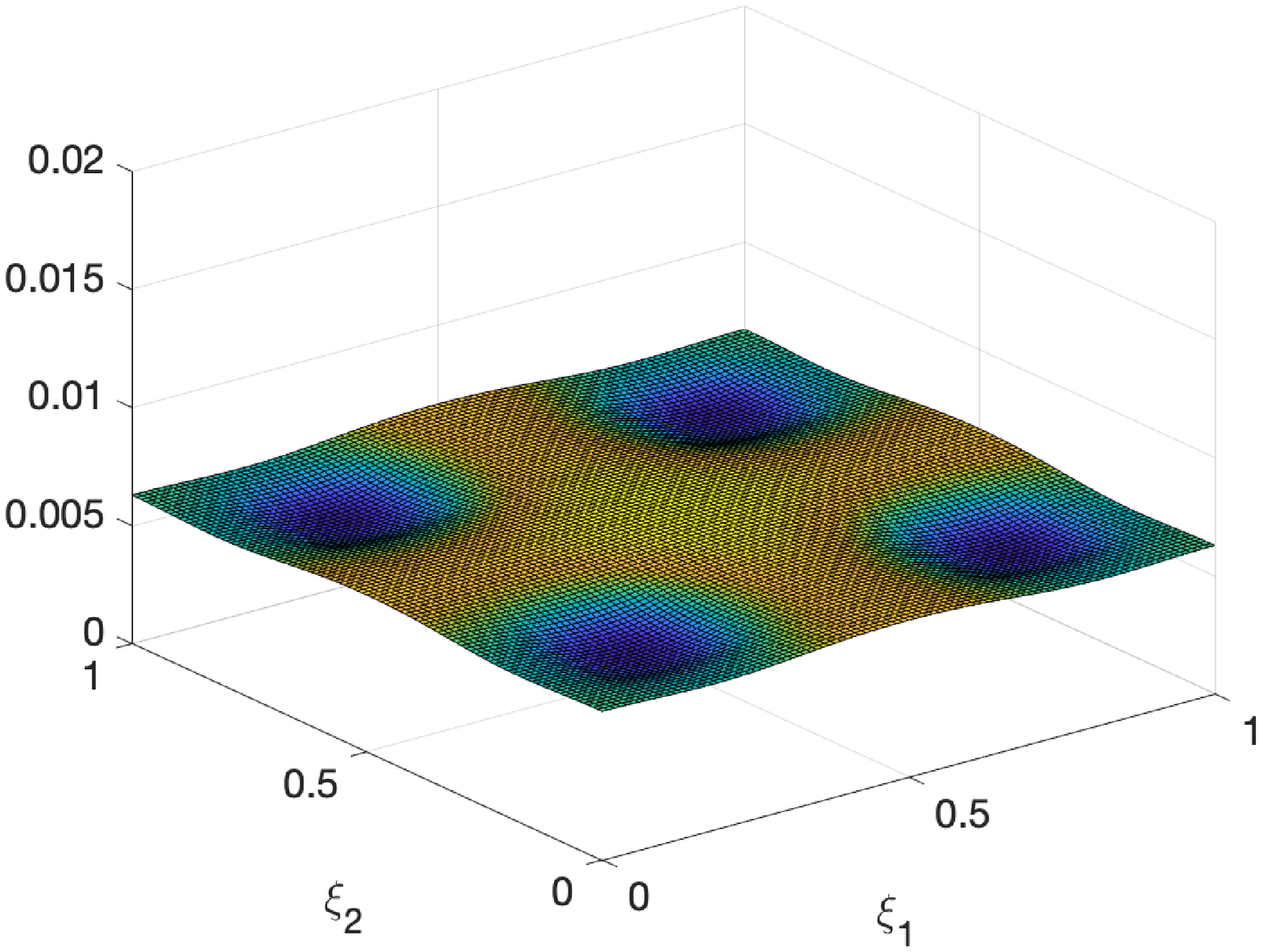}
\caption{Test \ref{sec:ex2}: state of the system at time $t=3$ with $\Htwo-$ and $\Hinf-$controls and disturbance $w(t)=0.1\sin(2t)$.
Top:  $\Htwo$-control using Algorithm \ref{Alg_sdre} (left) and Algorithm \ref{Alg_sdreon} (right). Bottom:  $\Hinf$-control using Algorithm \ref{Alg_sdre} (left) and Algorithm \ref{Alg_sdreon} (right).}  
\label{fig1:dist1}
\end{figure}

\begin{example}{\it Experiments for $\Hinf$-control.}\label{sec:ex2}
We next show the results of the optimal solution under disturbances with the following configuration:
 $$P = 1,\,\, \gamma = 0.5,\,\, w_1(t)=0.1\sin(40t),\,\, w_2(t)=0.1\sin(2t).$$ 
 We omit reporting the behavior of the LQR-based control as it fails to stabilize the dynamics. To compare the proposed approaches we compute the $\Htwo-$ and $\Hinf-$ controls with the same disturbance using both Algorithm \ref{Alg_sdre} and Algorithm \ref{Alg_sdreon}. The results are presented in Figure \ref{fig1:dist1} and Figure \ref{fig1:dist2}. In every test case, the SDRE-based methodologies effectively stabilize the perturbed dynamics to a small neighbourhood around $X \equiv 0$.

\begin{figure}[htbp]
\centering
\includegraphics[scale=0.29]{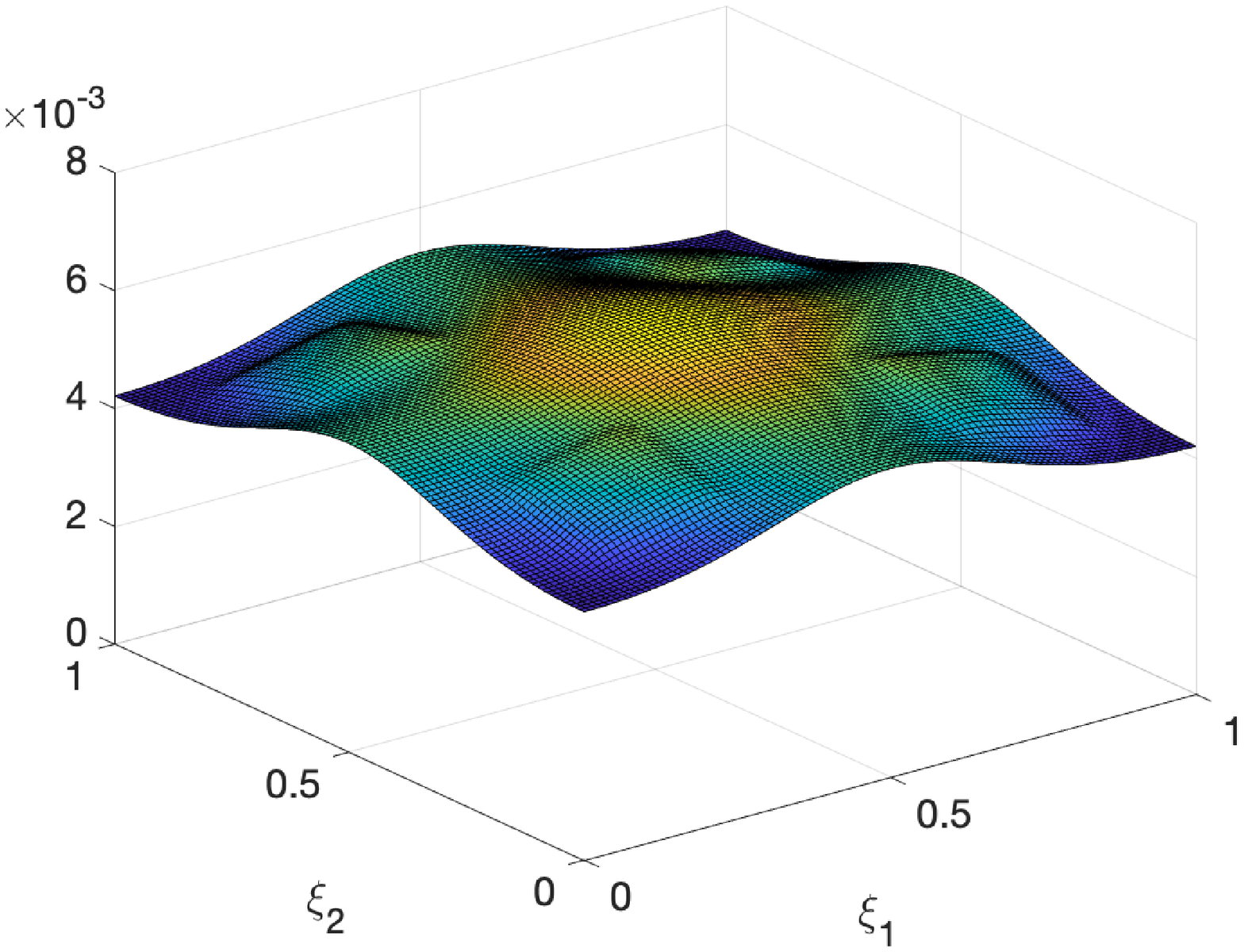}
\includegraphics[scale=0.29]{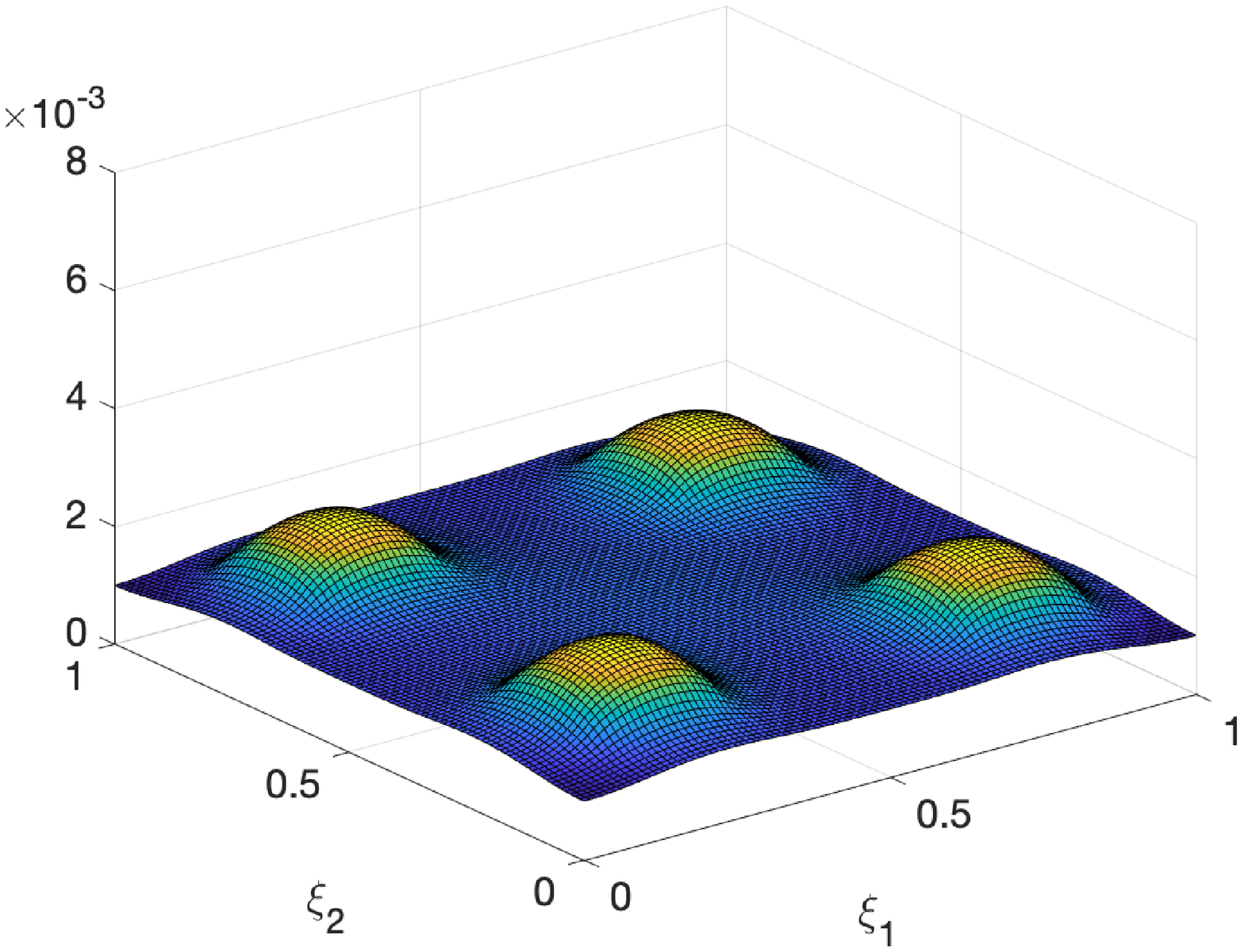}
\includegraphics[scale=0.29]{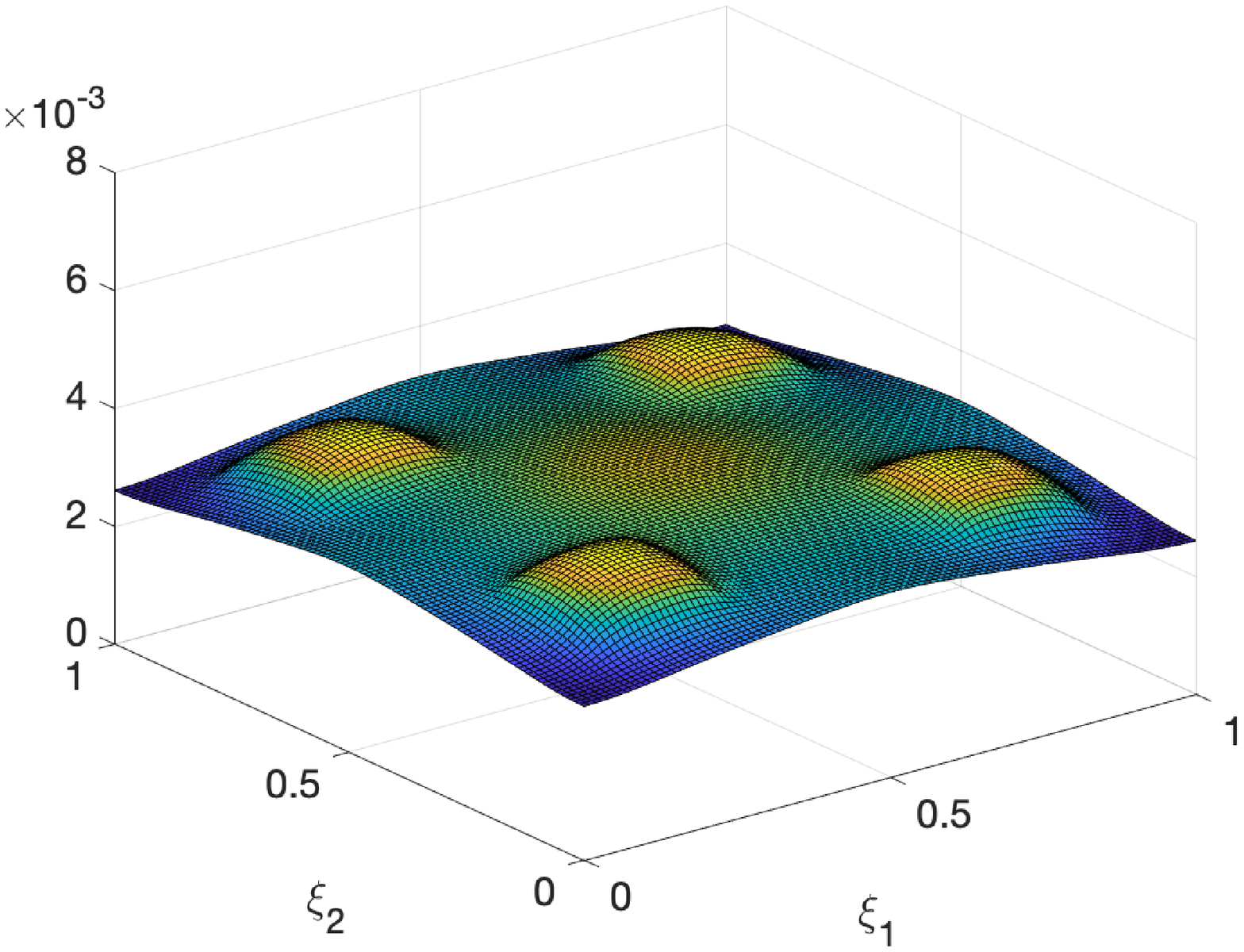}
\includegraphics[scale=0.29]{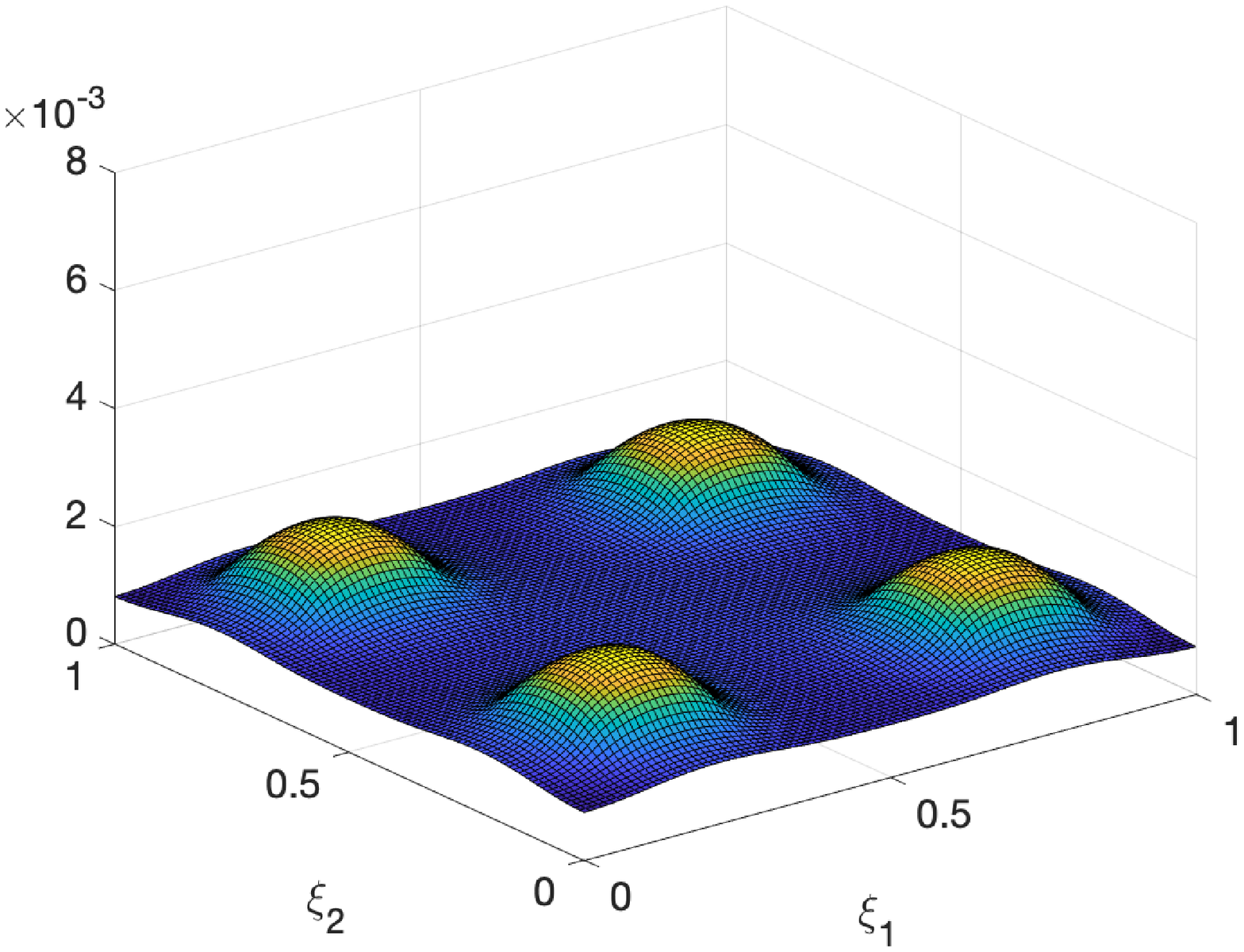}
\caption{Test \ref{sec:ex2}: Solutions at time $t=3$ with $\Htwo-$ and $\Hinf-$control and disturbance $w(t)=0.1\sin(40t)$.
Top:  $\Htwo$-control using Algorithm \ref{Alg_sdre} (left) and Algorithm \ref{Alg_sdreon} (right). Bottom:  $\Hinf$-control using Algorithm \ref{Alg_sdre} (left) and Algorithm \ref{Alg_sdreon} (right).}  
\label{fig1:dist2}
\end{figure}

A quantitative study is proposed in Figure \ref{fig1:an1} where we show the evaluation of the cumulative $\Htwo$ cost functional for both disturbances in the left panels. As expected, Algorithm \ref{Alg_sdre} with $\Hinf-$control exhibits the best performance, closely followed by  Algorithm \ref{Alg_sdreon}. The right panels of Figure \ref{fig1:an1} show the different control inputs, which reflect the observed differences in closed-loop performance.

\begin{figure}[htbp]
\centering
\includegraphics[scale=0.29]{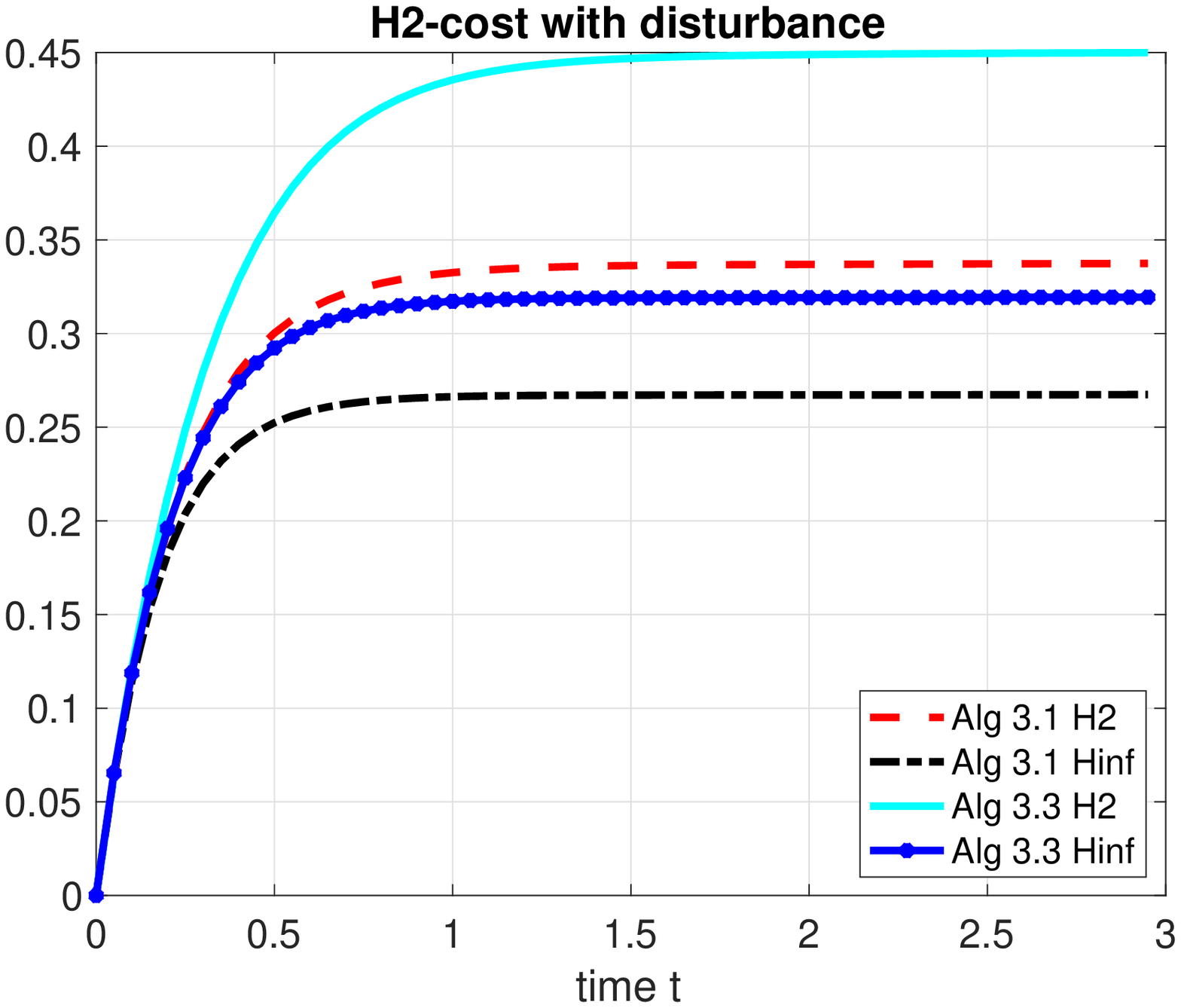}
\includegraphics[scale=0.29]{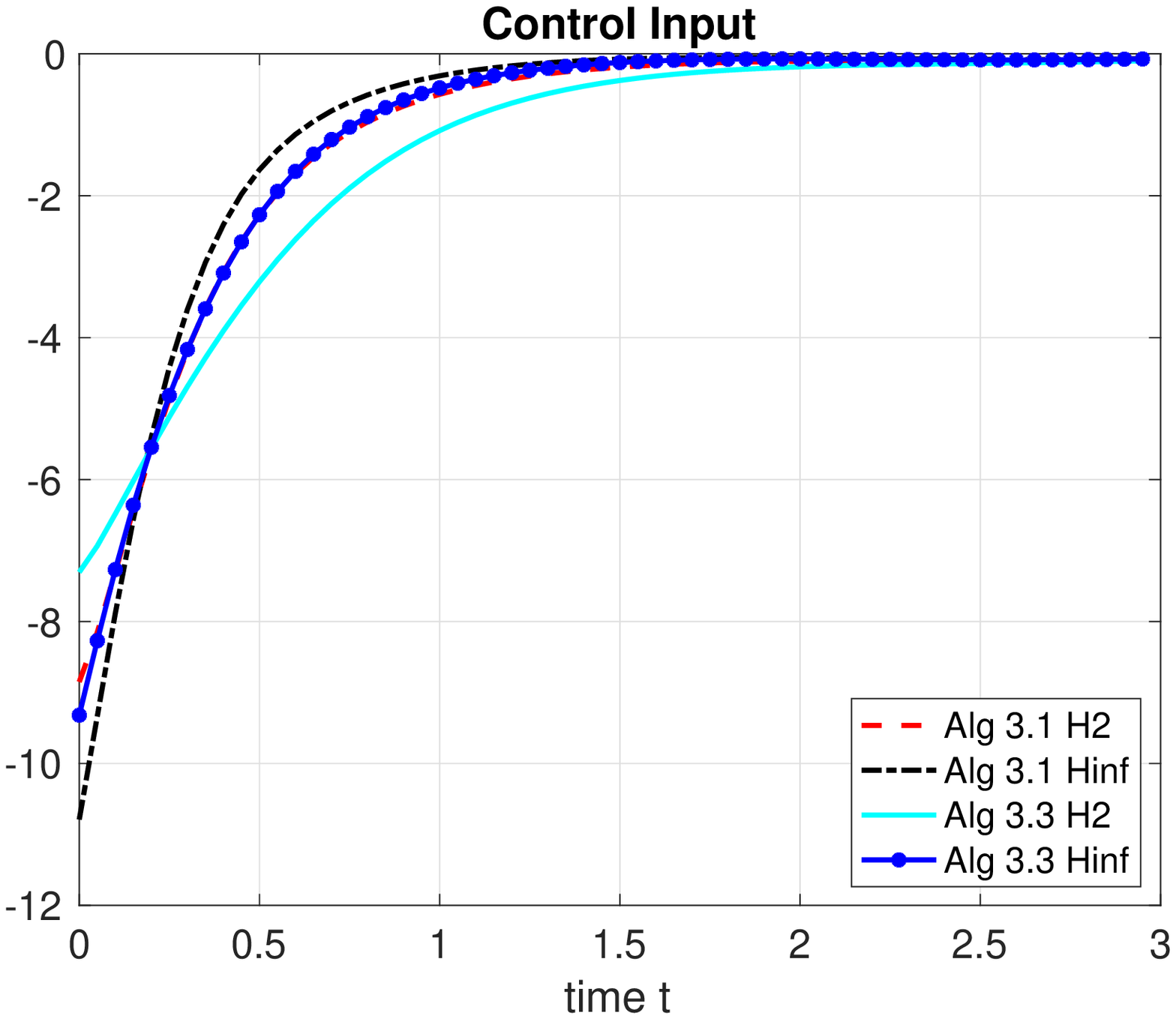}
\includegraphics[scale=0.29]{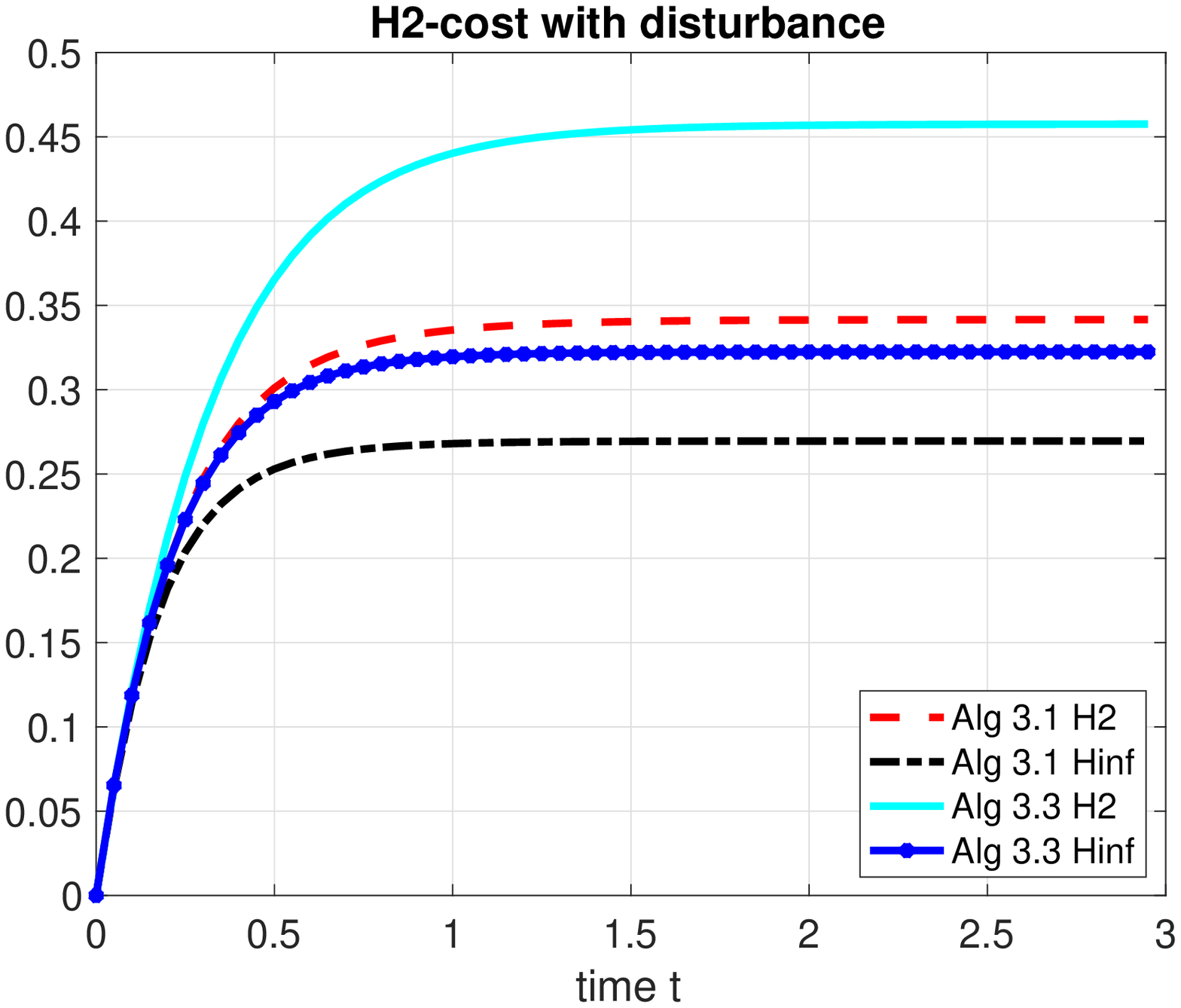}
\includegraphics[scale=0.29]{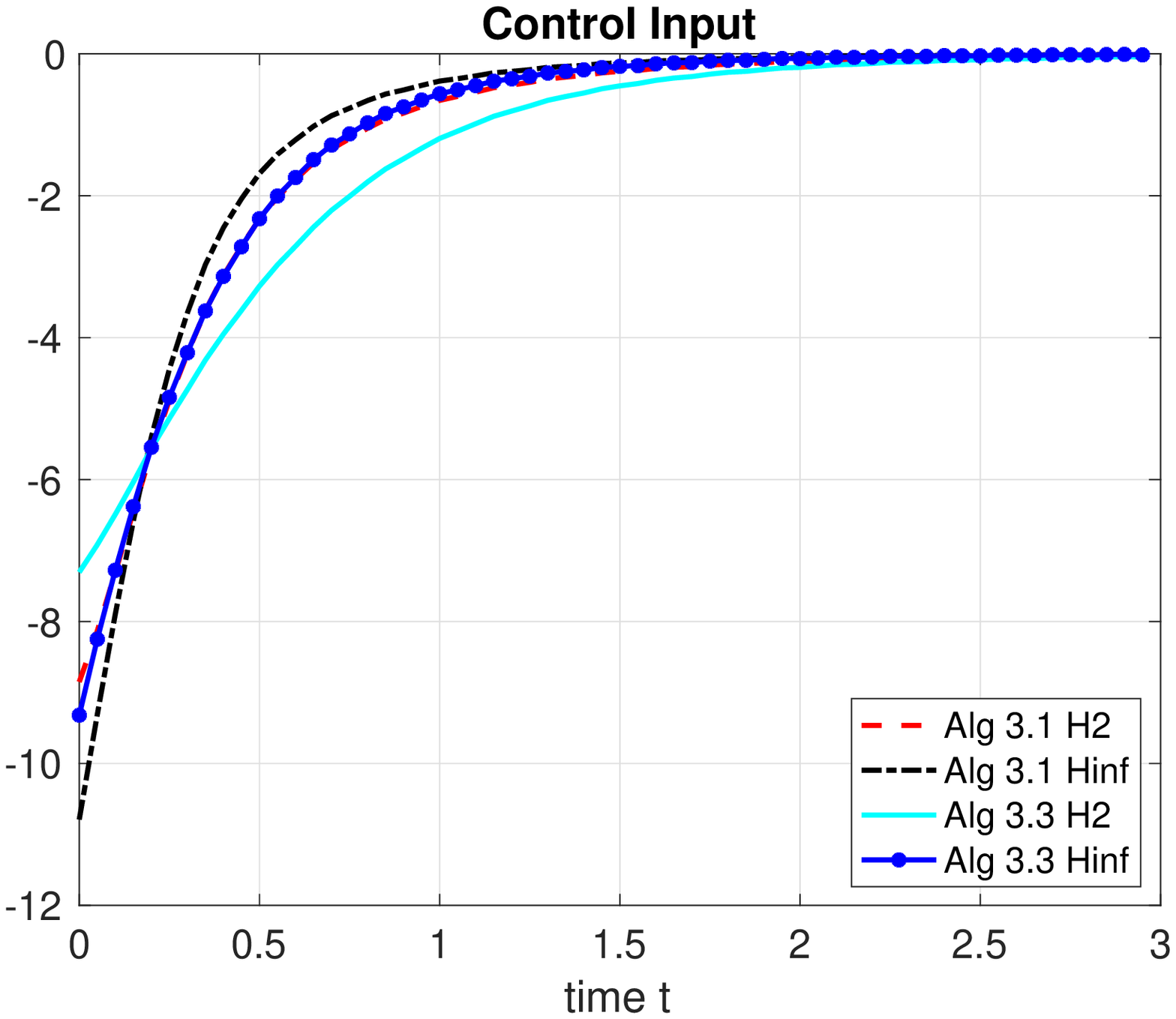}
\caption{Test \ref{sec:ex2}: evaluation of the cumulative cost functional $\Htwo$ (left) and control inputs (right) for $w(t)=0.1\sin(2t)$ (top) and $w(t)=0.1\sin(40t)$ (bottom). }
\label{fig1:an1}
\end{figure}
\end{example}
\begin{example}{\it On the use of a feedback matrix oriented implementation}\label{sec:ex3}

To conclude the first case study we provide a numerical example where we  synthesize the feedback operator $K(x)=-R^{-1}B^{\top}\Pi(x)$ directly, circumventing the computation of $\Pi(x)$, as explained in Section \ref{sec:feedback_matrix}. For this test we introduce the following changes: $\Omega = [0,1]\times[0,1],\,\, \epsilon=0.1,\,\, \nu = 0,\,\, \mu=8$,
and initial condition ${\mathbf x}_0(\xi) = \sin (\xi_1)\sin(\xi_2)$, on a discretized space grid of $n_{\xi_1}\times n_{\xi_2}$ nodes with $n_{\xi_1} = n_{\xi_2} = 101$. We replace Neumann boundary conditions with zero Dirichlet boundary conditions. For this test, we only consider the $\Htwo-$control case since many considerations are similar to the previous part of this section. {We use Algorithm \ref{Alg_sdre} with the Extended Krylov subspace as in Section \ref{sec:feedback_matrix} and compare the performances of Algorithm \ref{Alg_sdre} using the Rational Krylov subspaces and Extended Krylov subspaces.}

Figure \ref{fig3:sol} reports the solution at time $t=3$. In the top-left panel the uncontrolled state grows in time. In the top-right panel we show the results of the LQR control computed by linearizing equation \eqref{pol:sd} around the desired configuration. The control steers the solution to the origin at a very slow rate. 
The controlled solution with Algorithm \ref{Alg_sdre} {using an Extended Krylov subspace} and Algorithm \ref{Alg_sdreon} are shown at the bottom of the same figure. Both algorithms reach the desired configuration.

\begin{figure}[htbp]
\centering
\includegraphics[scale=0.29]{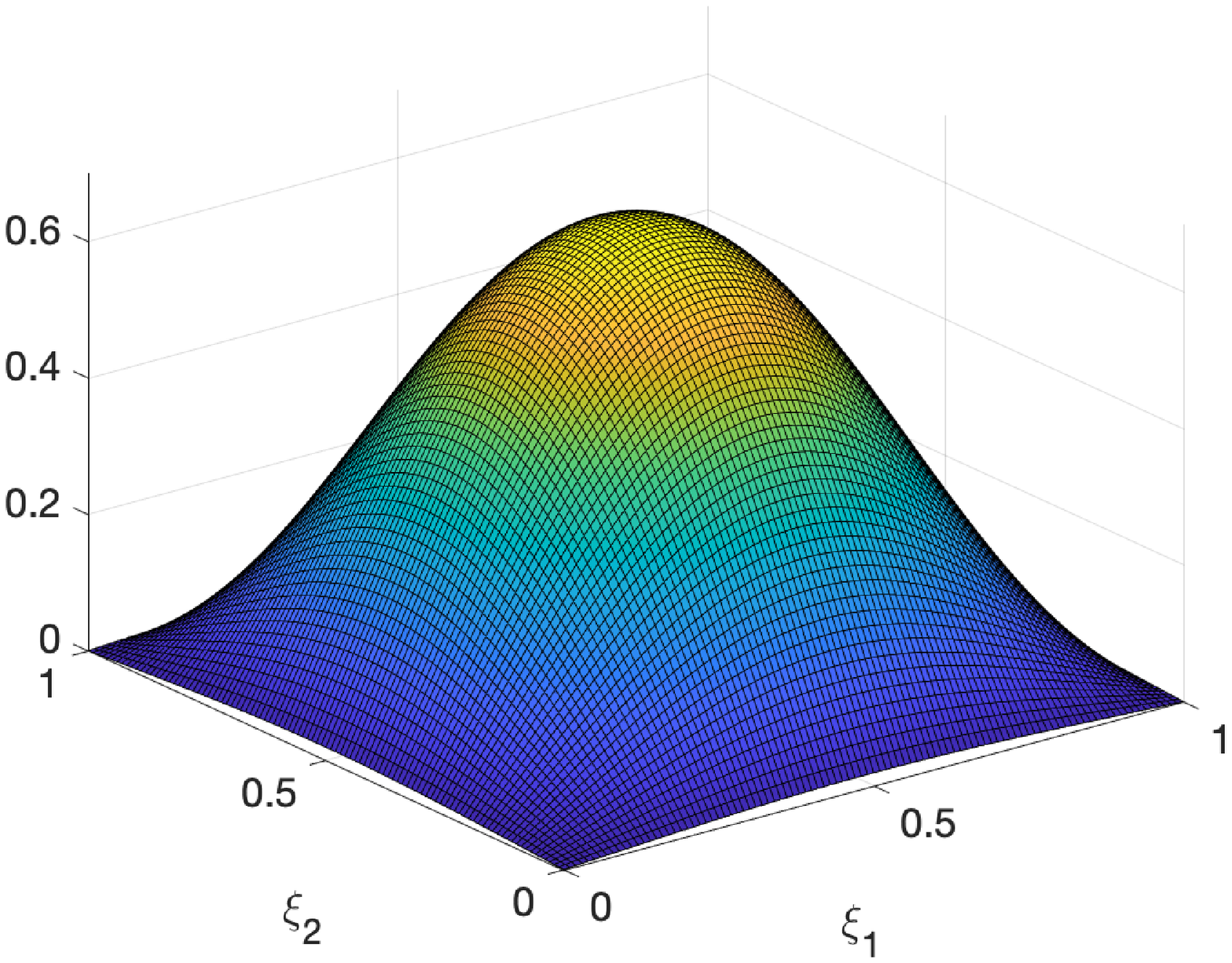}
\includegraphics[scale=0.29]{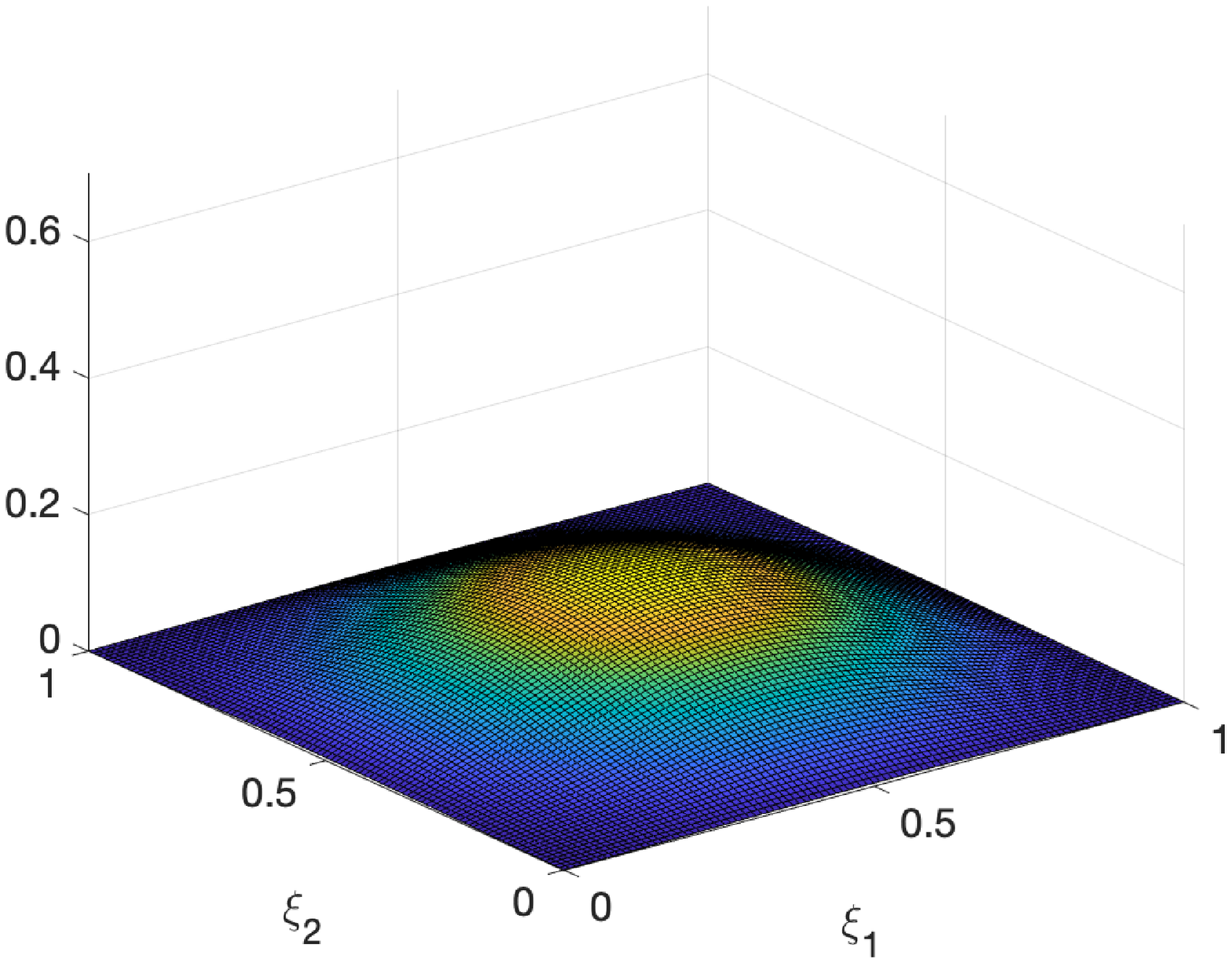}
\includegraphics[scale=0.29]{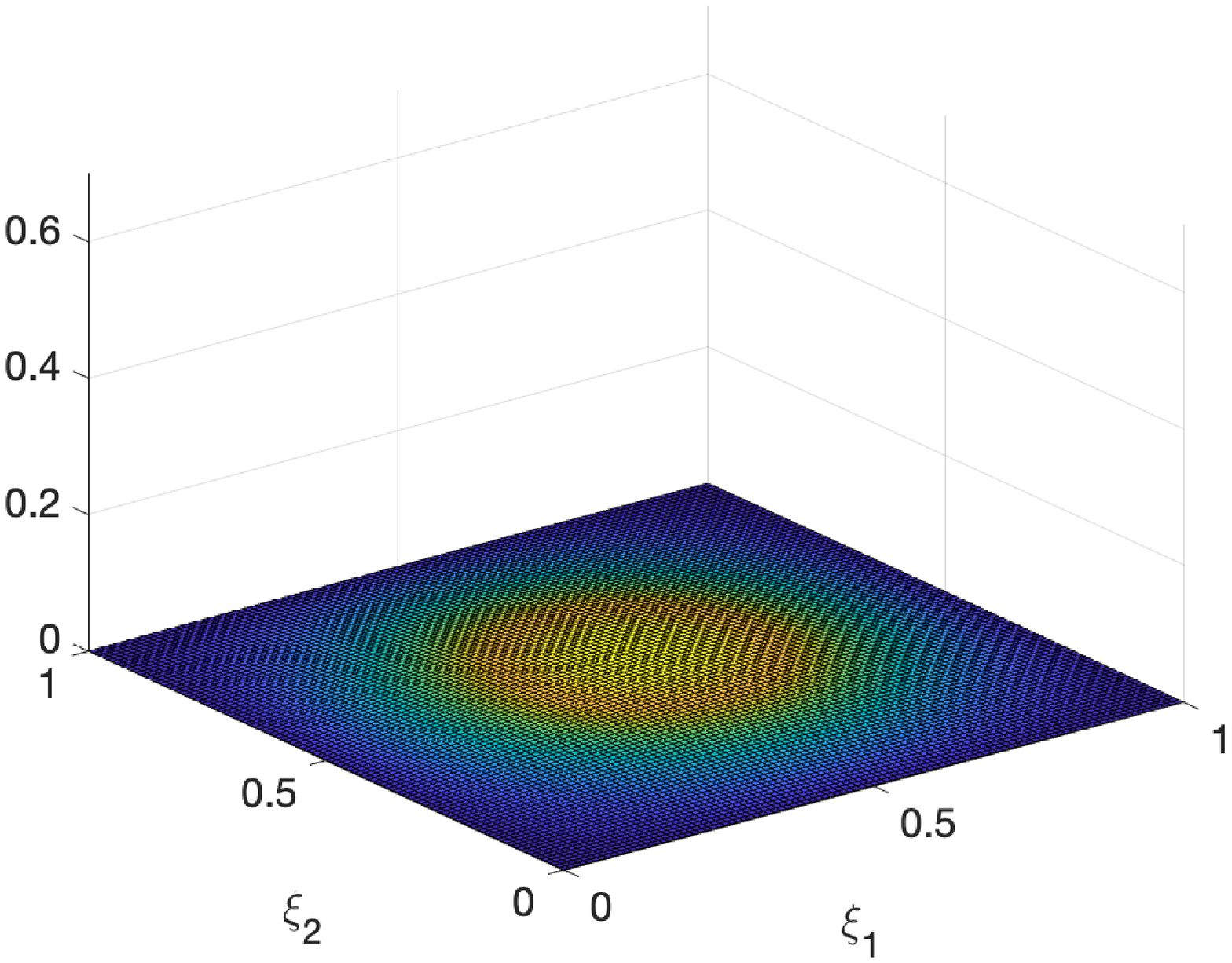}
\includegraphics[scale=0.29]{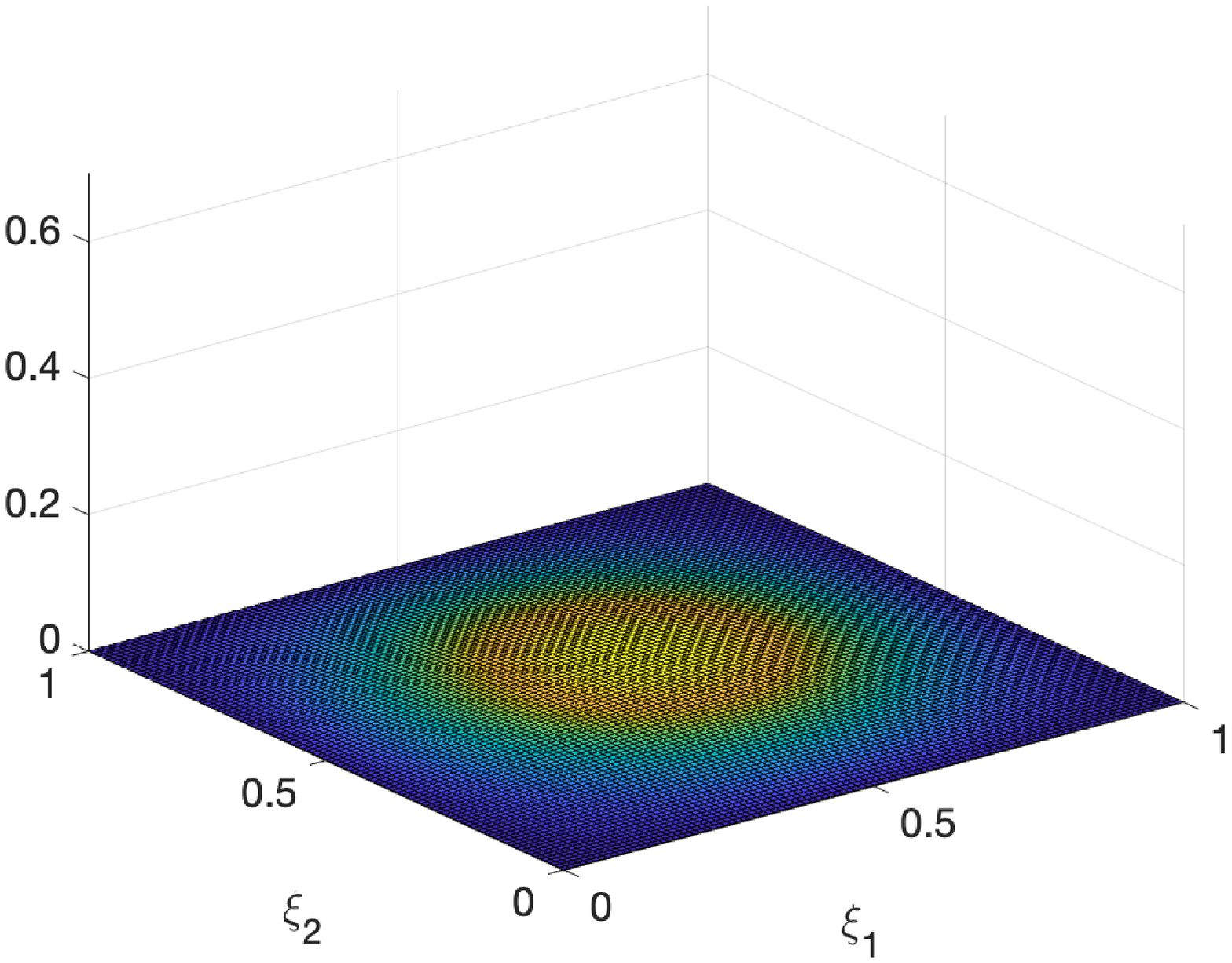}
\caption{Test \ref{sec:ex3}: controlled trajectories
\eqref{ex1} at time $t=3$ with zero Dirichlet boundary conditions. Top: Uncontrolled solution (left), LQR control (right). Bottom: $\mathcal{H}_2$-controlled solution with Algorithm \ref{Alg_sdre} with Extended Krylov subspaces (left), $\mathcal{H}_2$-controlled solution with Algorithm \ref{Alg_sdreon} (right).}
 \label{fig3:sol}
\end{figure}

\begin{figure}[htbp]
\centering
\includegraphics[scale=0.3]{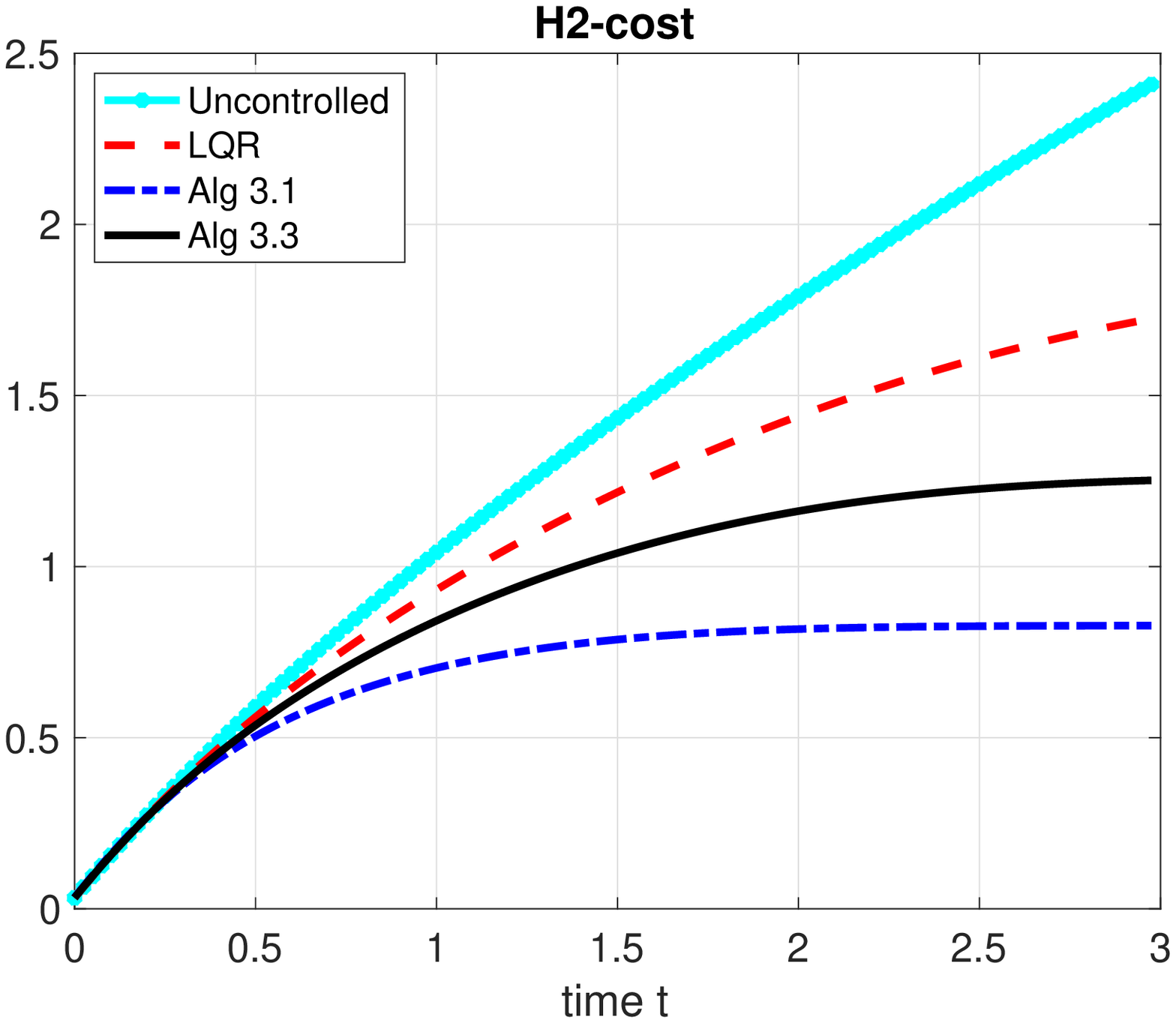}
\includegraphics[scale=0.3]{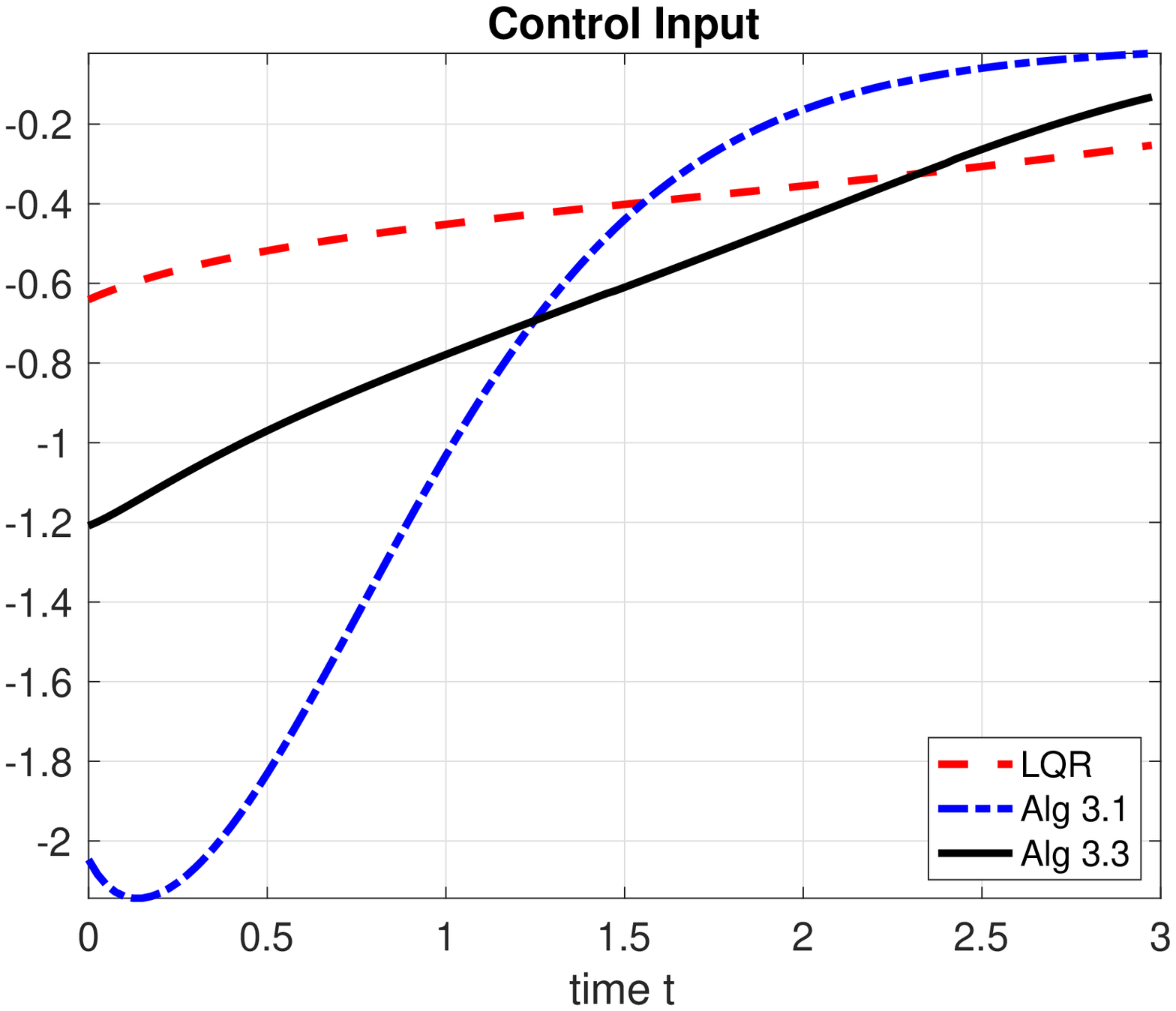}
\includegraphics[scale=0.3]{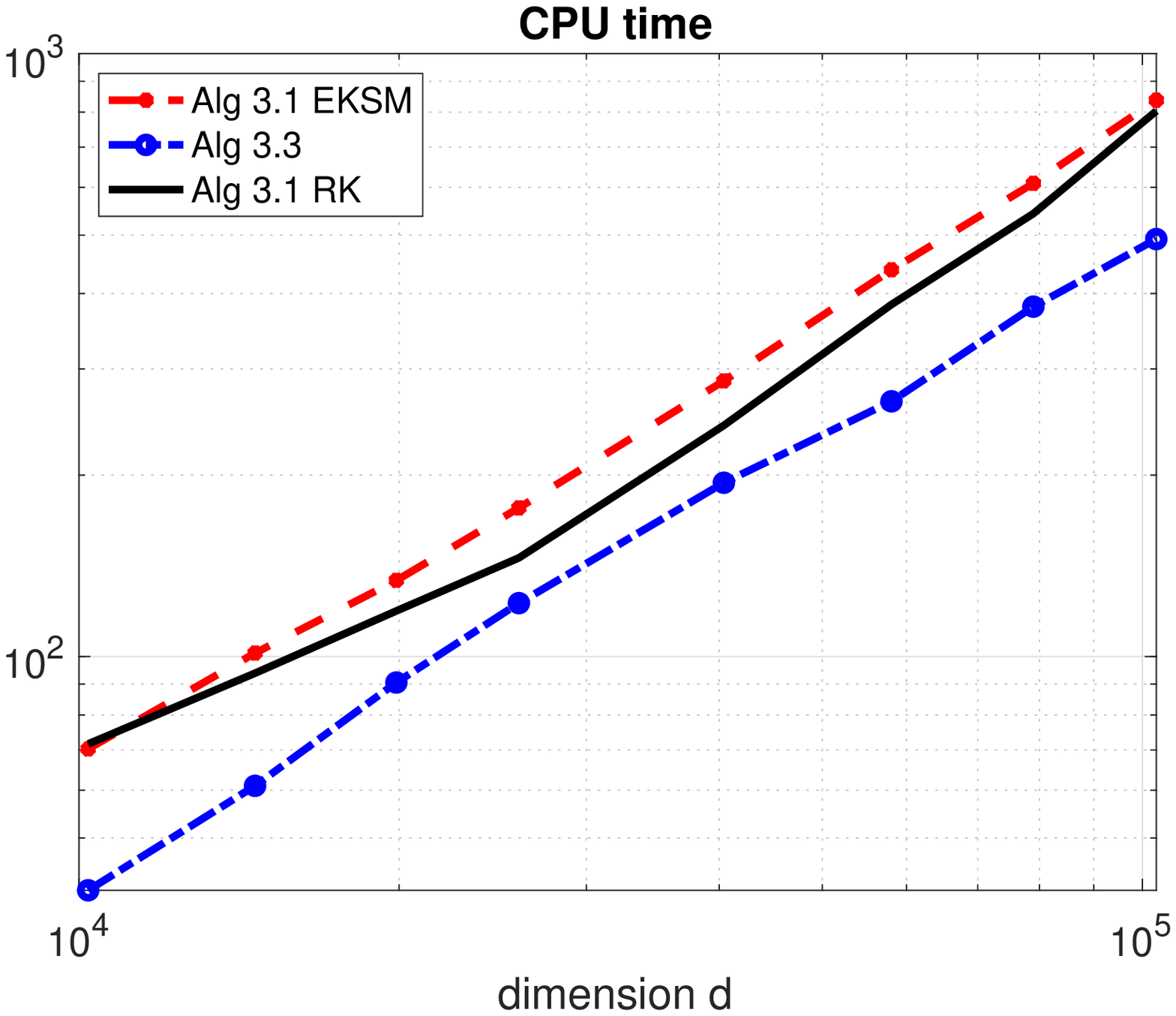}
\includegraphics[scale=0.3]{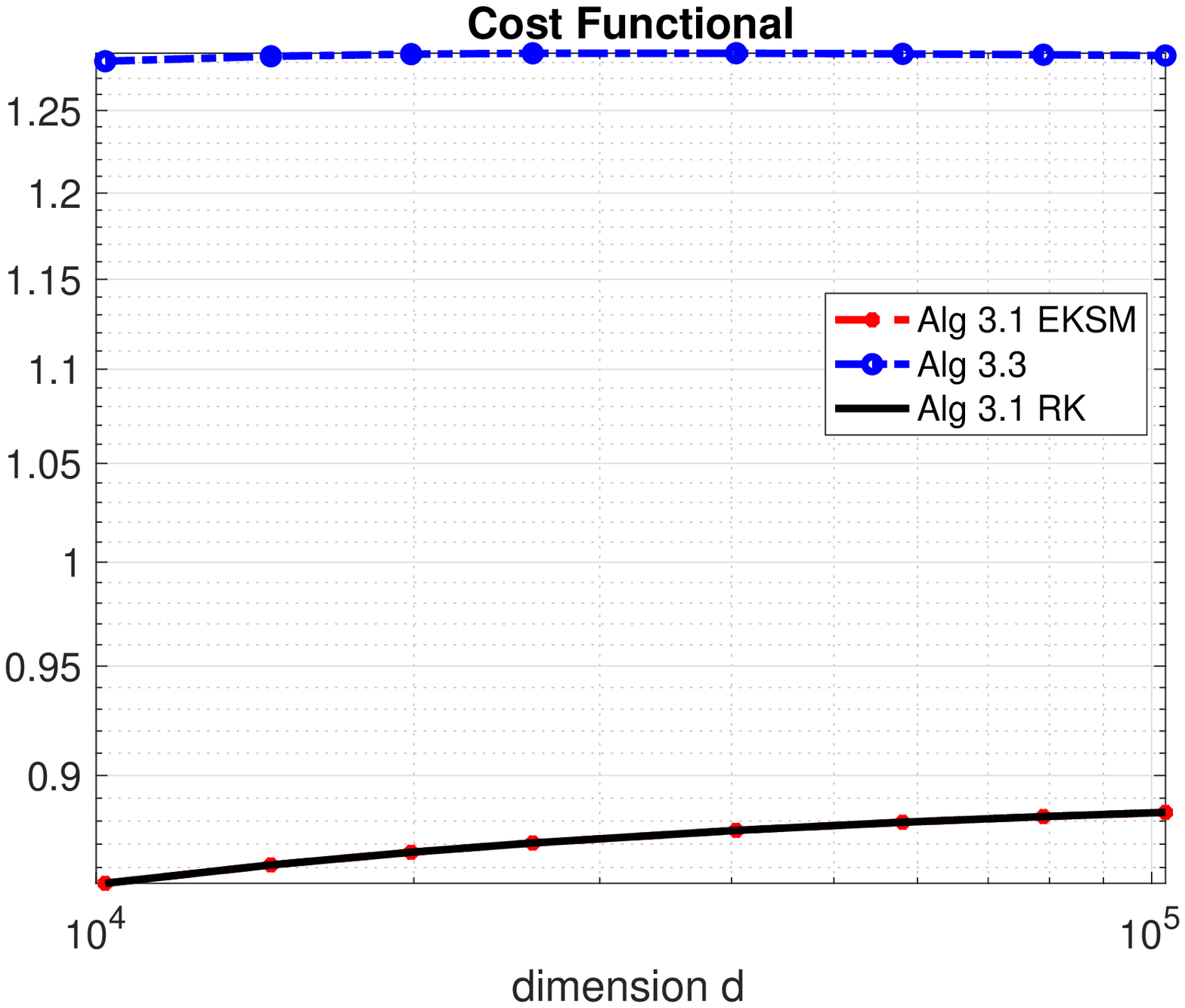}
\caption{Test \ref{sec:ex3}: cumulative functional (top-left) with $\Htwo-$control and corresponding control input (top-right). Evaluation of the cost functional with respect to the dimension of the problem $d$ (x-axis) (bottom-right) and CPU time (bottom-left) for both Algorithm \ref{Alg_sdre} {using Rational Krylov subspaces and Extended Krylov subspaces} and Algorithm \ref{Alg_sdreon}.}
\label{fig3:cost}
\end{figure}
\end{example}

The performance of the different feedback synthesis methods is presented in Figure \ref{fig3:cost}. Here we show the evaluation of the cumulative cost functional in the top-left panel. For completeness we provide the control inputs on the top-right panel. We compare the performance of Algorithm \ref{Alg_sdre} using Rational Krylov (RK) and Extended Krylov (EKSM) subspaces as the problem dimension increases. Specifically, we recall that the Extended Krylov subspace computes directly the gain matrix. As expected, Algorithm \ref{Alg_sdre} provides the best closed-loop performance among the proposed algorithms. However, in terms of CPU time Algorithm~\ref{Alg_sdreon} is faster than Algorithm~\ref{Alg_sdre} when increasing 
the problem dimension $d$ as shown in the bottom-left panel of Figure~\ref{fig3:cost}. 
When the problem dimension $d$ increases, the cost functional converges to $0.9$ for Algorithm~\ref{Alg_sdre} and to 
$1.3$ for Algorithm \ref{Alg_sdreon}. 
The two different Krylov subspace methods in Algorithm \ref{Alg_sdre} lead exactly to the same solution. 
{In the plot legend we refer to {\it Alg 3.1 RK} for Rational Krylov subspaces and to 
{\it Alg 3.1 EKSM} for Extended Krylov subspaces.}

\subsection{Case study 2: the viscous Burgers equation with exponential forcing term}

The second experiment deals with the control of a viscous Burgers equation with exponential forcing term over $\Omega\times\R^+_0$, with $\Omega\subset\R^2$ and Dirichlet boundary conditions:
 \begin{align}\label{ex2}
\begin{aligned}
\pt X(\xi,t)&=\epsilon\Delta X(\xi,t)- X(\xi,t)\cdot \nabla X(\xi,t)+1.5X(\xi,t)e^{-0.1X(\xi,t)}\\
&\qquad+\chi_{\omega_c}(\xi)u(t)+\chi_{\omega_d}(\xi)w(t)\\
X(\xi,t)&=0\,, \quad\xi\in\partial\Omega,\\
 X(\xi,0)&={\mathbf x}_0(\xi)\,,\quad\xi\in\Omega\,.
\end{aligned}
\end{align}
In this case, the scalar control and disturbance act, respectively, through the indicator function $\chi_{\omega_c}(\xi), \chi_{\omega_d}(\xi)$ with $\omega_c,\omega_d\subset\Omega$. A finite difference discretization of the space of the system dynamics leads to a state space representation of the form
\begin{equation}\label{burg:sd}
\dot X(t) = \Delta_d X(t)- X(t)\circ (DX(t) + 1.5e^{-0.1X(t)}) +B u(t) + H w(t)\,,
\end{equation}
where the matrices $\Delta_d,D\in\R^{d\times d}$ and $B,H\in \R^d$ are finite-dimensional approximations of the Laplacian, gradient, control and disturbance operators, respectively, and the exponential term is understood component-wise. {In particular, $D$ is obtained using a backward finite difference discretization
$$D:=-\Delta\xi^{-1}( B_{n_{\xi_2}}\otimes \mathbf{I}_{n_{\xi_1}} + \mathbf{I}_{n_{\xi_2}}\otimes B_{n_{\xi_1}})\,,$$
where  $B_n:=\texttt{tridiag}([\begin{array}{ccc}-1 &1& 0\end{array}],n)$, and
$\otimes$ denotes the Kronecker product.}
We proceed to express the semi-discretized dynamics in semilinear form. For this, we define
$$A(X) := \epsilon\Delta_d -\tilde D(X) + \texttt{diag}(1.5e^{-0.1X(t)})\,,\quad [\tilde D(X)]_{k,l}=D_{k,l}X_k\,,$$
and then
\begin{align*}
A_0 = \epsilon\Delta_d,\quad [A_j]_{k,l}=D_{j,l}\delta_{k,l},\quad f_j(X)=X_j,\quad j=1,\ldots,d\,,\\
[A_j]_{k,l}=\delta_{k,j-d}\delta_{l,j-d},\quad f_j(X)=1.5e^{-0.1X_{j-d}},\quad j=d+1,\ldots,2d\,.
\end{align*}
For our numerical experiments we set $\Omega = [0,1]\times[0,1]$, $\epsilon=0.1$, $R = 0.05$,
and initial condition ${\mathbf x}_0(\xi) = \sin (\xi_1)\sin(\xi_2)$, on 
a discretized space grid of $n_{\xi_1}\times n_{\xi_2}$ nodes with $n_{\xi_1} = n_{\xi_2} = 101$ ($d=10201$). 
The matrices $B$ and $C$ are defined as in the previous study case (see Figure \ref{fig:bc_coll}). 
In the following we discuss the results for $\Htwo$ and $\Hinf$ control, i.e., $P = 0$ and $P\neq 0$, respectively, in \eqref{costex1}.  

\begin{example}{\it Experiments for $\Htwo-$control.}\label{sec2:ex1}

The trajectories of the controlled system problem with $P = 0$ in \eqref{costex1} are shown in Figure \ref{fig2:sol}. The uncontrolled solution tends to move towards the top-right corner of $\Omega$. All algorithms tend to control the solution to zero for large time, but at a different rate.

\begin{figure}[htbp]
\centering
\includegraphics[scale=0.29]{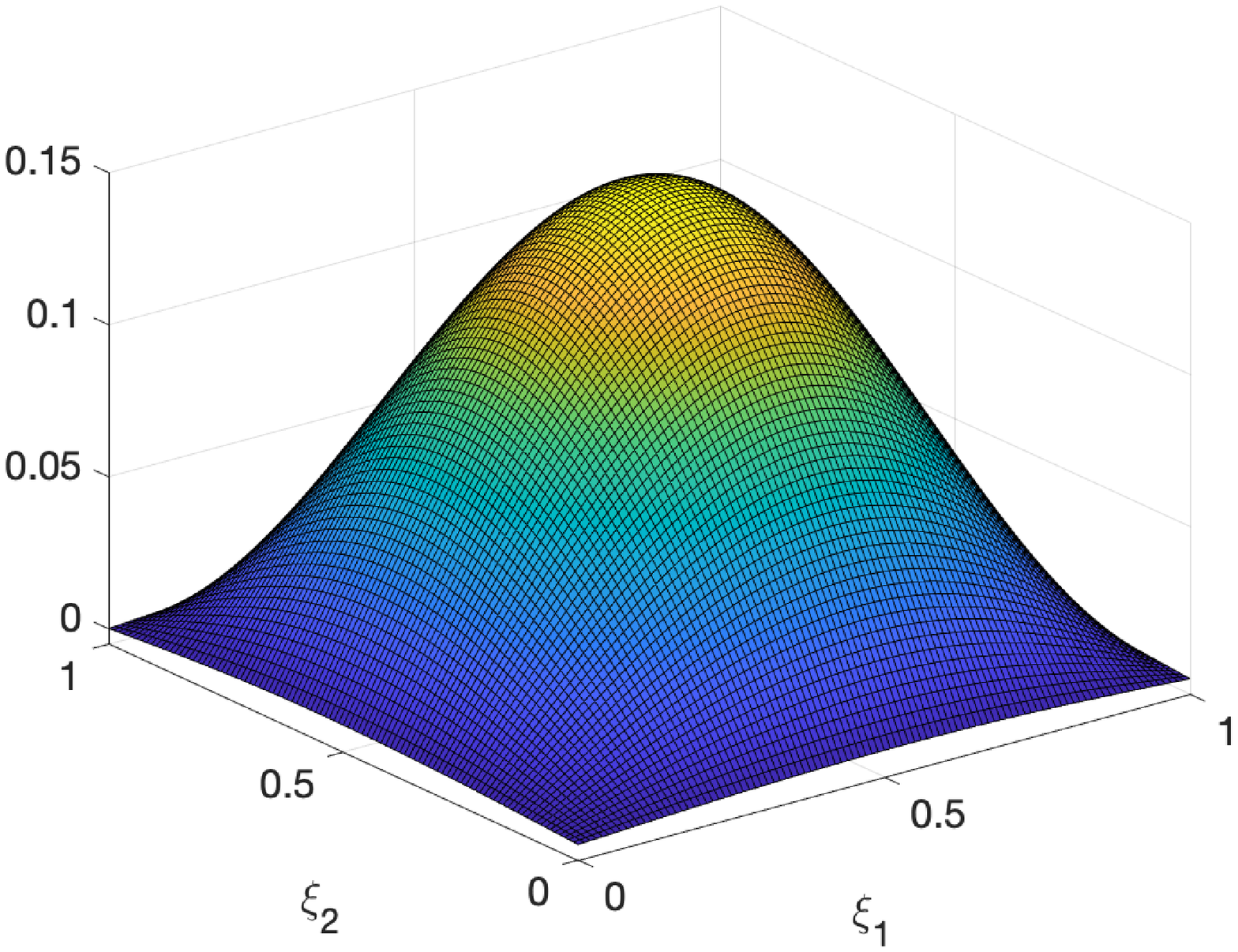}
\includegraphics[scale=0.29]{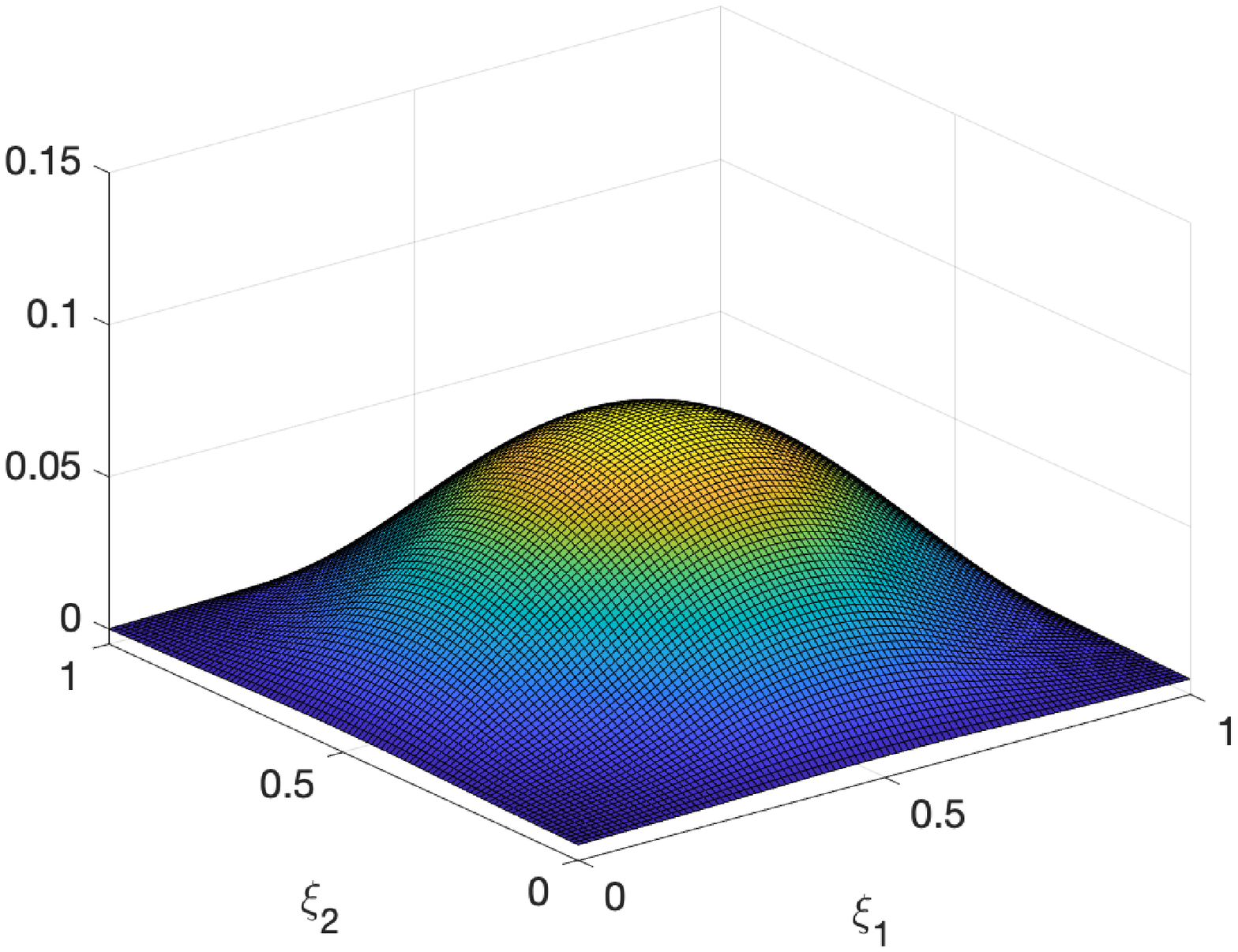}
\includegraphics[scale=0.29]{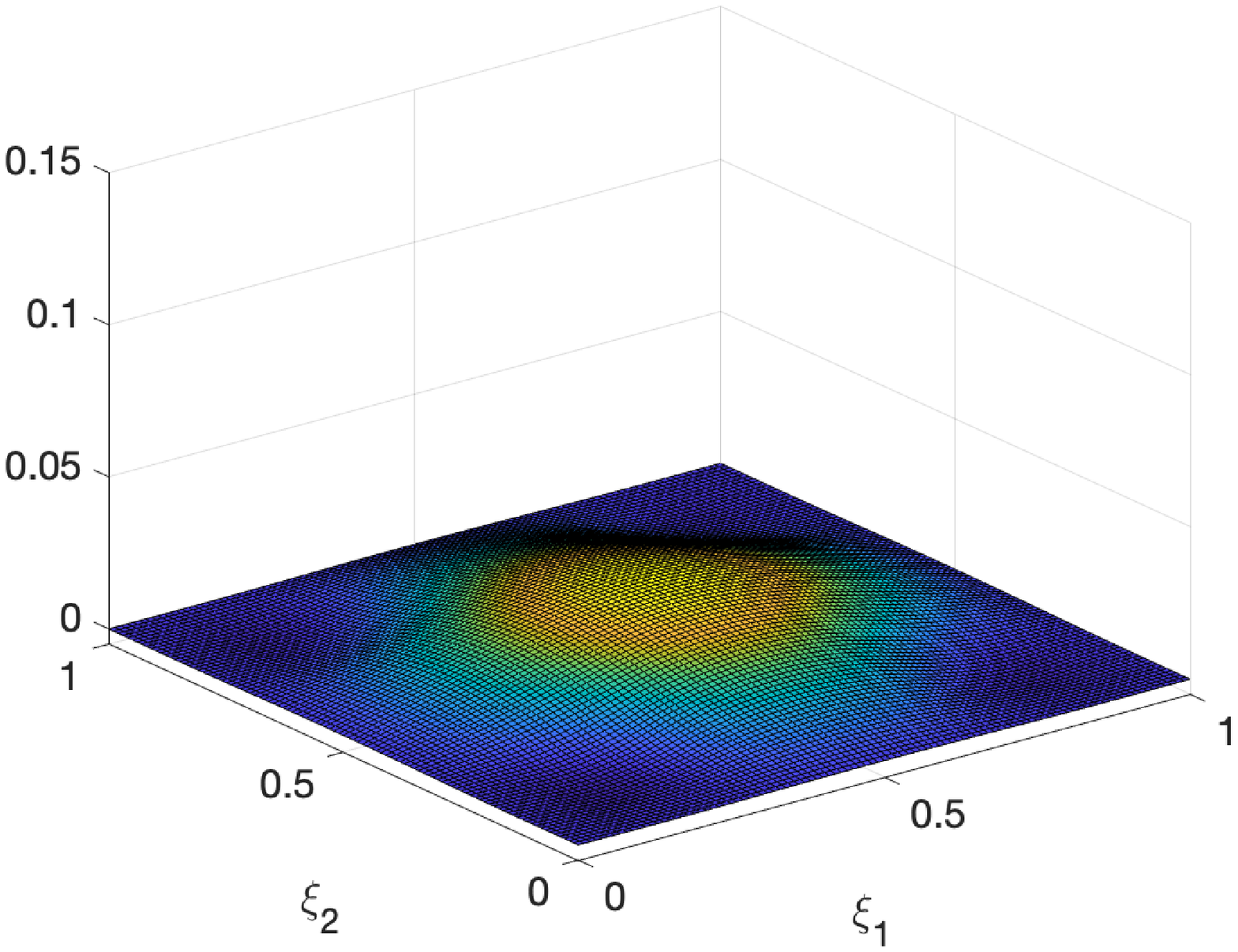}
\includegraphics[scale=0.29]{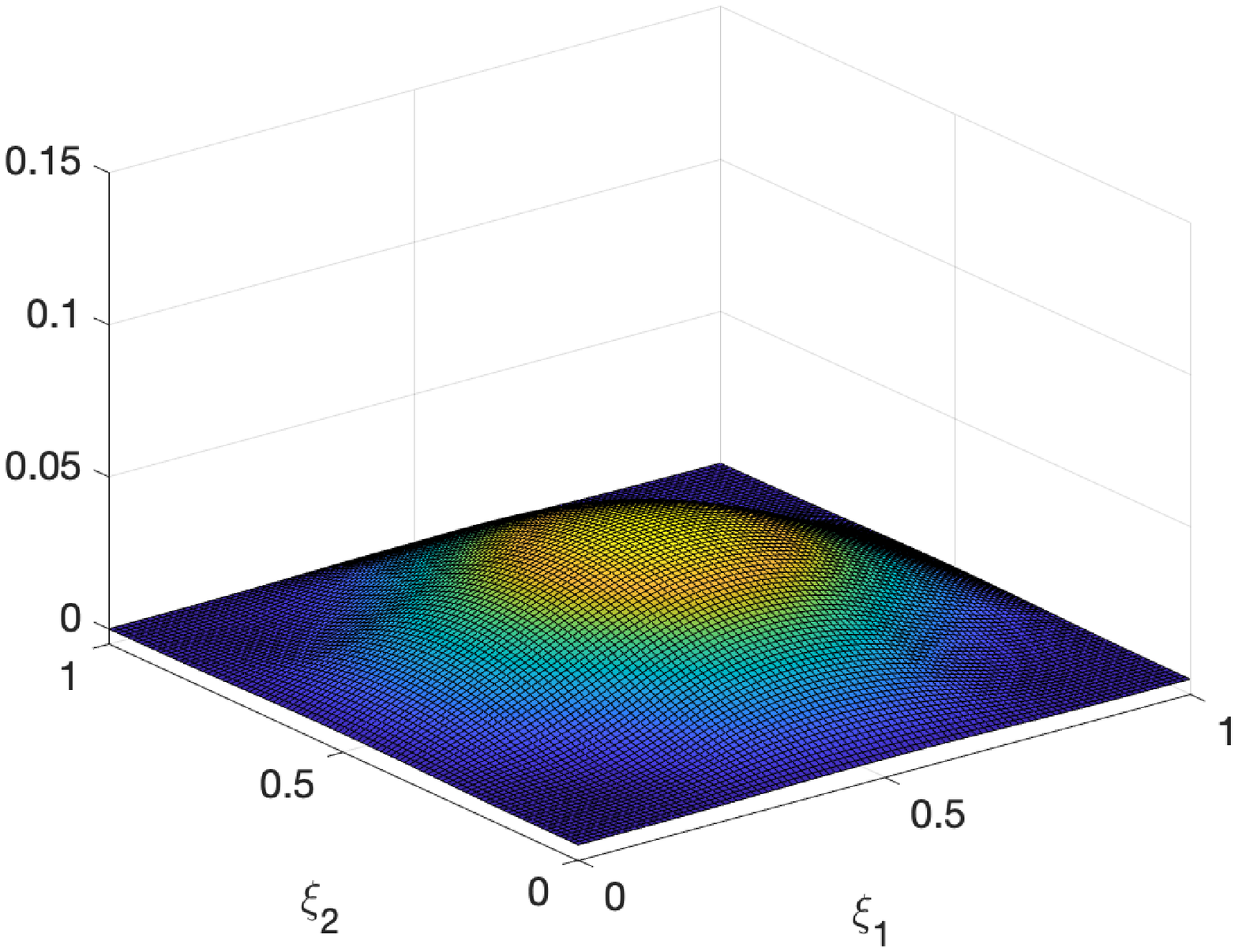}
\caption{Test \ref{sec2:ex1}: controlled dynamics
\eqref{ex2} at time $t=3$. Top: uncontrolled solution (left), $\mathcal{H}_2$-solution with LQR control (right). Bottom: $\mathcal{H}_2$-controlled solution with Algorithm \ref{Alg_sdre} (left), $\mathcal{H}_2$-controlled solution with Algorithm \ref{Alg_sdreon} (right). The controllers stabilize the dynamics to $X\equiv 0$ at a different rate.}
 \label{fig2:sol}
\end{figure}

 The control inputs are then shown in the top-right of Figure \ref{fig2:cost}. Algorithm \ref{Alg_sdre} has the largest control input, leading to the smallest cumulative cost functional as shown in the top-left panel of Figure \ref{fig2:cost}. 
We also observe that the values of the cost functional are very similar for Algorithm~\ref{Alg_sdre} 
and Algorithm~\ref{Alg_sdreon}. This is also confirmed for different discretization of increasing dimensions as shown in the bottom-left panel of Figure \ref{fig2:cost}. However, the time needed for Algorithm \ref{Alg_sdreon} to compute the solution is lower than for Algorithm~\ref{Alg_sdre} as depicted in the bottom-right panel of Figure \ref{fig2:cost}.

\begin{figure}[htbp]
\centering
\includegraphics[scale=0.29]{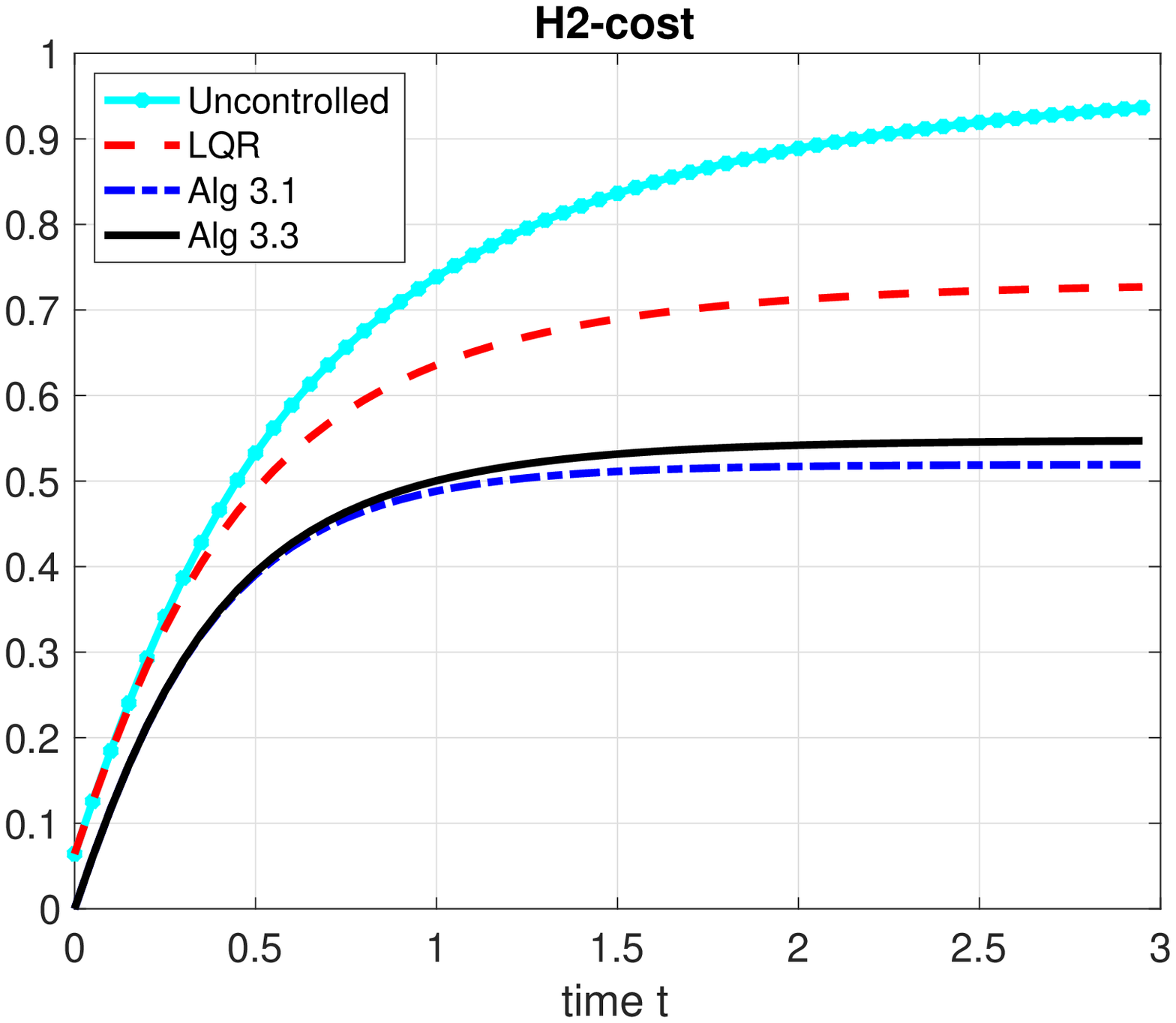}
\includegraphics[scale=0.29]{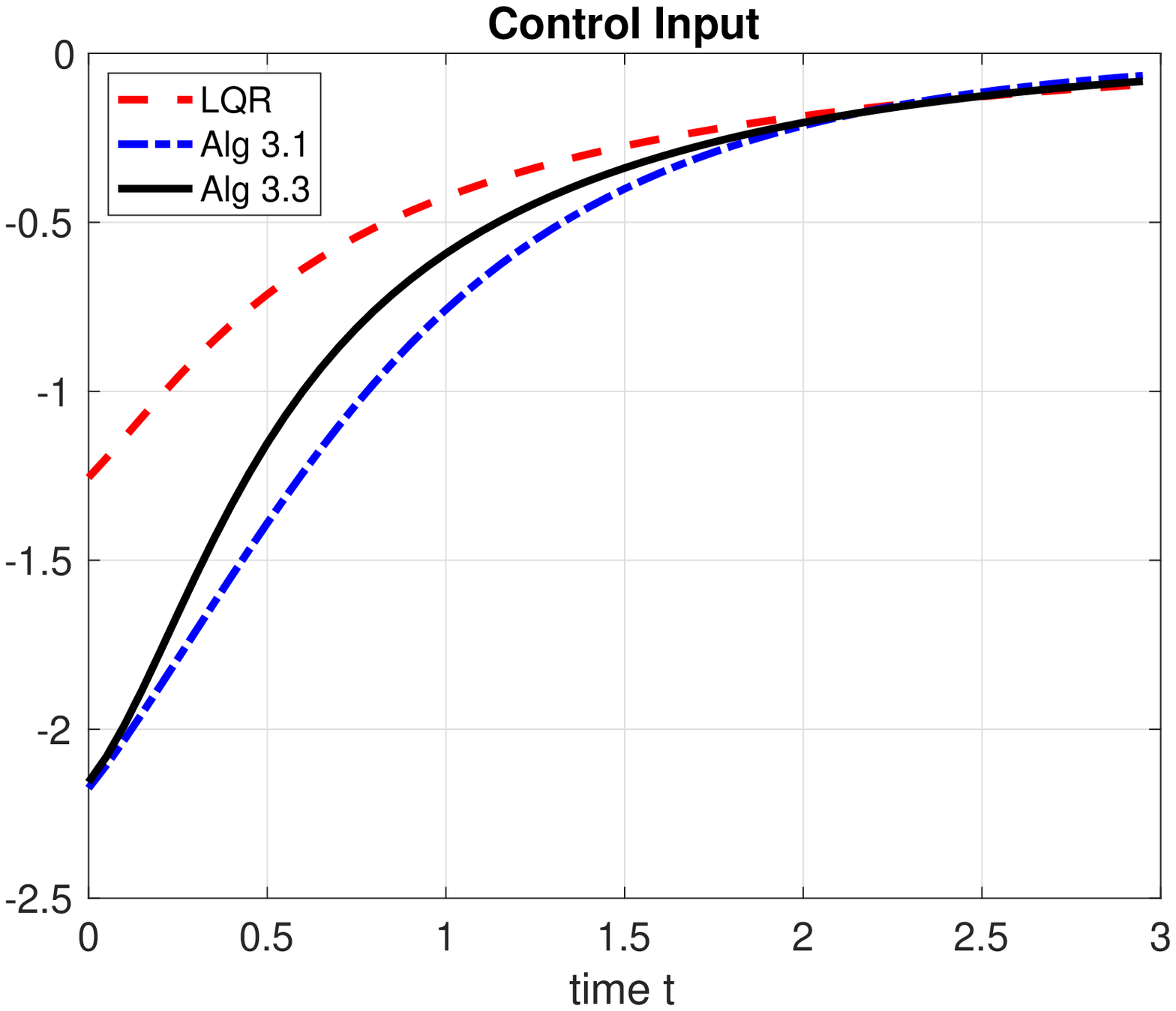}
\includegraphics[scale=0.29]{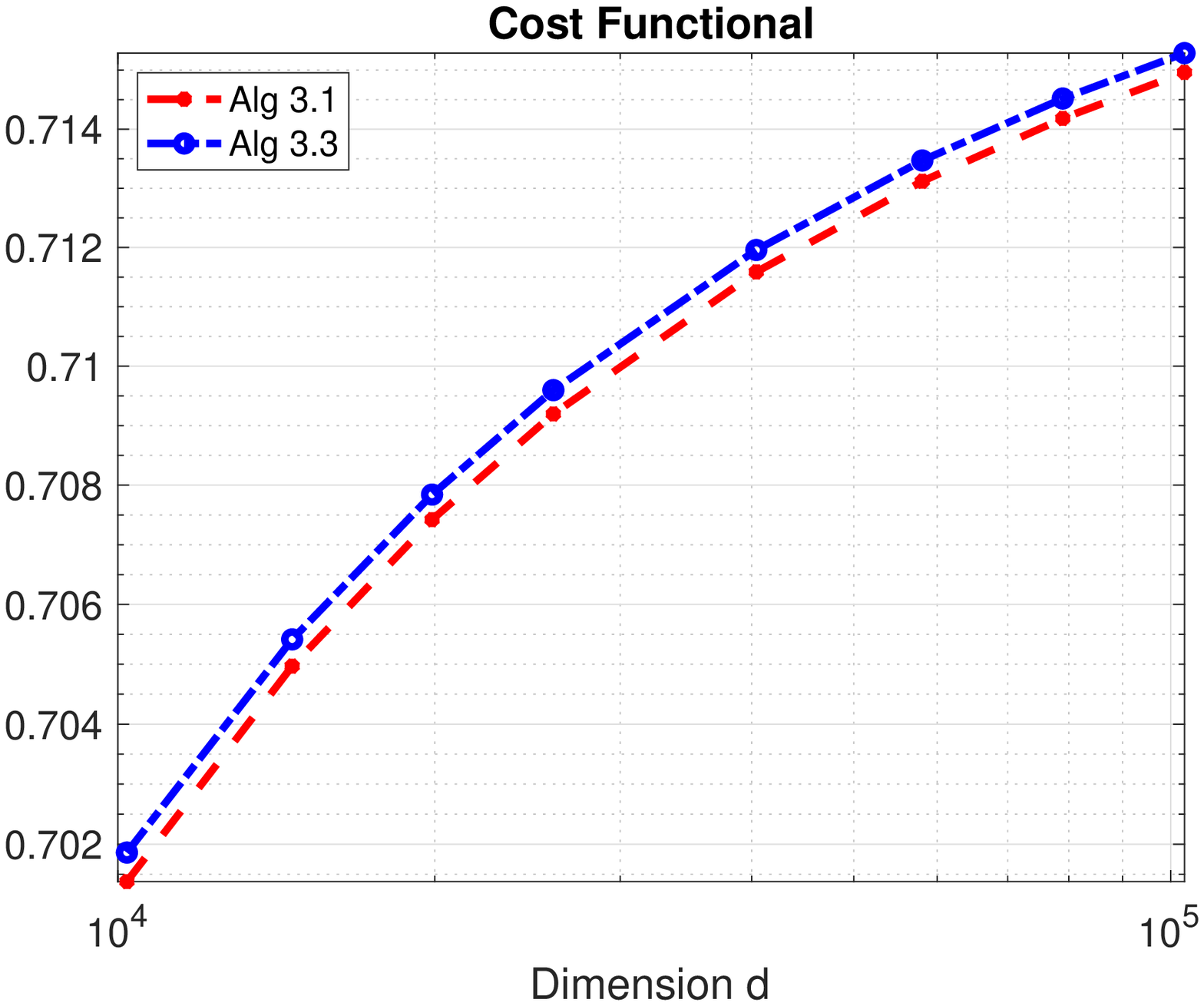}
\includegraphics[scale=0.29]{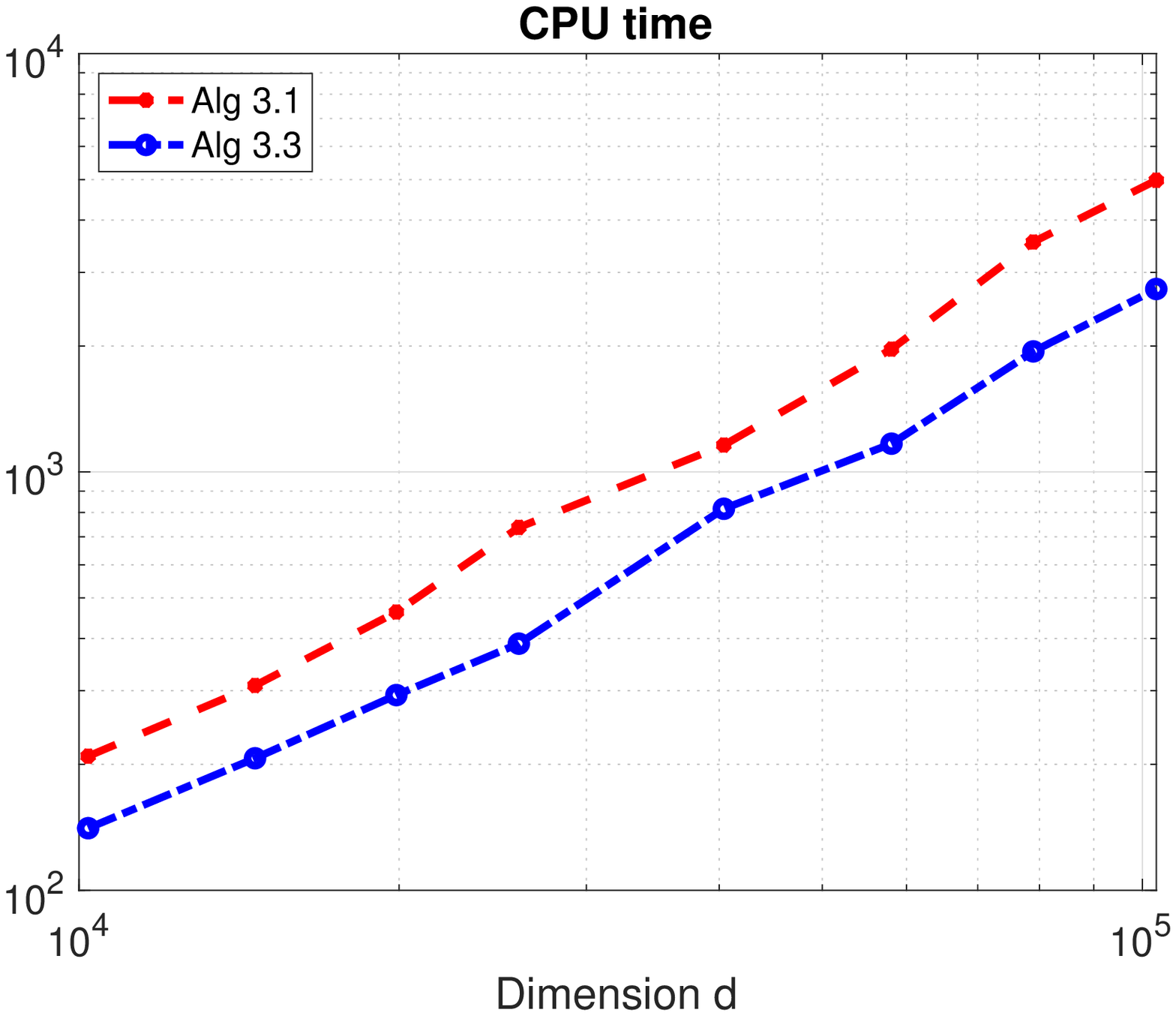}
\caption{Test \ref{sec2:ex1}. Top: evaluation of the cumulative cost functional with $\Htwo-$control (left) and control inputs (right). Bottom: cost functional for dynamics of increasing dimension for different algorithms (left) and CPU time for both Algorithm \ref{Alg_sdre} and Algorithm \ref{Alg_sdreon} (right).}
\label{fig2:cost}
\end{figure}
\end{example}

\begin{example}{\it Experiments for $\Hinf-$control.}\label{sec2:ex2}

Finally, we discuss the results for $P=1$, $\gamma=0.1$ and disturbance $w(t) = \{0.1\sin(2t)\}$ in \eqref{ex2}. The results presented in Figure \ref{fig2:an1} are in line with our first case study. In this example, Algorithm \ref{Alg_sdre} stabilizes the solution faster. Since it is difficult to visualize differences in the controlled state variables, we provide a qualitative analysis through Figure \ref{fig2:an1}. We show the evaluation of the cost functional \eqref{costex1} on the left panel and the control inputs on the right panel. Again, we find that Algorithm \ref{Alg_sdre} with $\gamma\neq0$ has the lowest values for the cost functional.

\begin{figure}[htbp]
\centering
\includegraphics[scale=0.29]{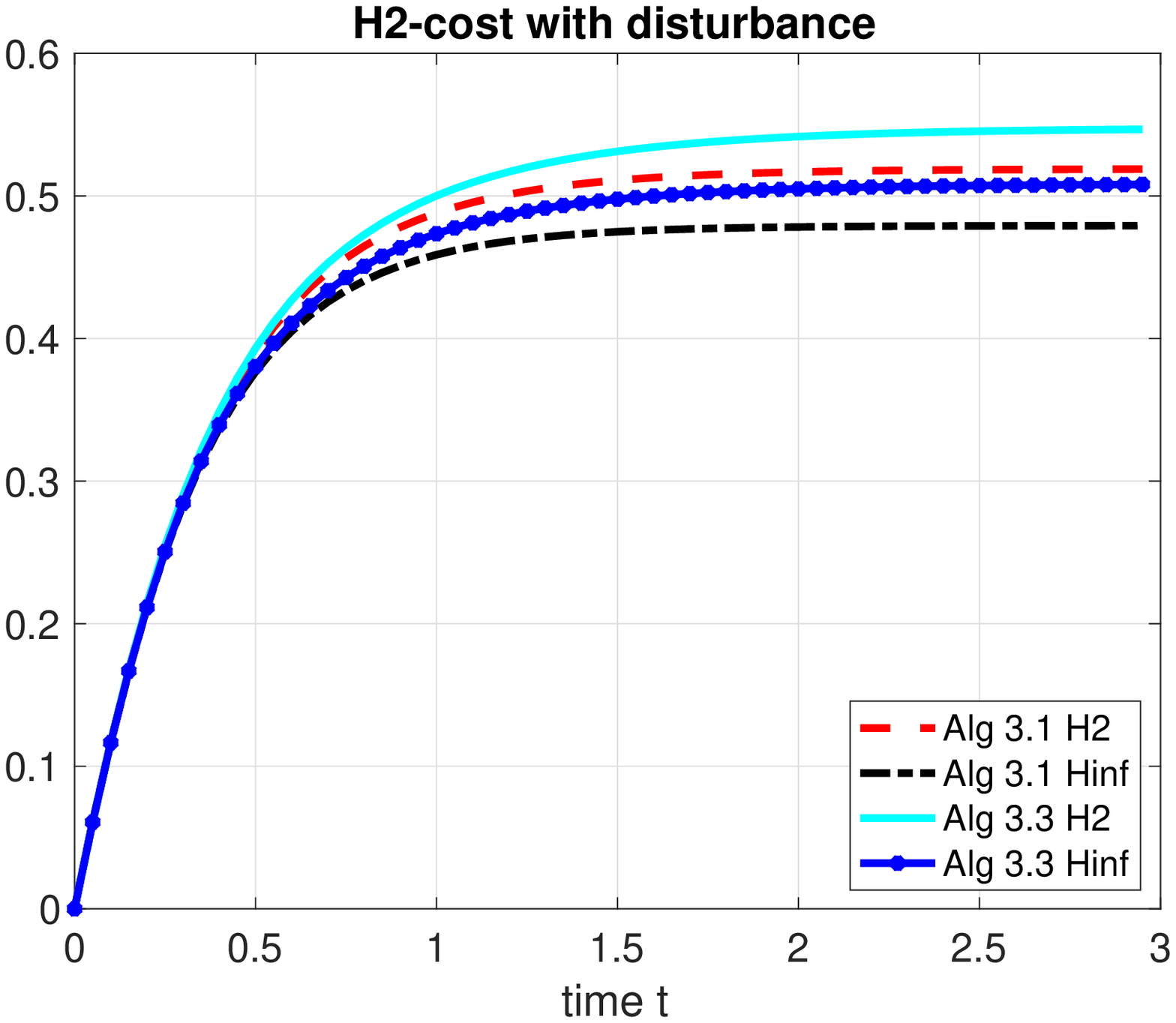}
\includegraphics[scale=0.29]{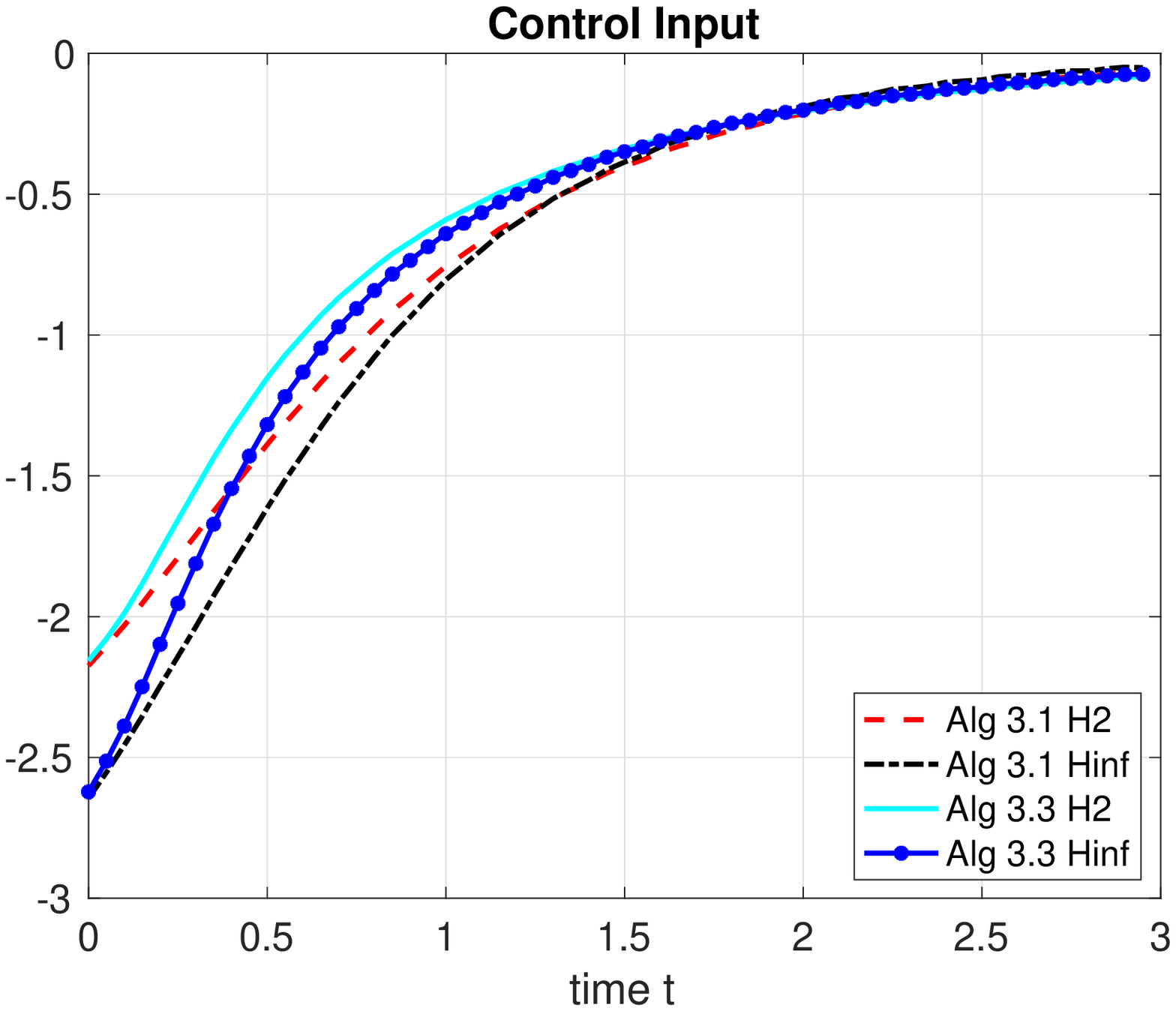}
\caption{Test \ref{sec2:ex2}: evaluation of the cost functional $\Htwo$ (left) and control input( right) for $w(t)=0.1\sin(2t)$.}
\label{fig2:an1}
\end{figure}
\end{example}

\section{Conclusions and future work}
In this work we have discussed different alternatives for the synthesis of feedback laws for stabilizing nonlinear PDEs. In particular, we have studied the use of state-dependent Riccati equation methods, both for $\Htwo$ and $\Hinf$ synthesis. Implementing an SDRE feedback law requires expressing the dynamics in semilinear form and the solution of algebraic Riccati equations at an arbitrarily high rate. This is a stringent limitation in PDE control, where high-dimensional dynamics naturally emerge from space discretization. Hence, we study offline and offline-online synthesis alternatives which circumvent or mitigate the computational effort required in the SDRE synthesis. Most notably, we have proposed an offline-online method which replaces the sequential solution of algebraic Riccati equations by Lyapunov equations. Through extensive computational experiments, including two-dimensional nonlinear PDEs, we have assessed that the SDRE offline-online method 
provides a reasonable approximation of purely online SDRE synthesis, yielding similar performance results at a reduced computational cost. Moreover, the nonlinearities arising in nonlinear reaction and nonlinear advection PDE models can be easily represented within the semilinearization framework required by SDRE methods. In conclusion, SDRE-based feedback laws constitute a reasonable alternative for suboptimal feedback synthesis for large-scale, but well structured, nonlinear dynamical systems. Future research directions include the study of the SDRE methodology for high-dimensional systems arising from interacting particle systems, and the interplay with deep learning techniques to lower the computational burden associated to a real-time implementation.

\bibliographystyle{spmpsci}
\bibliography{bibsdre}
\end{document}